\documentclass[11pt]{amsart}

\usepackage[british]{babel}
\usepackage[utf8]{inputenc}
\usepackage{latexsym}
\usepackage{color}
\usepackage{graphicx}
\usepackage{tikz}
\usetikzlibrary{arrows.meta,shapes.geometric}
\usepackage{amsfonts}
\usepackage{amssymb}
\usepackage{amsmath}
\usepackage{amsthm}
\usepackage{amscd}
\usepackage{mathtools}
\usepackage{hyperref}
\usepackage{cleveref}
\usepackage{mathrsfs}
\usepackage[margin=1.1in]{geometry}

\newtheorem{theorem}{Theorem}
\newtheorem{proposition}[theorem]{Proposition}
\newtheorem{lemma}[theorem]{Lemma}

\theoremstyle{definition}
\newtheorem{definition}[theorem]{Definition}

\theoremstyle{remark}
\newtheorem{remark}[theorem]{Remark}
\newtheorem{remarks}[theorem]{Remarks}

\numberwithin{equation}{section}
\newcommand{\Nt}{|\hskip-0.04cm|\hskip-0.04cm|}
\def\eps{{\varepsilon}}

\newcommand\N{{\mathbb N}}

\newcommand\R{{\mathbb R}}

\newcommand\Z{{\mathbb Z}}
\def\AA{{\mathcal A}}
\def\BB{{\mathcal B}}

\def\LL{{\mathcal L}}

\def\NN{{\mathcal N}}
\def\OO{{\mathcal O}}

\def\VV{{\mathcal V}}

\def\YY{{\mathcal Y}}
\def\ZZ{{\mathcal Z}}

\def\MMM{{\mathscr M}}

\def\fM{{\mathfrak M}}

\newcommand{\dd}{{\, \mathrm d}}

\newcommand{\fa}{\forall \,}

\newcommand{\dt}{\frac{{\rm d}}{{\rm d}t}}
\newcommand{\init}{\text{{\tiny {\rm in}}}}

\let\oldmarginpar\marginpar
\renewcommand\marginpar[1]{\-\oldmarginpar[\raggedleft\footnotesize #1]{\raggedright\footnotesize #1}}

\newcommand{\Blue}{\color{black}}
\newcommand{\blue}{\color{black}}

\newcommand{\sfL}{\mathsf{L}}

\newcommand{\cL}{\mathcal{L}}
\newcommand{\sfH}{\mathsf{H}}

\newcommand{\wangle}[1]{\lfloor #1 \rceil}

\newcommand{\beqn}{\begin{equation}}
\newcommand{\eeqn}{\end{equation}}
\newcommand{\bear}{\begin{eqnarray}}
\newcommand{\eear}{\end{eqnarray}}
\newcommand{\bean}{\begin{eqnarray*}}
\newcommand{\eean}{\end{eqnarray*}}

\allowdisplaybreaks

\title[Optimal decay rates for linear kinetic equations in the half-space]{Optimal decay rates for linear kinetic equations \\ in the half-space}

\author[E. Bouin]{Émeric Bouin} \address[E. Bouin]{Université Paris-Dauphine -- PSL Research University, France} \email{bouin@ceremade.dauphine.fr}

\author[S. Mischler]{Stéphane Mischler} \address[S. Mischler]{Université Paris-Dauphine -- PSL Research University \& IUF, France} \email{mischler@ceremade.dauphine.fr}

\author[C. Mouhot]{Clément Mouhot} \address[C. Mouhot]{University of Cambridge, UK} \email{c.mouhot@dpmms.cam.ac.uk}

\begin{document}

\begin{abstract}
  We prove that solutions to linear kinetic equations in a half-space with absorbing boundary conditions decay for large times like $t^{-\frac{1}{2}-\frac{d}{4}}$ in a weighted $\sfL^{2}$ space and like $t^{-1-\frac{d}{2}}$ in a weighted $\sfL^{\infty}$ space, i.e., faster than in the whole space and in agreement with the decay of solutions to the heat equation in the half-space with Dirichlet conditions. The class of linear kinetic equations considered includes the linear relaxation equation, the kinetic Fokker-Planck equation and the Kolmogorov equation with spherical velocities associated with the kinetic Brownian motion. 
\end{abstract}

\maketitle

\tableofcontents

\section{Introduction}

\subsection{The question at hand}

The optimal polynomial decay of {\blue kinetic} semigroups in unbounded domains---with or without boundary---is an interesting theoretical question, and it is important for many applications. In terms of spectral theory, it corresponds to understanding the role of the continuous spectrum in the time asymptotics.

While diffusion semigroups are well-understood, much less is known about hypocoercive semigroups of kinetic type that converge to diffusion semigroups in the parabolic time-space limit, particularly in unbounded domains with boundary. The standard elementary method to understand the polynomial decay of diffusion semigroups in the whole space is through self-similar changes of variables: it reduces the problem to establishing a spectral gap on a new unknown. However, the scaling structure of {\blue kinetic} semigroups seems too complicated for such a method to have any chance.

When it comes to the simplest unbounded domain with boundary, the half-space, the standard elementary method to understand the decay of diffusion semigroups is to extend them to the whole space by symmetrisation. Again, such methods do not seem available for the hypocoercive semigroups of kinetic type considered here, due to the more complicated interplay between spatial and kinetic variables.

\subsection{The setting}
\label{subsec:setting}

To study the interplay between hypodissipativity and unbounded domains with boundary, we consider the \emph{linear kinetic equation} of the form
\begin{equation}
  \label{eq:main}
  \partial_t f + v \cdot \nabla_x f = \cL f \quad \hbox{on}\quad (0,\infty) \times \Omega_x \times \Omega_v,
\end{equation}
where the unknown $f$ depends on the variables
\begin{equation*}
  t \in [0,+\infty), \quad x \in \Omega_x := \R_+^d \quad \text{and} \quad v \in \Omega_v,
\end{equation*}
with the velocity domain
\begin{equation*}
  \Omega_v := B(0,R) \subset \R^d \quad \text{with} \ R>0, \qquad \text{or} \qquad \Omega_v = \mathbb{S}^{d-1}, \qquad \text{or} \qquad \Omega_v = \R^d.
\end{equation*}
Here, $d \ge 1$ is the underlying space dimension and $\R_+^d$ is the open half-space
\begin{equation*}
  \R_+^d := \left\{ x \in \R^d \, : \, x_1 > 0 \right\}.
\end{equation*}

The unknown $f \ge 0$ is the density of particles and the linear operator $\mathcal L$ acts only on the kinetic variable $v$. We denote $\Omega := \Omega_x \times \Omega_v$. We complement the evolution equation with initial conditions $f(0,\cdot,\cdot)=f_{\init}$, and with the kinetic counterpart of the Dirichlet boundary conditions on the spatial density, i.e., the ``\emph{no-incoming particle conditions}'':
\begin{equation}
  \label{eq:mainBoundary}
  f = 0 \quad \hbox{on}\quad (0,\infty) \times \partial\Omega_+ ,
\end{equation}
where $\partial\Omega_+$ is the ingoing part of the boundary:
\begin{equation*}
  \partial\Omega_+ := \left(\{0\} \times \R^{d-1} \right) \times \Big(\Omega_v \cap \{v_1 > 0 \}\Big)
\end{equation*}
(the '+' sign refers to the sign of $v_1$).
%\textcolor{red}{[CM: I had not realised but we have chosen the opposite convention than $99\%$ of papers about the ``+'' which should be ``outgoing'' part of the boundary...]} 
In other words, $f(t,x,v)=0$ whenever $x_1=0$ and $v_1 > 0$, which means there are no \emph{incoming} particles from the boundary.

\smallskip
Regarding the collision operator $\LL$ in the right-hand side of~\eqref{eq:main}, we assume that it acts only on the velocity variable $v \in \Omega_v$ and {\Blue we} will consider the following three paradigmatic cases:
\begin{enumerate}
\item the \textit{linear relaxation operator}:
  \begin{equation*}
    \mathcal Lf := M \, \int_{\Omega_v} f(v_*) \dd v_* - f,
  \end{equation*}

\item the \textit{Fokker-Planck operator} (only with $\Omega_v = \R^d$):
  \begin{equation*}
    \mathcal Lf := \Delta_v f + \nabla_v \cdot ( {\Blue v f}),
  \end{equation*}

\item the \textit{Laplace-Beltrami operator} (only with  $\Omega_v = \mathbb{S}^{d-1}$):
  \begin{equation*}
    \mathcal Lf := \Delta^{\mathbb{S}^{d-1}}_v f.
  \end{equation*}
\end{enumerate}
When $\mathcal{L}$ is the linear relaxation operator, we assume that either 

\medskip

(a) $\Omega_v$ is bounded and $M$ is a uniformly bounded probability measure such that
\begin{equation}
  \label{eq:prop-M}
  %\int_{\Omega_v}  |v|^2 M(v) \dd v < \infty, \quad
  \int_{\Omega_v} v M(v)  \dd v =0 \qquad \text{and} \qquad 
  \int_{\Omega_v} (v \otimes v) \, M(v)  \dd v = \text{Id},
\end{equation}
or 

(b) $\Omega_v := \R^d$ and $M$ is the normalised gaussian
\begin{equation*}
  M(v) = \MMM_{d}(v) := (2\pi)^{-\frac{d}{2}} e^{-\frac{|v|^2}2}.
\end{equation*}
Observe that in both scenarios, $M$ is a local equilibrium profile, i.e., $\LL M = 0$. We also denote $M :=|\mathbb{S}^{d-1}|^{-1}$ the local equilibrium profile when $\Omega_v = \mathbb{S}^{d-1}$ and $\mathcal{L}$ is the Laplace-Beltrami operator. Finally in the case when $\Omega_v = \R^d$ and $\mathcal{L}$ is the Fokker-Planck operator, the local equilibrium profile $M:= \MMM_d$ is the same as for the linear relaxation. We emphasise that the moment conditions \eqref{eq:prop-M} hold true in all cases. 

\medskip

To proceed, we first introduce the functional spaces where our analysis will be carried out. We then define a collection of shorthand notation as follows. Given a weight function on the velocity variable $\omega : \Omega_v \to (0,\infty)$, $k \in \{-1,0,+1\}$, and $p \in [1,\infty]$ we define the weighted Lebesgue space
\begin{equation}
  \label{eq:defLpkomega}
  \sfL^p _{k,\omega}(\Omega) := \sfL^p\bigl(\Omega;\wangle{x_1,v}^k\omega\bigr), \qquad
  \| f \|_{\sfL^{p} _{k,\omega}(\Omega)} := \left( \int_{\Omega} |f|^p \wangle{x_1,v}^{kp} \omega^p \dd x \dd v \right)^{\frac{1}{p}}.
\end{equation}
When the domain is clear from the context, we simply write $\sfL^p_{k,\omega}$. Furthermore, we omit the index $k$ when $k = 0$, writing $\sfL^p_{\omega} := \sfL^p_{0,\omega}$, and we omit $\omega$ when $\omega \equiv 1$, writing $\sfL^p_k := \sfL^p_{k,1}$. When $k=0$ and $\omega=1$, we simply write $\sfL^p$. We have used the notation $\wangle{\cdot} := (1 + |\cdot|^2)^{\frac{1}{2}}$ and its generalised version
\begin{equation*}
  \wangle{x_1,v}:= (1+|x_1|^2+|v|^2)^{\frac{1}{2}},
\end{equation*}
and the following decomposition of the spatial variable
\begin{equation*}
    x = (x_1,x'), \quad  x_1 \ge 0, \quad x' \in \R^{d-1}.
\end{equation*}
When $\Omega_v$ is bounded, we set $\omega := 1$ and replace $\wangle{x_1,v}$ with $\wangle{x_1}$ in the polynomial pre-factor weight.

\subsection{Main results}
\label{sec:results}

To present our main results, we first introduce the class of \emph{admissible} weight functions that will be employed in our decay estimates.

\begin{definition}[Admissible weight functions]
  \label{def:admissible}
  When $\Omega_v$ is bounded, the only permissible weight function is $\omega = 1$. For the case $\Omega_v = \mathbb{R}^d$, the admissible weight functions depend on the choice of collision operator, as described below:
  \begin{itemize}
 
  \item[(i)] For the linear relaxation operator, we call \emph{admissible weight function} any function $\omega : \R^d \to \R_+$ such that
    \begin{equation*}
      %\label{eq:adm-relax}
      \forall \, v \in \Omega_v, \quad \wangle{v}^d \lesssim \omega(v) \lesssim \frac{1}{\MMM_d(v)}.
    \end{equation*}
    \smallskip

  \item[(ii)] For the Fokker-Planck operator, we call \emph{admissible weight functions} the functions
    \begin{equation*}
      %\label{eq:adm-fp}
      \forall \, v \in \Omega_v, \quad \omega(v) := e^{r \wangle{v}^s},
    \end{equation*}
    with $s=1$ and $r > 1$, or $s \in (1,2)$ and $ r > 0$, or $s=2$ and $r \in (0,\frac{1}{2})$.
  \end{itemize}
\end{definition}

\smallskip
Observe that in both cases, $M^{-\frac{1}{2}}$ is an admissible weight function. We always denote
\begin{equation}
  \label{eq:notation-weights}
  \omega_0 := 1, \qquad \omega_1 := M^{-\frac{1}{2}}, \qquad \omega_2 := M^{-1}.
\end{equation}

Our main results are as follows.

\begin{theorem}[Integral and pointwise decay estimates]
  \label{theo:main}
  Let $f_{\init} \in \sfL^2_\omega$ non-negative, and let $f_t \in \sfL^2_\omega$ be the corresponding unique global solution to~\eqref{eq:main}-\eqref{eq:mainBoundary} (see Lemma~\ref{lem:exist} for the exact notion of solution). Let $\omega$ be an admissible weight function in the sense of Definition~\ref{def:admissible}. Then we have:

  \begin{itemize}
  \item[(i)] \textbf{Integral decay}: there exists $\mathfrak a > 0$ such that if one assumes additionally $f_{\init} \in \sfL^1_1 \cap \sfL^2_\omega$, then
    \begin{equation}
      \label{eq:decay-int}
      \forall \, t \ge 0, \quad \| f_t \|_{\sfL^{2}_\omega} \lesssim \wangle{t}^{-\frac{1}{2}-\frac{d}4} \| f_{\init} \|_{\sfL^{1}_{1}} + e^{-\mathfrak a t} \|  f_{\init} \|_{\sfL^{2}_\omega}.
    \end{equation}
  \item[(ii)] \textbf{Pointwise decay}: there exists $\mathfrak a > 0$ such that if one assumes additionally $f_{\init} \in \sfL^1_1 \cap \sfL^\infty _{-1,\omega}$ then
    \begin{equation}\label{eq:decay-pw}
      \forall \, t \ge 0, \quad \left\Vert f_t \right\Vert_{\sfL_{-1,\omega}^\infty} \lesssim \wangle{t}^{-1-\frac{d}2} \| f_{\init} \|_{\sfL^{1}_{1}} + e^{-\mathfrak a t} \| f_{\init} \|_{\sfL^{\infty}_{-1,\omega}}.
    \end{equation}
  \end{itemize}
\end{theorem}
\begin{remarks}
  \begin{enumerate}
  \item Observe that the pointwise decay estimate is stronger than the integral estimate. In particular,~\eqref{eq:decay-int} follows directly from~\eqref{eq:decay-pw} together with the uniform-in-time $\sfL^1$ moment bound established in Subsections~\ref{subsec:Moment-bddvelocities}--\ref{subsec:Moment-unbddvelocities-d1}--\ref{subsec:Moment-unbddvelocities}. However, our proof is based on \textit{first} establishing~\eqref{eq:decay-int} (see Section~\ref{sec:kin-int}) and \textit{then} deducing~\eqref{eq:decay-pw} (see Section~\ref{sec:kin-pw}).
  \item An analysis of the proof reveals that the weight $\wangle{x_1, v}$ in the polynomially decaying terms on the right-hand side of these decay estimates can, in fact, be replaced by the weaker weight $\wangle{x_1}$ without loss of validity.
  \item In the case when $\Omega_x = \R^d$ without boundary condition, the decay rates of the squared $\sfL^{2}_\omega$ norm and the $\sfL^{\infty}_\omega$ norm are only $t^{-\frac{d}2}$. See~\cite[Theorem 1]{zbMATH07239825} for the $\sfL^{2}_\omega$ decay, and the next section for a discussion of the pointwise decay. These rates are optimal in general: see the discussion in~\cite{zbMATH07239825}.
  \item This result gives a unified treatment to all operators described above.
  \end{enumerate}
\end{remarks}
We define, given $f_{\init} \in \sfL^{1} _{1}(\Omega)$, the \emph{total first-component average}
\begin{equation}
  \label{eq:defMMM0}
  \NN_0 := \int_\Omega (x_1+v_1) f_{\init}(x,v) \dd x \dd v.
\end{equation}

\begin{theorem}[Localization of the mass in dimension $d=1$]
  \label{theo:main2}
  Let $d=1$, let $\LL$ be the linear relaxation or Fokker-Planck operator, let $\omega$ be an admissible weight function in the sense of Definition~\ref{def:admissible} so that $\wangle{v}^3 \omega^{-2} \in L^1(\Omega_v)$, and let the initial data $f_{\init}$ be non-negative and belong to $\sfL^{1}_{1,\omega} \cap \sfL^{\infty}_{-1,\omega}$, with $\mathcal{N}_0 > 0$ where $\mathcal{N}_0$ is defined in~\eqref{eq:defMMM0}. Let $f_t \in \sfL^2_{\omega}$ be the corresponding unique global solution to~\eqref{eq:main}-\eqref{eq:mainBoundary} (see Lemma~\ref{lem:exist} for the exact notion of solution). 
    
  Then there are $t_0 > 0$, $0<a<b<+\infty$ and $0<A<B<+\infty$ so that
  \begin{equation}
    \label{eq:loc-mass-int}
    \forall \, t \ge t_0, \quad \frac{A}{t^{\frac{1}{2}}} \le  \int_{a \sqrt{t}}^{b\sqrt{t}} \left( \int_{\Omega_v} f_t(x,v) \dd v \right)\dd x \le  \frac{B}{t^{\frac{1}{2}}}.
  \end{equation}

  Under the same assumptions, there is also $x_t \in [a \sqrt{t},b \sqrt{t}]$ and $r>0$ so that
  \begin{equation}
    \label{eq:loc-mass}
    \forall \, t \ge t_0, \quad x \in \left[x_t-1,x_t+1\right], \ v \in \Omega_v \cap [-r,r], \quad \frac{A} {t} \le f(t,x,v) \le \frac{B}{t},
  \end{equation}
  where $a,b,A,B$ depend on $\NN_0$, $r$ and on the $\sfL^{1}_{1,\omega} \cap \sfL^{\infty}_{-1,\omega}$ norm of $f_{\init}$. 
\end{theorem}

\begin{remarks}
  \begin{enumerate}
  \item Estimate~\eqref{eq:loc-mass-int} shows, in particular, that the decay estimates \eqref{eq:decay-int} and \eqref{eq:decay-pw} are optimal in dimension $d=1$; see Subsection~\ref{ss:ioptimality} for details.
  \item The rate of decay in the first integral estimate coincides with that of the heat equation in the half-space (see the next section). In fact, this first integral estimate can easily be strengthened so that the region $[a \sqrt{t},b \sqrt{t}] \times \Omega_v$ captures the majority of the mass. By contrast, our second pointwise estimate is still optimal in terms of the dependence in $t$ of the lower and upper bounds, but the length of the spatial interval $[x_t-1,x_t+1]$ is too small for it to capture a majority of the mass.
  \item It is likely that the same result holds for more general equilibrium states $M$ in the case of the linear relaxation operator. However, in this work, we do not aim for full generality and instead focus on the paradigmatic cases. For example, more general relaxation operators, with non-degenerate cross-sections, may be considered. Fokker-Planck operators, with various types of confinement forces (leading to different decays for the local equilibrium) may also be considered.
  \item It is likely that~\eqref{eq:loc-mass} may be improved and that the following estimate can be established:  there are $t_0,\delta,C> 0$, $\nu \ge 0$, $0<a<b<+\infty$, $0<A<B<+\infty$ and $x_t \in \left[a \sqrt{t},b \sqrt{t}\right]$ so that
   \begin{equation}
     \label{eq:loc-mass-improved}
     \forall \, t \ge t_0, \ x \in \left[x_t-1,x_t+1\right], \ v \in \Omega_v, \quad \frac{A M(v)^{1+\nu}}{t} \le f(t,x,v) \le \frac{B M(v)}{t} + \frac{C}{t^{1+\delta}},
   \end{equation}
   for initial data $f_{\init}$ satisfying the same assumption as in Theorem~\ref{theo:main2}.  Incidentally, a careful reading of the proof shows that the lower bound in~\eqref{eq:loc-mass-improved} is established with $\nu=0$ for the linear relaxation equation.
   \item The extension of this paper to linear models with a collision operator preserving more than one invariant will be considered in a forthcoming work. 
 \end{enumerate}
\end{remarks}

\subsection{Structure of the paper}

In Section~\ref{sec:heat}, we review known results on the heat equation as a warm-up. It allows us to introduce a key differential inequality we will then aim at obtaining in the kinetic case, and to recall the decay rates for the heat flow that are to be matched. In Section~\ref{sec:kin-int}, we prove a differential inequality~\eqref{eq:lastdiffineq} for a modified norm~\eqref{eq:mod-norm} which is equivalent to $\sfL^2_{\omega_1}$: the proof combines hypocoercivity estimates specifically developed to include and take advantage of the partial coercivity at the boundary, trace estimates, and an improved Nash inequality. This differential inequality involves the first spatial moment~\eqref{eq:def-mt}, whose estimation is discussed in Section~\ref{sec:local}. Together, these results lead to the first integral decay estimate obtained at the end of Section~\ref{sec:local}: specifically, estimate~\eqref{eq:decay-int-1} for bounded velocities or in dimension $d=1$, and the slightly weaker version~\eqref{eq:decay-int-add} otherwise. In Section~\ref{sec:kin-pw}, we strengthen these decay estimates to establish the integral decay estimate~\eqref{eq:decay-int} as well as the pointwise decay estimate~\eqref{eq:decay-pw} in all cases. Our approach refines and generalizes the semigroup factorization techniques introduced in~\cite{MR3779780, MisMou}, enabling us to employ an approximate duality argument for a component of the semigroup. The case $\Omega_v=\R^d$ with $d\ge2$ requires additional technical intermediate estimates on a marginal of the solution. We finally deduce a precise localization of the mass in Section~\ref{sec:loc-mass} by using iterated Duhamel principles and Harnack inequalities. In the follow-up paper~\cite{BMM_delay}, we use these decay results on the initial boundary value problem to study the Bramson correction for kinetic reaction-transport equations.

\begin{figure}[htbp]
  \centering
 % \resizebox{\textwidth}{!}{%
    \begin{tikzpicture}[
    font=\fontsize{11}{13}\selectfont,
      >=Latex,
      x=0.72cm,
      ingredient/.style={rectangle, draw=red, thick, align=center, text=black, text width=4.1cm, minimum height=1.25cm, inner sep=5pt},
      decay/.style={rectangle, draw=red, thick, align=center, text=black, text width=4.9cm, minimum height=1.25cm, inner sep=6pt},
      pointwise decay/.style={rectangle, draw=red, thick, align=center, text=black, text width=5.2cm, minimum height=1.25cm, inner sep=6pt},
      arrow/.style={->, thick, draw=blue}
    ]
      \node[ingredient] (hypocoercivity) at (-6.5,2.6) {Hypocoercive estimates\\with boundary terms\\(Section~3)};
      \node[ingredient] (nash) at (0,2.6) {Improved spatial\\Nash inequality\\(Section~3)};
      \node[ingredient] (moments) at (6.5,2.6) {Spatial moment\\estimates\\(Section~4)};
      \node[decay] (integral) at (0,0) {Weighted $\sfL^1 \cap \sfL^2 \to \sfL^2$ decay\\(Section~4)};
      \node[pointwise decay] (pointwise) at (0,-2.7) {Weighted $\sfL^1 \cap \sfL^\infty \to \sfL^\infty$ decay\\(Section~5)};

      \draw[arrow] (hypocoercivity) -- (integral);
      \draw[arrow] (nash) -- (integral);
      \draw[arrow] (moments) -- (integral);
      \draw[arrow] (integral) -- node[right, align=center, text=black] {Semigroup factorization} (pointwise);
    \end{tikzpicture}%
 % }
  \caption{Structure of the proof.}
  \label{fig:proof-structure}
\end{figure}

\subsection{Acknowledgements and AI tool disclosure}

Codex was utilized to help proofread typos and to assist in generating the code for the TikZ figure. Besides those tool usages, the text in this paper was fully human generated.

\section{Reminder on the heat equation}
\label{sec:heat}

\subsection{The whole space}

Consider the heat equation
\begin{equation}
  \label{eq:heat}
  \partial_t \rho = \Delta_x \rho
\end{equation}
on the unknown $\rho = \rho(t,x)$ for $t >0$ and $x \in \R^d$, associated to an  initial datum $\rho(0, \cdot)=\rho_{\init} \in \sfL^{1}(\R^d)$. We recall Nash's argument for proving the $\sfL^2$ and $\sfL^\infty$ norms decay. Using Kato's inequality
\begin{equation*}
  \text{sign}(\rho) \Delta_x \rho \le \Delta_x |\rho|,
\end{equation*}
we have
\begin{align}
  \nonumber
  & \dt \| \rho_{t} \|_{\sfL^{1}(\R^d)} = \int_{\R^d} (\text{sign } \rho) \Delta_x \rho_{t} \dd x \le \int_{\R^d} \Delta_x |\rho_{t}| \dd x =0, \\
  \label{eq:heat-energy}
  & \dt \| \rho_{t} \|_{\sfL^{2}(\R^d)}^{2} = - 2 \| \nabla_x \rho_{t} \|_{\sfL^{2}(\R^d)}^{2},
\end{align}
so that the $\sfL^{1}$ and $\sfL^{2}$ norms are non-increasing. Combined with Nash's inequality~\cite{Nash}
\begin{equation}
  \label{eq:nash}
  \| \rho_{t} \|_{\sfL^{2}(\R^d)} ^{1+\frac{2}{d}} \lesssim \| \rho_{t} \|_{\sfL^{1}(\R^d)}^{\frac{2}{d}} \| \nabla_x \rho_{t} \|_{\sfL^{2}(\R^d)}, 
\end{equation}
the previous $\sfL^{2}$ estimate implies
\begin{equation*}
  \dt \| \rho_{t} \|_{\sfL^{2}(\R^d)}^2 \lesssim - \| \rho_{t} \|_{\sfL^{1}(\R^d)} ^{-\frac{4}{d}} \| \rho_{t} \|_{\sfL^{2}(\R^d)}^{2 + \frac{4}{d}}.
\end{equation*}
Using now $\| \rho_{t} \|_{\sfL^{1}(\R^d)} \le \| \rho_{\init} \|_{\sfL^{1}(\R^d)}$ due to the $\sfL^{1}$ decay estimate, we get
\begin{equation*}
  \dt \| \rho_{t} \|_{\sfL^{2}(\R^d)}^2 \lesssim - \| \rho_{\init} \|_{\sfL^{1}(\R^d)} ^{-\frac{4}{d}} \| \rho_{t} \|_{\sfL^{2}(\R^d)}^{2 + \frac{4}{d}}.
\end{equation*}
Given $K,\alpha>0$, $u := (\alpha  K t)^{-\frac{1}{\alpha}}$ solves the ordinary differential equation $u' = - K u^{1+\alpha}$ and using a  standard comparison principle for differential equations, we deduce  the decay
\begin{equation*}
  \| \rho_{t} \|_{\sfL^{2}(\R^d)}^2 \lesssim (1+t)^{-\frac{d}2} \| \rho_{\init} \|_{\sfL^{1}(\R^d)}^2.
\end{equation*}
To obtain the pointwise decay, one simply observes that due to the $\sfL^2$ symmetry of the semigroup, the last inequality implies by duality
\begin{equation*}
  \| \rho(t,\cdot) \|_{\sfL^{\infty}(\R^d)}^2 \lesssim t^{-\frac{d}2} \| \rho_{\init} \|_{\sfL^{2}(\R^d)}^2
\end{equation*}
and by combining the last two inequalities we get
\begin{equation*}
  \| \rho(t,\cdot) \|_{\sfL^{\infty}(\R^d)}^2 \lesssim \left(\tfrac{t}{2}\right)^{-\frac{d}2} \left\| \rho\left(\tfrac{t}{2}, \cdot\right) \right\|_{\sfL^{2}(\R^d)}^2 \lesssim \left(\tfrac{t}{2}\right)^{-d} \left\| \rho_\init \right\|_{\sfL^{1}(\R^d)}^2 \lesssim t^{-d} \left\| \rho_\init \right\|_{\sfL^{1}(\R^d)}^2.
\end{equation*}
These integral and pointwise decay rates are optimal as can be easily seen from the formula giving the fundamental solution of~\eqref{eq:heat}.

\subsection{The half-space with Dirichlet conditions}

We now consider the heat equation~\eqref{eq:heat} on $\rho=\rho(t,x)$, in the domain $x \in \R_+^d$ with the initial condition $\rho(0,\cdot)=\rho_\init \in \sfL^{1}(\R_+^d)$ and with the Dirichlet boundary conditions
\begin{equation*}
  \rho(t,x)=0 \quad \mbox{on} \quad \partial \R_+^d := \left\{ x \in \R^d \, : \, x_1=0 \right\}.
\end{equation*}
We adapt to this new setting the previous Nash type arguments. By linearity and preservation of non-negativity, we may assume without loss of generality that the initial data $\rho_\init$ is non-negative. This ensures that the solution $\rho$ remains non-negative for all times. Then, the $\sfL^{1}$ norm is non-increasing since
\begin{equation*}
  \dt \| \rho_t \|_{\sfL^{1}(\R_+^d)} = - \int_{\partial \R_+^d} \partial_{x_1} \rho(t,x) \dd x' \le 0,
\end{equation*}
where $x = (0,x')$ on the boundary. The integrand $\partial_{x_1} \rho(t,x)$ on the right-hand side is non-negative, since the outward normal derivative has the appropriate sign: $\rho \ge 0$ for $x \ge 0$ and $\rho = 0$ at the boundary.

Assuming the first moment of the initial data $\rho_\init$ is finite, a double integration by parts demonstrates that the first moment is conserved and provides information about the evolution of the mass distribution over time:
\begin{equation*}
  \dt \| x_1\rho_t \|_{\sfL^{1}(\R_+^d)} = \int_{\R_+^d} x_1 \Delta_x \rho(t,x) \dd x = - \int_{\R_+^d} \partial_{x_1} \rho(t,x) \dd x =0.
\end{equation*}

The energy estimate~\eqref{eq:heat-energy} remains unchanged, but the Nash inequality~\eqref{eq:nash} can now be improved. Given $\rho$ defined on $  \R_+^d$, we define its anti-symmetric extension $\varrho$ on $  \R^d$:
\begin{equation*}
  \forall \, x_1 \le 0, \ x' \in \R^{d-1}, \quad \varrho(x_1,x') = - \rho (-x_1,x').
\end{equation*}
The function $\varrho$ is defined on the whole space $\R^d$, has finite first spatial moment and satisfies:
\begin{equation}
  \label{eq:zero-averages}
  \forall \, x' \in \R^{d-1}, \quad \int_\R \varrho(x_1,x') \dd x_1 =0.
\end{equation}
Therefore, its Fourier transform $\hat \varrho$ satisfies $\hat \varrho(0,\xi')=0$ for all $\xi' \in \R^{d-1}$ and, given $R>0$,
\begin{equation*}
  \forall \, t \ge 0, \ \xi \in \R^d \, : \, |\xi_1| \le R, \quad |\hat \varrho(t,\xi)| \le R \left\| \partial_{\xi_1} \hat \varrho_t \right\|_{\sfL^{\infty}(\R^d)} \le R \int_{\R^d} \varrho(t,x) |x_1| \dd x.
\end{equation*}
We therefore estimate the $\sfL^{2}$ norm of $\varrho$ as follows
\begin{align*}
  \left\| \varrho \right\|_{\sfL^{2}(\R^d)}^2  = \left\| \hat \varrho \right\|^2_{\sfL^{2}(\R^d)} = \int_{|\xi_1| \le R} |\hat \varrho|^2 + \int_{|\xi_1| > R} |\hat \varrho|^2 
  \lesssim R^{d+2} \left\| x_1 \varrho \right\|_{\sfL^{1}(\R^d)} ^2 + \frac{1}{R^2} \left\| \nabla_x \varrho \right\|_{\sfL^{2}(\R^d)}^2,
\end{align*}
and by optimising $R$, we deduce the following improved variant of the Nash inequality
\begin{equation*}
  \| \varrho \|_{\sfL^{2}(\R^d)} ^{1+\frac{2}{d+2}} \lesssim \| x_1 \varrho \|_{\sfL^{1}(\R^d)}^{\frac{2}{d+2}} \| \nabla_x \varrho \|_{\sfL^{2}(\R^d)}.
\end{equation*}
Because of the anti-symmetry of $\varrho$, the same holds on $\R^d_+$ rather than $\R^d$:
\begin{equation}
  \label{eq:nash-improved-half}
  \| \rho \|_{\sfL^{2}(\R^d_+)} ^{1+\frac{2}{d+2}} \lesssim \| x_1 \rho \|_{\sfL^{1}(\R^d_+)}^{\frac{2}{d+2}} \| \nabla_x \rho \|_{\sfL^{2}(\R^d_+)}.
\end{equation}
It is worth noting that the anti-symmetrization $\rho \to \varrho$ can be applied directly to the entire time-dependent solution $\rho(t, \cdot)$ of the half-space heat equation, yielding a solution $\varrho(t, \cdot)$ to the heat equation in the whole space. However, such an anti-symmetrisation of the equation itself does not appear to be feasible in the kinetic setting discussed in the following section.

We finally plug this improved Nash inequality into the energy estimate~\eqref{eq:heat-energy} to deduce
\begin{equation*}
  \dt \| \rho_t \|_{\sfL^{2}(\R_+^d)}^2 \lesssim - \| x_1 \rho_t \|_{\sfL^{1}(\R_+^d)} ^{-\frac{4}{d+2}} \| \rho_t \|_{\sfL^{2}(\R_+^d)}^{2 + \frac{4}{d+2}}.
\end{equation*}
Using the conservation
\begin{equation*}
  \| x_1 \rho_t \|_{\sfL^{1}(\R_+^d)} = \| x_1 \rho_\init \|_{\sfL^{1}(\R_+^d)},
\end{equation*}
we obtain
\begin{equation}
  \label{eq:heat-L2-improved}
  \dt \| \rho_t \|_{\sfL^{2}(\R_+^d)}^2 \lesssim - \| x_1 \rho_\init \|_{\sfL^{1}(\R_+^d)} ^{-\frac{4}{d+2}} \| \rho_t \|_{\sfL^{2}(\R_+^d)}^{2 + \frac{4}{d+2}},
\end{equation}
which implies, by standard comparison principles again, the improved decay
\begin{equation}
  \label{eq:decay-heat-L2}
  \| \rho_t \|_{\sfL^{2}(\R_+^d)}^2 \lesssim (1+t)^{-1- \frac{d}2} \| x_1 \rho_\init \|_{\sfL^{1}(\R_+^d)}^2.
\end{equation}

To obtain the pointwise decay we use the same duality argument as in the case of the whole space, which yields
\begin{equation*}
  \left\| x_1^{-1} \rho_t \right\|_{\sfL^{\infty}(\R_+^d)} \lesssim t^{-1- \frac{d}2} \| x_1 \rho_\init \|_{\sfL^{1}(\R_+^d)}.
\end{equation*}

Again, the decay rates obtained both in $\sfL^{2}$ and $\sfL^{\infty}$ are optimal. This can be deduced from the formula on the fundamental solution; the latter being obtained from the fundamental solution in $\R^d$ by the same anti-symmetrisation procedure we have used above. Note that this simple anti-symmetrisation procedure does not seem applicable to the case of kinetic equations, and we develop another method in Section~\ref{sec:kin-int}.

\subsection{Localization of the mass in dimension $1$}
\label{subsec:localiz-parabolic}

We expect the solution to asymptotically approach the eigenmode with the lowest eigenvalue that is present in the orthogonal decomposition of the initial data with respect to the eigenbasis. In the case of the half-space, the antisymmetrization procedure enables the construction of a solution in the whole space that is odd in $x_1$ and thus has no component along the first eigenmode, which is even. We therefore anticipate that the mass will localize in a manner similar to the second eigenmode,
\begin{equation*}
  (t, x) \mapsto x\, t^{-\frac{3}{2}} e^{-\frac{x^2}{t}},
\end{equation*}
which, for any $0 < a < b$, satisfies
\begin{equation*}
  \forall\, t \ge 0, \qquad \int_{at^{\frac{1}{2}}}^{bt^{\frac{1}{2}}} x\, t^{-\frac{3}{2}} e^{-\frac{x^2}{t}} \, \mathrm{d}x \approx \frac{1}{t^{\frac{1}{2}}}
\end{equation*}
and
\begin{equation*}
  \forall \, x \in \left[ a t^{\frac{1}{2}}, b t^{\frac{1}{2}} \right], \qquad \frac{\min_{z \in [a,b]} z \,e^{-z^2}}{t} \le x \, t^{-\frac32} \, e^{-\frac{x^2}{t}} \le \frac{\max_{z \in [a,b]} z \,e^{-z^2}}{t}.
\end{equation*}
Note that the initial data cannot be orthogonal to this second eigenmode if it is non-zero and built by anti-symmetrisation from $\rho$ non-negative on $x \ge 0$. 

Although such behavior could be deduced from spectral arguments, we will instead establish these bounds directly using energy estimates. This approach will also prepare us for the kinetic case, where explicit formulas for the fundamental solutions and explicit eigenbases are not available. Standard interpolation yields (remember that in this subsection the dimension $d=1$),
\begin{equation}
  \label{eq:interpol-L1L2L11}
  \forall \, t \ge 0, \quad \| \rho_t \|_{\sfL^{1}(\R_+)} \lesssim \| \rho_t \|_{\sfL^{2}(\R_+)} ^{\frac23} \| x \rho_t \|_{\sfL^{1}(\R_+)} ^{\frac13}.
\end{equation}
Therefore, by using the conservation of the first moment and the decay estimate~\eqref{eq:decay-heat-L2}, we deduce
\begin{equation}
  \label{eq:rateL1-heat}
  \| \rho_t \|_{\sfL^{1}(\R_+)} \lesssim t^{-\frac12} \, \| x \rho_\init \|_{\sfL^{1}(\R_+)}.
\end{equation}

We then compute the time evolution of the second moment (assuming it is initially finite):
\begin{equation*}
  \dt \int_{\R_+} x^2 \rho(t,x) \dd x = 2 \int_{\R_+} \rho(t,x) \dd x \lesssim t^{-\frac12}  \, \| x \rho_\init \|_{\sfL^{1}(\R_+)}
 \end{equation*}
 which implies, after integration in time,
 \begin{equation*}
  \| x^2 \rho_t \|_{\sfL^{1}(\R_+)} \lesssim t^{\frac12}  \, \| x \rho_\init \|_{\sfL^{1}(\R_+)} + \| x^2 \rho_{\init} \|_{\sfL^{1}(\R_+)}
\end{equation*}

The Cauchy-Schwarz inequality then implies
\begin{eqnarray*}
  \| x \rho_\init \|^2 _{\sfL^{1}(\R_+)}
  &=& \| x \rho_t \|^2 _{\sfL^{1}(\R_+)} \\
  &\le& \| \rho_t \|_{\sfL^{1}(\R_+)} \| x^2 \rho_t \|_{\sfL^{1}(\R_+)} \\
  &\lesssim& t^{\frac12} \| \rho_t \|_{\sfL^{1}(\R_+)}  \| x \rho_\init \|_{\sfL^{1}(\R_+)} + \| \rho_t \|_{\sfL^{1}(\R_+)}  \| x^2 \rho_{\init} \|_{\sfL^{1}(\R_+)},
\end{eqnarray*}
which implies the lower bound
\begin{equation}
  \label{eq:heat-rho-estimate-from-below}
  \| \rho_t \|_{\sfL^{1}(\R_+)} \gtrsim t^{-\frac12} \frac{\| x \rho_\init \|^2 _{\sfL^{1}(\R_+)}}{\| x \rho_\init \|_{\sfL^{1}(\R_+)} + \| x^2 \rho_{\init} \|_{\sfL^{1}(\R_+)}}.
\end{equation}
for $t \ge 1$. Together with the upper bound~\eqref{eq:rateL1-heat}, this shows
\begin{equation*}
  \forall \, t >0, \quad \| \rho_t \|_{\sfL^{1}(\R_+)} \approx \wangle{t}^{-\frac{1}{2}}.
\end{equation*}

To more precisely localize the mass, we use the $\sfL^{2}$ estimate~\eqref{eq:decay-heat-L2} to prevent concentration near $x = 0^+$, and the estimate on the first moment to rule out loss of mass at infinity:
\begin{align*}
  & \int_0 ^{a \sqrt{t}} \rho(t,x) \dd x \lesssim \sqrt{a} t^{\frac14} \| \rho_t\|_{\sfL^{2}(\R_+)} \lesssim \sqrt{a} t^{-\frac12} \| x \rho_\init \|_{\sfL^{1}(\R_+)} \\
  & \int_{b \sqrt{t}} ^{+\infty} \rho(t,x) \dd x \lesssim b^{-1} t^{-\frac12} \| x \rho_t \|_{\sfL^{1}(\R_+)} \lesssim b^{-1} t^{-\frac12} \| x \rho_\init \|_{\sfL^{1}(\R_+)}.
\end{align*}
Therefore by taking $a>0$ small enough and $b>0$ large enough, we deduce
\begin{equation*}
  \int_{a \sqrt{t}} ^{b\sqrt{t}} \rho(t,x) \dd x \gtrsim \wangle{t}^{-\frac12}.
\end{equation*}
This proves the integral lower bound. The integral upper bound readily follows from the $\sfL^{1}$ estimate. It is then possible to translate this integral estimate into pointwise bounds without using the fundamental solution: the elliptic regularisation $\sfL^{1} \to \sfL^{\infty}$ implies the pointwise upper bound, and the Harnack inequalities imply the pointwise lower bound. For kinetic equations, we will use, in a similar way, the  hypoelliptic regularisation and hypoelliptic Harnack inequalities in the case of models with diffusion in the kinetic variable, but we will use a different approach based on the transport characteristics and the associated dispersion property in the case of non-diffusive integral collision operators.

\section{Hypocoercive estimates}
\label{sec:kin-int}

In this section, we establish hypocoercive estimates that incorporate boundary contributions and combine them with the improved Nash inequality to derive the differential inequality~\eqref{eq:lastdiffineq} for the modified norm~\eqref{eq:mod-norm}, which is equivalent to $\sfL^2_{\omega_1}$ (we recall that $\omega_1 := M^{-1/2}$, see~\eqref{eq:notation-weights}). This differential inequality involves the first spatial moment~\eqref{eq:def-mt} in the direction orthogonal to the boundary, whose estimation will be addressed in the following section.

\subsection{Preliminary}
\label{subsection-setting}

We study the linear kinetic equation~\eqref{eq:main} on $(0,\infty) \times \Omega$ with the absorbing boundary condition~\eqref{eq:mainBoundary} imposed on $(0,\infty) \times \partial\Omega_+$. Here, $\Omega$, $\mathcal{L}$, and $M$ are as specified in the introduction; in particular, we recall that when $\Omega_v = \mathbb{R}^d$, the function $M$ denotes the normalized Gaussian. Without loss of generality, we may assume that the unknown function $f$ is non-negative, since the evolution is linear and preserves the sign. % In this section, we prove \EB{a slightly weaker version of the first integral estimate \eqref{eq:decay-int}, which is a key step establishing \eqref{eq:decay-int}}.
Following \cite{zbMATH07239825}, our proof employs a modification of the systematic approach introduced in~\cite{dolbeault2015hypocoercivity}, which forms part of a broader framework for establishing hypocoercivity estimates~\cite{Herau-Nier,Herau-2006,Mouhot-Neumann,Mem-villani,MR2679358,MR4581432}. The key innovation here is to exploit coercivity at the boundary to obtain (weak) decay; see also \cite{MR2679358,MR4179249,MR4581432} for similar ideas. Throughout, we use the notation~\eqref{eq:notation-weights} for the weights $\omega_0$, $\omega_1$, and $\omega_2$.
\begin{lemma}[Properties of the collision operator]
  \label{lem:obs}
  We gather in this lemma a few useful properties on the collision operators of the introduction.
  \begin{itemize}
  \item[(i)] \textbf{Calculations of observables.} Let us denote by $\LL^\sharp$ the adjoint of $\LL$ in the space $\sfL^2(\Omega_v)$ (without weight). Then there exists $c_\LL>0$ such that
    \begin{equation}\label{eq:prop-Lsharp}
      \LL^\sharp 1 = 0 \qquad \mbox{and}\qquad \LL^\sharp( v_1) = - c_\LL v_1.
    \end{equation}
    For further reference, we denote $\phi(v) := c_\LL^{-1} v_1$.

  \item[(ii)] \textbf{Spectral gap.} $\LL$ is symmetric in $\sfL^{2}_{\omega_1}(\Omega_v)$, it has a one-dimensional kernel  $\hbox{Ker} \,\LL = \hbox{Span}(M)$ and there exists $\lambda_m > 0$ such that
    \begin{equation}
      \label{eq:sg}
      \left\langle \mathcal Lf , f \right\rangle_{\sfL^2_{\omega_1}(\Omega_v)} \le - \lambda_m \| f - \rho[f] M \|_{\sfL^2_{\omega_1}(\Omega_v)}^2,
    \end{equation}
    for $f=f(v)$ in the $\sfL^2_{\omega_1}(\Omega_v)$ domain of $\mathcal{L}$, where we use here and in the sequel the following notation for the local density and momentum
    \begin{equation}
      \label{eq:rf}
      \rho[f] := \int_{\Omega_v} f(v) \dd v, \quad
      \iota[f] := \int_{\Omega_v} v f(v) \dd v.
    \end{equation}
  \end{itemize}
\end{lemma}

\begin{proof}[{\bf Proof of Lemma~\ref{lem:obs}}]
  (i) It is immediate to check the identities with $c_\LL=1$ in the case of the linear relaxation and the Fokker-Planck operators. Let us prove them with $c_\LL = d-1$ in the case of the Laplace-Beltrami operator (assuming $d \ge 2$). We introduce the coordinate system $v=(\cos \theta,[\sin \theta] \, w)$ with $w \in \mathbb S^{d-2}$, so that
    \begin{align*}
      \Delta^{\mathbb S^{d-1}}_v = (\sin \theta)^{2-d} \partial_{\theta} \left( (\sin \theta)^{d-2} \partial_\theta \right)+ \Delta^{\mathbb S^{d-2}}_w .
    \end{align*}
  Observing then that
    \begin{align*}
      (\sin \theta)^{2-d} \partial_{\theta} \left( (\sin \theta)^{d-2} \partial_\theta \left( \frac{\cos(\theta)}{d-1} \right)\right) = - \cos(\theta),
    \end{align*}
  we deduce that the function $\phi(v) = v_1/(d-1)$ satisfies $\mathcal L^\sharp \phi = \mathcal L \phi = -v_1$.
  \smallskip

  \noindent
  (ii) It is a classical result that our three operators satisfy the spectral gap property. We refer for instance to \cite{MR1803225,MR2162224} for the linear relaxation operator,
  to \cite{MR2386063,MR3155209} for the Fokker-Planck operator and to \cite{MR3155209} for the Laplace-Beltrami operator.
  In the two last cases, the spectral gap inequality is nothing but the Poincaré inequality.
\end{proof}

We conclude this preliminary subsection by clarifying the notion of solution employed throughout the paper. To this end, for any admissible weight function $\omega$ according to Definition~\ref{def:admissible} and any $k \in \mathbb{Z}$, we define the \emph{dissipation of energy} functional $\mathcal{E}_{k,\omega}$ as follows:
\begin{enumerate}
\item[(1)] For the linear relaxation operator, we set:
  \begin{equation*}
      \mathcal E_{k,\omega}(f)=0.
  \end{equation*}
  \item[(2)] For the Fokker-Planck operator, we set: 
    \begin{equation*}
      \mathcal E_{k,\omega}(f) :=  \int_{\R_+^d \times \Omega_v} \left| \nabla_v f \right|^2  \wangle{x_1,v}^{2k} \omega^2 \dd x \dd v.
    \end{equation*}
  \item[(3)] For the Laplace-Beltrami operator (in which case we take $\omega=1$), we set:
    \begin{equation*}
      \mathcal E_{k,\omega}(f) := \int_{\R_+^d \times \Omega_v} \left| \nabla_v ^{\mathbb S^{d-1}} f \right|^2  \wangle{x_1}^{2k} \dd x \dd v,
    \end{equation*}
    where $\nabla_v ^{\mathbb S^{d-1}}$ is the orthogonal projection of the gradient on the tangent plane to the sphere.
  \end{enumerate}
  Given any smooth renormalizing function $\beta : \mathbb{R} \to \mathbb{R}$, we define the \emph{$\beta$-carré-du-champ chain-rule defect} $\mathcal{F}^\beta[f]:=\LL \beta(f) - \beta'(f) \LL f$. For the Fokker-Planck operator, it takes the explicit form
  \begin{equation*}
    \mathcal{F}^\beta[f] = \beta''(f) \left|\nabla_v f\right|^2 + d \left[ \beta(f) - f \beta'(f) \right],
  \end{equation*}
  while for the Laplace-Beltrami operator it simplifies to
  \begin{equation*}
    \mathcal{F}^\beta[f] = \beta''(f) \left|\nabla_v^{\mathbb{S}^{d-1}} f\right|^2.
  \end{equation*}
  
  \begin{lemma}
    \label{lem:exist}
    For any admissible weight function $\omega$,  any $k \in \Z$, any $f_\init \in \sfL^2 _{k,\omega}$ and any $T>0$, there is a unique function $f \in \mathcal{C}([0,T];\sfL^{2}_{k,\omega})$ 
    such that $t \mapsto \mathcal E_{k,\omega}(f(t,\cdot)) \in L^1(0,T)$ and
    \begin{align*}
      \int_0^T \!\! \int_{\Omega} f \big( \partial_t \varphi + v \cdot \nabla_x  \varphi + \mathcal{L}^\sharp  \varphi \big) \dd t \dd x\dd v
      = \int_{\Omega} f_\init \varphi(0,\cdot) \dd x\dd v,
    \end{align*}
    for any test function $\varphi \in \mathcal{C}^2_c([0,T) \times \bar\Omega)$ such that $\varphi = 0$ on $\partial\Omega_-$, where $\mathcal{L}^\sharp$ denotes the adjoint of $\mathcal{L}$ for the flat measure. It is in particular a distributional solution   to~\eqref{eq:main}-\eqref{eq:mainBoundary}.

    Furthermore, $f \in \sfL^{2}((0,T) \times \partial\Omega;\omega \wangle{x_1,v}^k|v_1|)$ and $f$ is a \emph{renormalised solution}, i.e.
    \begin{multline}
      \label{eq:renormalization}
      \int_\Omega \beta(f_T)  \varphi \dd x\dd v + 
      \int_0^T \!\! \int_{\Omega} \left[ \beta(f) \big( \partial_t \varphi + v \cdot \nabla_x  \varphi + \mathcal{L}^\sharp  \varphi \big) - \mathcal F^\beta[f] \varphi \right] \dd t \dd x \dd v \\ 
      =   \int_0^T \!\! \int_{\partial\Omega_-} |v_1| \beta(f) \varphi \dd t \dd x' \dd v +  \int_\Omega \beta(f_\init) \varphi(0,\cdot) \dd x\dd v, 
    \end{multline}
    for any renormalization function $\beta \in \mathcal{C}^2(\R)$ so that $\beta(0) = 0$ and $\beta'' \in L^\infty(\R)$, and any test function
    $\varphi \in \mathcal{C}^2([0,T] \times \bar\Omega)$ so that  $\varphi \omega^{-1} \in L^\infty((0,T) \times \Omega)$. 
  \end{lemma}
  
We omit the proof of this lemma, since the strategy is standard, and refer for instance to \cite{sanchez2024KR}, \cite[Proof of Proposition 3.3]{carrapatoso2024KR} and~\cite[Section~2]{carrapatoso2024landau} for details. For the existence of variational solutions, one may use Lions' variant of the Lax-Milgram theorem~\cite[Chap~III, \textsection 1]{MR0153974}, as in~\cite{MR875086,carrapatoso2024KR}. Continuity in time and uniqueness follow directly from the renormalised formulation~\eqref{eq:renormalization}. The derivation of the renormalised formulation is classical and relies on DiPerna-Lions techniques for transport and kinetic equations in the whole space~\cite{MR1022305,MR1014927,MR972541}. The extension to equations posed in domains, together with the associated trace theory, was developed in~\cite{MR1765137,MR2721875} and can be applied, or adapted in a straightforward way, in the present setting.

\subsection{The hypocoercivity method}

The essence of hypocoercivity techniques is the use of commutators to control what is ``missed'' by the coercivity of the operator. Here
\begin{align*}
  \dt \Vert f_t \Vert_{\sfL^2_{\omega_1}}^{2}
  & = - \int_\Omega v \cdot \nabla_x \left( f_t^2 \right) M^{-1} \dd x \dd v + 2 \left\langle \mathcal Lf_t , f_t \right\rangle_{\sfL^2_{\omega_1}} \\
  & = - \int_{\partial \R_+^d} \int_{\Omega_v} |v_1| f(t,0,x',v)^2 M^{-1} \dd v \dd x' + 2 \left\langle \mathcal Lf_t , f_t \right\rangle_{\sfL^2_{\omega_1}}.
\end{align*}
The first term on the right-hand side is non-positive thanks to the boundary conditions; this provides some coercivity at the boundary that we will exploit later. The second term on the right-hand side satisfies the spectral gap inequality~\eqref{eq:sg}: however the resulting energy estimate fails to provide coercivity on the component $\rho[f]$ of the unknown, which is ``missed'' by~\eqref{eq:sg}:
\begin{align}
  \nonumber  
  \dt \Vert f_t \Vert_{\sfL^2_{\omega_1}}^{2}
  &\le - \int_{\partial \R_+^d} \int_{\Omega_v} |v_1| f(t,0,x',v)^2 M^{-1} \dd v \dd x'  \\
  \label{eq:hypoco}
  &\qquad- 2 \lambda_m  \int_{\R_+^d}   \| f_t - \rho[f_t] M \|_{\sfL^2_{\omega_1}}^2\dd x.
 \end{align}

\subsection{Modifying the norm}

We introduce a new norm on $\Omega$, equivalent to the weighted norm $\sfL^{2}_{\omega_1}$, defined as follows. The standard Sobolev spaces $\sfH^1_0$, $\sfH^2$, and $\sfH^{-1}$, all without weights, will also be employed in the proofs. Let $\mathbf{R}$ denote the operator such that $u = \mathbf{R}\eta$, where $u$ is the unique solution to the following elliptic boundary value problem on the half-space with Dirichlet boundary conditions:
\begin{equation}
  \label{eq:ellipticRu1}
  u \in \sfH^{1}_{0}(\mathbb{R}_+^d), \qquad (1-\Delta_x)u = \eta \quad \text{in} \quad \mathbb{R}_+^d.
\end{equation}
More precisely, $\mathbf R\eta \in \sfH^{1}_{0}(\R_+^d)$ satisfies
\begin{equation}
  \label{eq:ellipticRu2}
  \forall \, \varphi \in \sfH^{1}_{0}(\R_+^d), \quad \int (\mathbf R\eta) \varphi + \int \nabla (\mathbf R\eta) \cdot \nabla \varphi = \int \eta \varphi.
\end{equation}
Then, given $\epsilon>0$ (to be chosen below), we define
\begin{equation}
  \label{eq:mod-norm}
  \Nt f \Nt^2 := \| f \|_{\sfL^2_{\omega_1}} ^2 + \epsilon  \int_{\R_+^d} \left( \nabla_x \mathbf R \rho[f] \right) \cdot \iota[f] \dd x
\end{equation}
(where $\rho[f]$ and $\iota[f]$ are defined in~\eqref{eq:rf}). 
Standard ellipticity estimates (see for instance \cite[Theorem 9.25]{MR2759829} for the second estimate) yield
\begin{equation}
  \label{eq:ellip}
  \begin{cases}
    \| \mathbf R \eta \|_{\sfH^{1}_{0}(\R^d_+)} \lesssim \| \eta \|_{\sfH^{-1}(\R^d_+)} \\[2mm]
    \| \mathbf R \eta \|_{\sfH^{1}_{0}(\R^d_+)} + \| \nabla_x^2 \mathbf R \eta \|_{\sfL^{2}(\R_+^d)} \lesssim \| \eta \|_{\sfL^{2}(\R^d_+)},
  \end{cases}
\end{equation}
where $\nabla_x^2$ denotes the matrix of second-order derivatives. In other words, $\mathbf R$ is bounded from $\sfL^2$ to $\sfH^1_0 \cap \sfH^2$ and from $\sfH^{-1}$ to $\sfH^1_0$.

By the Cauchy-Schwarz inequality, using either the boundedness of $\Omega_v$ or the decay of $M$ when $\Omega_v$ is unbounded, we have
\begin{equation}
  \label{eq:borne-rho&iota}
  \begin{cases}
    \| \rho[f] \|_{\sfL^{2}(\R^d_+)} ^2 \lesssim \| f \|_{\sfL^2_{\omega_1}} ^2 \\[2mm]
    \| \iota[f] \|_{\sfL^{2}(\R^d_+)} ^2 \lesssim \| f \|_{\sfL^2_{\omega_1}}^2,
  \end{cases}
\end{equation}
so that, applying the first estimate in~\eqref{eq:ellip},
\begin{equation*}
  \left| \int_{\R_+^d} \left( \nabla_x \mathbf R \rho[f] \right) \cdot \iota[f] \dd x \right| \lesssim \| f \|_{\sfL^2_{\omega_1}} ^2.
\end{equation*}
The norm $\Nt \cdot \Nt$ is thus equivalent to the $\sfL^{2}_{\omega_1}$ norm for $\epsilon >0$ small enough:  we next choose $\epsilon > 0$ such that
\begin{equation}
  \label{eq:equivNormeps}
  \frac12  \| f \|_{\sfL^2_{\omega_1}} \le \Nt f \Nt  \le 2  \| f \|_{\sfL^2_{\omega_1}}.
\end{equation}

\subsection{Time-evolution of the modified norm}

We differentiate in time the modified norm along the flow and we use~\eqref{eq:hypoco} in order to get
\begin{align}
  \nonumber
  & \dt  \Nt f_t \Nt^2 = \dt \| f_t \|_{\sfL^2_{\omega_1}} ^2 + \epsilon \dt \int_{\R_+^d} \left( \nabla_x \mathbf R \rho[f_t] \right) \cdot \iota[f_t] \dd x \\
  \nonumber
  & \le - \int_{\partial \R_+^d} \int_{\Omega_v} |v_1| f(t,0,x',v)^2 M^{-1} \dd v \dd x' - 2 \lambda_m  \| f_t - \rho[f_t] M \|_{\sfL^2_{\omega_1}}^2 \\ \label{eq:estim-hypoco}
  & \qquad + \epsilon \int_{\R_+^d} \left( \nabla_x \mathbf R \rho[(\LL-v \cdot \nabla_x)f_t] \right) \cdot \iota[f_t] \dd x +  \epsilon \int_{\R_+^d} \left( \nabla_x \mathbf R \rho[f_t] \right) \cdot \iota[(\LL-v \cdot \nabla_x)f_t] \dd x.
\end{align}
Although somewhat tedious, these calculations adhere to standard procedures. Specifically, the computation is reformulated using weak (integrated) differentiation with respect to time. It is essential to ensure that this weak differentiation in $t$ is justified. This justification relies on the definition of a solution, which ensures that $f$ is weakly differentiable in $t$ when viewed as an $\sfL^2$-valued function.

Denoting the \emph{microscopic part} of $f$ by
\begin{equation}
  \label{eq:fperp-def}
  f^\bot_t := f_t - \rho[f_t] M,
\end{equation}
 we then calculate
\begin{align*}
  & \rho[(\mathcal L-v \cdot \nabla_x)f_t] = - \nabla_x \cdot \iota[f^\bot_t], \\
  & \iota[(\mathcal L-v \cdot \nabla_x)f_t] = \iota[\mathcal Lf^\bot_t] - \nabla_x : \int_{\Omega_v} (v \otimes v) \, f_t \dd v,
\end{align*}
where we have used the second property in~\eqref{eq:prop-M} and the fact that $M$ is a local equilibrium to deduce
\begin{equation}
  \label{eq:prop-Mbis}
  \iota[f_t] = \iota[f^\bot_t], \quad \mathcal L[f_t] = \mathcal L [f^\bot_t].
\end{equation}

As a consequence, we find
\begin{align*}
  & \int_{\R_+^d} \Big[ \nabla_x \mathbf R \Big( \rho[(\mathcal L-v \cdot \nabla_x)f_t] \Big) \Big] \cdot \iota[f_t] \dd x
    = - \int_{\R_+^d} \left[ \nabla_x \mathbf R \left( \nabla_x \cdot \iota [f^\bot_t] \right) \right] \cdot \iota[f^\bot_t] \dd x \\
  & \quad = \int_{\R_+^d} \left[ \mathbf R \left( \nabla_x \cdot \iota [f^\bot_t] \right)\right] \cdot \left( \nabla_x \cdot \iota[f^\bot_t] \right) \dd x = \left\| \mathbf R \left( \nabla_x \cdot \iota [f^\bot_t] \right) \right\|^2 _{\sfH^{1}_{0}(\R^d_+)} \\
  & \quad \lesssim \big\| \nabla_x \cdot \iota[f^\bot_t] \big\|_{\sfH^{-1}(\R_+^d)} ^2 \lesssim \big\| \iota[f^\bot_t] \big\|_{\sfL^{2}(\R_+^d)} ^2 \lesssim \big\| f^\bot _t\big\|_{\sfL^2_{\omega_1}}^2,
\end{align*}
where we have used~\eqref{eq:ellipticRu2} (with $\eta_t := \nabla_x \cdot \iota [f^\bot_t]$ and $\varphi_t :=  \mathbf R \eta_t$) in the last equality and~\eqref{eq:ellip}  in the first inequality.
In order to control the last term on the right-hand side of~\eqref{eq:estim-hypoco}, we write
\begin{align*}
  \int_{\R_+^d} \big[ \nabla_x \mathbf R \left( \rho[f_t]\right) \big] \cdot \iota[(\mathcal L-v \cdot \nabla_x)f_t] \dd x
  & = \int_{\R_+^d} \left[ \nabla_x \mathbf R \left( \rho[f_t] \right) \right] \cdot \iota[\mathcal L f_t] \dd x \\
  &\quad - \int_{\R_+^d} \left[ \nabla_x \mathbf R \left( \rho[f_t]\right) \right] \cdot \left( \nabla_x : \int_{\Omega_v} (v \otimes v) \, f_t \dd v \right) \dd x,
\end{align*}
and we deduce
\begin{align}
  \nonumber
  & \int_{\R_+^d} \big[ \nabla_x \mathbf R \left( \rho[f_t]\right) \big] \cdot \iota[(\mathcal L-v \cdot \nabla_x)f_t] \dd x
    \lesssim \big\| \nabla_x \mathbf R \left( \rho[f_t] \right) \big\|_{\sfL^{2}(\R_+^d)} \| f^\bot_t \|_{\sfL^2_{\omega_1}} \\
  \label{eq:estim-hypo-inter}
  & \qquad + \int_{\R_+^d} \left[ \nabla^2 _x \mathbf R \left( \rho[f_t] \right) \right] : \left( \int_{\Omega_v} (v \otimes v) \, f_t \dd v \right) \dd x + \int_{\partial \R_+^d} \partial_{x_1} \mathbf R \left( \rho[f_t] \right) \left( \int_{\Omega_v} v_1^2 f_t \dd v \right) \dd x',
\end{align}
where we have used  the identity $\iota [ \mathcal L f_t] = - c_\LL \iota[f^\bot_t]$ (obtained by combining~\eqref{eq:prop-Lsharp} and~\eqref{eq:prop-Mbis}), the Cauchy-Schwarz inequality and the second estimate in~\eqref{eq:borne-rho&iota} in order to estimate the first term and one integration by parts for dealing with the second term.
We then decompose $f_t$ along~\eqref{eq:fperp-def} in the second term of the right-hand side of~\eqref{eq:estim-hypo-inter}:
\begin{align*}
  & \int_{\R_+^d} \big[ \nabla_x \mathbf R \left( \rho[f_t] \right) \big] \cdot \iota[(\mathcal L-v \cdot \nabla_x)f_t] \dd x
   \lesssim \big\| \nabla_x \mathbf R \left( \rho[f_t] \right) \big\|_{\sfL^{2}(\R_+^d)} \big\| f^\bot _t \big\|_{\sfL^2_{\omega_1}} \\
  & \hspace{2cm} + \int_{\R_+^d} \big[ \Delta_x \mathbf R \left( \rho[f_t] \right) \big] \rho[f_t] \dd x + \int_{\R_+^d} \left[ \nabla^2 _x \mathbf R \left( \rho[f_t] \right) \right] : \left( \int_{\Omega_v} (v \otimes v) f^\bot_t \dd v \right) \dd x  \\
  & \hspace{2cm} + \int_{\partial \R_+^d} \big[ \partial_{x_1} \mathbf R \left( \rho[f_t] \right) \big] \left( \int_{\Omega_v} v_1^2 f_t \dd v \right) \dd x'.
\end{align*}
To estimate the boundary term, we apply the trace inequality and the Cauchy--Schwarz inequality as follows:
\begin{align*}
  & \int_{\partial \R_+^d} \left[ \partial_{x_1} \mathbf{R} \left( \rho[f_t] \right) \right]^2 \, \dd x' \lesssim \left\| \nabla_x \mathbf{R} \left( \rho[f_t] \right) \right\|_{\sfH^1(\R_+^d)}^2, \\
  & \int_{\partial \R_+^d} \left( \int_{\Omega_v} v_1^2 f_t \dd v \right)^2 \dd x' \lesssim \int_{\partial \R_+^d} \int_{\Omega_v} |v_1| f(t,0,x',v)^2 M^{-1} \dd v \dd x',
\end{align*}
where we have used the boundedness of the third moment of $M$ which holds in our context since $M$ is Gaussian when $\Omega_v=\R^d$ is unbounded. We deduce overall
\begin{align*}
  & \int_{\R_+^d} \big[ \nabla_x \mathbf R \left( \rho[f_t] \right) \big] \cdot \iota[(\mathcal L-v \cdot \nabla_x)f_t] \dd x  \\
  & \lesssim \big\| \nabla_x \mathbf R \left( \rho[f_t] \right) \big\|_{\sfH^1(\R_+^d)} \left( \| f^\bot \|_{\sfL^2_{\omega_1}}^2 +  \int_{\partial \R_+^d} \int_{\Omega_v} |v_1| f(t,0,x',v)^2 M^{-1} \dd v \dd x' \right)^{\frac12} \\
  &\qquad\qquad- \big\| \nabla_x \mathbf R \left( \rho[f_t] \right) \big\|_{\sfH^1(\R_+^d)}^2 \\
  & \lesssim \| f^\bot  \|_{\sfL^2_{\omega_1}}^2 + \int_{\partial \R_+^d} \int_{\Omega_v} |v_1| f(t,0,x',v)^2 M^{-1} \dd v \dd x' 
    - \big\| \nabla_x \mathbf R \left( \rho[f_t] \right) \big\|_{\sfH^1(\R_+^d)}^2.
\end{align*}
Combining all terms, we deduce for $\epsilon$ small enough
\begin{equation*}
  \dt  \Nt f_t \Nt^2
  \lesssim - \| f_t ^\bot \|_{\sfL^2_{\omega_1}}^2  - \int_{\partial  \Omega} |v_1| f(t,0,x',v)^2 M^{-1} \dd v \dd x'
  - \epsilon \big\| \nabla_x \mathbf R \left( \rho[f_t] \right) \big\|_{\sfH^1(\R_+^d)}^2.
\end{equation*}
Let us introduce the notation
\begin{equation}
  \label{eq:notationsZY}
  \ZZ_t := \Nt f_t \Nt ^2 \quad \mbox{and} \quad \YY_t := \| f_t^\bot \|_{\sfL^2_{\omega_1}}^2 +  \| \nabla_x \mathbf R \rho[f_t] \|_{\sfH^1(\R_+^d)}^2.
\end{equation}
Then, throwing away the boundary contribution, we obtain
\begin{align}
  \label{eq:hypo}
  \dt  \ZZ_t  \lesssim - \YY_t.
\end{align}

\subsection{Improved Nash inequality}

We now use a Nash inequality adapted to our setting in order to deduce from~\eqref{eq:hypo} a decay estimate on $\Nt f_t \Nt$, and therefore on $\| f_t \|_{\sfL^2_{\omega_1}}^2$. First, observe that since $u_t = \mathbf{R} \rho[f_t]$ (where $\rho$ is defined in~\eqref{eq:rf}), Kato's inequality yields
\begin{equation*}
  -\Delta_x |u_t| + |u_t| \leq \operatorname{sign}(u_t)\left(-\Delta_x u_t + u_t\right) = \operatorname{sign}(u_t)\, \rho_t \leq |\rho_t|.
\end{equation*}
Moreover, the boundary conditions allow us, via integration by parts, to compute
\begin{align*}
  \int_{\mathbb{R}^d_+} x_1 \Delta_x |u_t|\, \mathrm{d}x
  &= \int_{\mathbb{R}^d_+} x_1 \partial_{x_1x_1} |u_t|\, \mathrm{d}x 
     + \underbrace{\sum_{i=2}^d \int_{\mathbb{R}^d_+} x_1 \partial_{x_i x_i} |u_t|\, \mathrm{d}x}_{=0} \\
  &= -\int_{\mathbb{R}^d_+} \partial_{x_1} |u_t|\, \mathrm{d}x = 0.
\end{align*}
Here, $\Delta_x |u_t|$ is interpreted in the sense of distributions. For brevity, we omit the standard but lengthy justification of these manipulations.

Therefore, we deduce
\begin{equation}
  \label{eq:estimL11-KineticEq}
  \int_{\R^d_+} x_1 |u_t| \dd x = \int_{\R^d_+} x_1 \left( - \Delta_x |u_t| + |u_t| \right) \dd x \le \int_{\R^d_+} x_1 |\rho_t| \dd x.
\end{equation}
Meanwhile, we have
\begin{equation*}
  \int_{\R^d_+} \rho_t^2 \dd x = \int_{\R_+^d} \left( u_t - \Delta_x u_t \right)^2 = \int_{\R_+^d} \left[ u_t^2 + (\Delta_x u_t)^2 + 2 |\nabla_x u_t|^2 \right] \dd x.
\end{equation*}
Plugging the improved Nash inequality~\eqref{eq:nash-improved-half} on $u$ in the previous line, we get
\begin{equation*}
  \| \rho_t \|_{\sfL^{2}(\R^d_+)}^2
  \le 2 \| \nabla_x u_t \|_{\sfH^1(\R^d_+)}^2  + \| x_1 u_t \|_{\sfL^{1}(\R^d_+)} ^{\frac{4}{d+4}} \| \nabla_x u_t \|_{\sfL^{2}(\R^d_+)} ^{2-\frac{4}{d+4}}.
\end{equation*}
This last estimate and the Kato estimate~\eqref{eq:estimL11-KineticEq} together imply
\begin{align}
  \label{eq:Nash+Kato}
  \| \rho_t \|_{\sfL^{2}(\R^d_+)}^2
  \lesssim \| \nabla_x \mathbf R \rho_t \|_{\sfH^1(\R_+^d)}^2 + \|x_1 \rho_t \|_{\sfL^{1}(\R^d_+)} ^{\frac{4}{d+4}} \| \nabla_x \mathbf R \rho_t \|_{\sfL^{2}(\R^d_+)} ^{2-\frac{4}{d+4}},
\end{align}
for any $\rho \in \sfL^1_1 (\R^d_+) \cap \sfH^{-1}(\R^d_+)$.
Recalling the notations~\eqref{eq:notationsZY}, using the equivalence~\eqref{eq:equivNormeps}
together with the notation and the elementary identities
\begin{equation*}
 \| f_t \|_{\sfL^2_{\omega_1}}^2  =  \| f^\bot _t\|_{\sfL^2_{\omega_1}}^2 +  \| \rho[f_t] \|_{\sfL^{2}(\R^d_+)}^2
\end{equation*}
and
\begin{equation}
  \label{eq:def-mt}
  {\mathfrak M}_t := \| x_1 f_t \|_{\sfL^{1}(\Omega)} =
  \| x_1 \rho[f_t]\|_{\sfL^{1}(\Omega_x)},
\end{equation}
the improved Nash type estimate~\eqref{eq:Nash+Kato} applied to $\rho_t := \rho[f_t]$ implies
\begin{equation}
  \label{eq:Nash+lesautres}
  \ZZ_t \lesssim \YY_t  + \mathfrak M_t ^{\frac{4}{d+4}} \YY_t^{1-\frac{2}{d+4}}.
\end{equation}

\smallskip
For $m \ge 0$, we define
\begin{equation*}
  \Phi_{m} (y) := y +  m^{ \frac{4}{d+4}}  y ^{\frac{d+2}{d+4}} \ge y \quad \text{ and } \quad \Psi_m := \left( \Phi_m \right)^{-1},
\end{equation*}
which is well defined because $y \mapsto \Phi_{m} (y)$ is strictly increasing. In particular, $\Psi_m$ is strictly increasing and $\Psi_m(z) \le z$ for any $z \ge 0$.
We observe that for $z \le z_0$, we have $\Psi_m(z) \le \Psi_m(z_0) \le z_0$, next
\begin{align*}
  z &= \Psi_m (z) +  m^{\frac{4}{d+4}}  \Psi_m(z) ^{1 - {\frac{2}{d+4}}} \\
    &\le  \Psi_m (z_0)^{\frac{2}{d+4}}   \Psi_m(z)^{1 - {\frac{2}{d+4}}}   +  m^{\frac{4}{d+4}}   \Psi_m(z)^{1 - {\frac{2}{d+4}}} \\
    &= \big(\Psi_m (z_0)^{\frac{2}{d+4}} + m^{\frac{4}{d+4}}\big) \Psi_m(z)^{ \frac{d+2}{d+4}}
\end{align*}
and thus
\begin{align*}
  \big(z_0^{\frac{2}{d+4}} + m^{\frac{4}{d+4}} \big)^{-\frac{d+4}{d+2}} z^{1 +\frac{2}{d+2}} \le \Psi_m(z).
\end{align*}
Together with the elementary inequality
\begin{align*}
  \big(z_0 + m^{2} \big)^{-\frac{2}{d+2}} \lesssim 
  \big(z_0^{\frac{2}{d+4}} + m^{\frac{4}{d+4}} \big)^{-\frac{d+4}{d+2}},
\end{align*}
we have established
\begin{align*}
   (z_0     +  m^{2}  )^{-\frac{2}{d+2}}  z^{1 +\frac{2}{d+2}}  \lesssim   \Psi_m(z) .
\end{align*}
From that estimate and~\eqref{eq:Nash+lesautres}, we obtain the functional inequality
\begin{align}\label{eq:Nash-lastformulation}
   \left(\ZZ_0 + \mathfrak M_t ^{2} \right)^{-\frac{2}{d+2}}  \ZZ_t^{1 +\frac{2}{d+2}} \lesssim \YY_t, 
\end{align}
conditionally to $\ZZ_t \le \ZZ_0$.
Coming back to the differential inequality~\eqref{eq:hypo}, which in particular implies $\ZZ_t \le \ZZ_0$, 
and using~\eqref{eq:Nash-lastformulation}, we conclude with
\begin{equation}
  \label{eq:lastdiffineq}
  \dt \ZZ_t \lesssim - \left(\ZZ_0 + \mathfrak M^{2}_t \right)^{-\frac{2}{d+2}}  \ZZ_t ^{1 +\frac{2}{d+2}}.
\end{equation}
At that stage, it becomes clear that one needs an estimate on the first moment $\mathfrak M_t$, which is the object of the next section.

\section{Localization estimates and integral decay}
\label{sec:local}

In this section, we derive estimates for the term $\mathfrak{M}_t$ as defined in~\eqref{eq:def-mt}. Utilizing these estimates in conjunction with~\eqref{eq:lastdiffineq}, we establish the integral decay estimates~\eqref{eq:decay-int-1} and~\eqref{eq:decay-int-add} that capture the optimal rate of polynomial decay. The analysis is organised into three subsections, addressing first the case where $\Omega_v$ is bounded, followed by the case of unbounded $\Omega_v$ in dimension $d=1$, and finally the case of unbounded $\Omega_v$ in dimensions $d \ge 2$. Note that the mass of a  non-negative solution $f$ to the linear kinetic equation~\eqref{eq:main} is always non-increasing (for our three types of collision operators $\LL$):
\begin{eqnarray*} 
  \dt \int_\Omega f(t,x,v) \dd x \dd v
  &=&   \int_{\Omega}  \left( \mathcal{L} - v \cdot \nabla_x\right) f(t,x,v) \dd x \dd v  \\
  &=& \int_{\partial \Omega_-} v_1  f(t,0,x',v) \dd v \dd x'\le 0,
\end{eqnarray*}
where we have used that $\LL^\sharp 1 = 0$, the notation $\partial\Omega_- := \partial\Omega \setminus \partial\Omega_+$, the non-incoming particle condition~\eqref{eq:mainBoundary}, and $f \ge 0$. This proves that
\begin{equation}
  \label{eq:mass-decay}
  \forall \, t \ge 0, \quad \| f_t \|_{\sfL^{1}} \le \| f_{\init} \|_{\sfL^{1}}.
\end{equation}

In the cases where the velocity domain $\Omega_v$ is bounded or the dimension $d=1$, we establish the following integral decay estimate:
\begin{equation}
  \label{eq:decay-int-1}
  \forall\, t \geq 0, \quad \|f_t\|_{\sfL^{2}_{\omega_1}} \lesssim \wangle{t}^{-\frac{1}{2}-\frac{d}{4}} \left( \|f_{\init}\|_{\sfL^{1}_{1}} + \|f_{\init}\|_{\sfL^{2}_{\omega_1}} \right).
\end{equation}

To present the weaker decay estimate obtained in this section for unbounded velocities in dimension $d \geq 2$, we introduce the following \emph{mixed} weighted Lebesgue spaces
\begin{equation*}
  \sfL^2_{\bar\omega_1}(\Omega_1; \sfL^{1}_{\bar \omega'}(\Omega')) \qquad \text{and} \qquad \sfL^2_{\bar\omega_1}(\Omega_1; \sfL^\infty_{\bar \omega'}(\Omega')).
\end{equation*}
Here, the domain $\Omega$ is given by the product $\Omega = \Omega_1 \times \Omega'$, with $\Omega_1$ corresponding to the variables $x_1 \in \R_+$ and $v_1 \in \R$, and $\Omega'$ corresponding to the remaining variables $x' \in \R^{d-1}$ and $v' \in \R^{d-1}$. The weight function $\bar\omega_1$ is associated with the variables in $\Omega_1$, while the weight $\bar\omega'$ corresponds to the variables in $\Omega'$. The norms are given by
\begin{equation}
  \label{eq:mixed-space}
  \begin{cases}
    & \displaystyle
    \|f\|_{\sfL^2_{\bar\omega_1}(\Omega_1; \sfL^{1}_{\bar\omega'}(\Omega'))} := \left( \int_{\Omega_1} \left[ \int_{\Omega'} |f| \, \bar\omega' \, \mathrm{d}x' \mathrm{d}v' \right]^2 \bar\omega_1^2 \, \mathrm{d}x_1 \mathrm{d}v_1 \right)^{\frac{1}{2}}, \\[3mm]
    & \displaystyle
    \|f\|_{\sfL^2_{\bar\omega_1}(\Omega_1; \sfL^\infty_{\bar\omega'}(\Omega'))} := \left( \int_{\Omega_1} \left[ \sup_{x',v' \in \R^{d-1}} |f| \, \bar\omega' \right]^2 \bar\omega_1^2 \, \mathrm{d}x_1 \mathrm{d}v_1 \right)^{\frac{1}{2}}.
  \end{cases}
\end{equation}
The subscript $\bar\omega'$ is omitted when $\bar\omega' = 1$.

For unbounded velocity domains $\Omega_v$ in dimension $d \geq 2$, we establish the following integral decay estimate:
\begin{equation}
  \label{eq:decay-int-add}
  \forall\, t \geq 0, \quad \|f_t\|_{\sfL^{2}_{\omega_1}} \lesssim \wangle{t}^{-\frac{1}{2}-\frac{d}{4}} \left( \|f_{\init}\|_{\sfL^{1}_{1}} + \|f_{\init}\|_{\sfL^{2}_{\omega_1}} + \|f_{\init}\|_{\sfL^{2}_{\bar\omega_1}(\Omega_1; \sfL^{1}(\Omega'))} \right).
\end{equation}

\subsection{The case of bounded velocities} \label{subsec:Moment-bddvelocities}

We first consider the cases where the velocities are bounded: $\Omega_v = B(0,1)$ with the linear relaxation operator, or $\Omega_v = \mathbb S^{d-1}$ with the Laplace-Beltrami operator.

\begin{lemma}\label{lem:moments_bdd_vel}
  Suppose $f_\init \in \sfL_1^{1} \cap \sfL_{\omega_1}^{2}$, and let $f$ denote the non-negative solution in $\sfL_{\omega_1}^2$ to the kinetic equation~\eqref{eq:main} provided by Lemma~\ref{lem:exist}, in any of the scenarios above where $\Omega_v$ is bounded. Then, for all $t \ge 0$,
  \begin{equation}
    \label{eq:bddrhox1-casborne}
    \| x_1 f_t \|_{\sfL^{1}} \lesssim \| (1+x_1) f_\init \|_{\sfL^{1}_1}.
  \end{equation}
\end{lemma}

\begin{proof}[{\bf Proof of Lemma~\ref{lem:moments_bdd_vel}}]
  We use the  function $\phi(v) = c_\LL^{-1} v_1$ introduced after~\eqref{eq:prop-Lsharp}, which satisfies $\phi \in \sfL^{\infty}$  and $\mathcal{L}^\sharp \left(\phi\right) (v) = -v_1$. Let $c_0 := \left\| \phi \right\|_\infty$.

  We compute
  \begin{align*}
    & \dt \int_\Omega \left( x_1 + \phi(v) + c_0 \right) f(t,x,v) \dd x \dd v \\
    &= \int_{\Omega} \left( x_1 + \phi(v) + c_0 \right) \left( \mathcal{L} - v \cdot \nabla_x\right) f(t,x,v) \dd x \dd v  \\
    &= \int_{\Omega} \phi(v) \mathcal{L}f \dd v \dd x - \int_{ \Omega} \left( x_1 + \phi(v) + c_0 \right) \left(  v \cdot \nabla_x  f(t,x,v) \right)\dd x \dd v \\
    & =  \int_{\partial \Omega_x \times \Omega_v} \left(\phi(v) + c_0 \right) v_1 f(t,x,v) \dd x' \dd v ,
  \end{align*}
  by using
  \begin{equation*}
    \int_\Omega (x_1+c_0) \mathcal Lf =0, \quad \int_\Omega \phi(v) \mathcal Lf = - \int_\Omega v_1 f, \quad \text{and } - \int_\Omega x_1 v \cdot \nabla_x f = \int_\Omega v_1 f.
  \end{equation*}
  Since $\phi+c_0 \ge 0$ everywhere and $f=0$ for $v_1 \ge 0$ and $x_1=0$, we get
  \begin{equation*}
    \dt \int_\Omega \left(x_1 + \phi(v) + c_0  \right) f(t,x,v) \dd x \dd v \le 0,
  \end{equation*}
  so that
  \begin{align*}
    \int_\Omega x_1 f(t,x,v) \dd x \dd v
    & \le \int_\Omega \left(x_1 + \phi(v) + c_0  \right) f(t,x,v) \dd x \dd v \\
    & \le \int_\Omega \left(x_1 + \phi(v) + c_0  \right) f_\init(x,v) \dd x \dd v \le \int_\Omega \left( x_1 + 2c_0 \right) f_\init(x,v) \dd x \dd v.
  \end{align*}
  This concludes the proof of~\eqref{eq:bddrhox1-casborne}.
\end{proof}

By combining~\eqref{eq:lastdiffineq} with Lemma~\ref{lem:moments_bdd_vel}, and using the notation from~\eqref{eq:notationsZY}, we can now establish
\begin{equation}
  \label{eq:diffineq-bddcase}
  \dt \ZZ_t \lesssim - \left(\ZZ_0 + \| f_\init \|_{\sfL^{1}_1}^2 \right)^{-\frac{2}{d+2}} \ZZ_t^{1 + \frac{2}{d+2}},
\end{equation}
which closely resembles the differential inequality~\eqref{eq:heat-L2-improved} arising in the analysis of the heat equation. We therefore deduce~\eqref{eq:decay-int-1}.

\subsection{The case of unbounded velocities in dimension $d=1$} \label{subsec:Moment-unbddvelocities-d1}

Let $\Omega_v = \R$ and let $\mathcal L$ be the linear relaxation operator or the Fokker-Planck operator. The local equilibrium $M$ is the normalised Gaussian on $\R$, i.e., $M(v_1) =  \MMM_1(v_1)$. We start by establishing a time dependent control on the first moment. In this subsection, we write $(x,v) = (x_1,v_1)$ to stress the fact we are working in dimension $d=1$.

\begin{lemma}\label{lem:moments_unbdd_vel}
  Let $f$ be a non-negative solution in $\sfL^2_{\omega_1}$ to the kinetic equation~\eqref{eq:main} associated to  the relaxation operator or the Fokker-Planck operator with velocity domain $\Omega_v = \R$. Let $u$ be a non-negative $C^1$ function on $\R_+$ increasing fast enough so that
  \begin{equation}
    \label{eq:hyp-u}
    t \mapsto \MMM_1(u(t))^{\frac{1}{2}} \in \sfL^1(\R_+),
  \end{equation}
  then we have
  \begin{equation}\label{eq:lem:moment2}
    \forall \, t \ge 0, \quad \| x_1 f_t \|_{\sfL^{1}} \lesssim
    \| f_\init \|_{\sfL^1_1} +  \| f_\init \|_{\sfL^2_{\omega_1}} + u(t) \| f_\init \|_{\sfL^1}.
  \end{equation}
\end{lemma}

\begin{proof}[{\bf Proof of Lemma~\ref{lem:moments_unbdd_vel}}]
  We calculate by integration by parts and~\eqref{eq:prop-Lsharp}
  \begin{equation*}
    \int_\Omega x_1 (- v_1 \partial_{x_1} f )+ \int_\Omega (x_1+v_1+u)  \mathcal Lf =  \int_\Omega v_1  f + \int_\Omega f  \mathcal \LL^\sharp v_1 = 0,
  \end{equation*}
($c_\LL=1$ for the linear relaxation and Fokker-Planck operators in~\eqref{eq:prop-Lsharp}) which implies
  \begin{equation*}
    \dt \int_\Omega \left(x_1 + v_1 + u\right) f(t,x_1,v_1) \dd x_1 \dd v_1 = \int_{\R} \left( v_1^2 + u(t) v_1\right)   f(t,0,v_1) \,   \dd v_1 + u'(t) \| f(t) \|_{\sfL^{1}}.
  \end{equation*}
  Integrating in time, we deduce
  \begin{multline} \label{eq:V(s)contribution}
    \| x_1 f_t \|_{\sfL^{1}} \le - \int_\Omega v_1 f(t,x_1,v_1) \dd x_1 \dd v_1 \\ + \int_\Omega (x_1  + u(0) + |v_1|) f_\init(x_1,v_1) \dd x_1 \dd v_1
    \\
    + \int_0^t   u'(s) \| f_s \|_{\sfL^1} \dd s + 
    \int_0^t \int_{\R} \left( v_1^2 + u(s) v_1 \right) f(s,0,v_1) \,    \dd v_1 \dd s,
  \end{multline}
  where we have used $uf$ non-negative to discard a term. We estimate each term separately.

  To control the first term, we estimate its time evolution:
  \begin{align*}
    \dt \int_{\R_+ \times \R} v_1 f(t,x_1,v_1) \dd x_1 \dd v_1 = \int_{\R} v_1^2 f(t,0,v_1) \dd v_1 - \int_{\R_+ \times \R} v_1 f(t,x_1,v_1) \dd x_1 \dd v_1,
  \end{align*}
  from which we deduce
  \begin{multline*}
    \int_{\R_+ \times \R} v_1 f(t,x_1,v_1) \dd x_1 \dd v_1 = e^{-t} \int_{\R_+\times \R} v_1 f_\init(x_1,v_1) \dd x_1 \dd v_1 \\
    + \int_0^t \left( \int_{\R} v_1^2 f(s,0,v_1) \dd v_1 \right) e^{-(t-s)} \dd s.
  \end{multline*}
  The Cauchy-Schwarz inequality implies
  \begin{equation*}
    \left( \int_{\R} v_1^2 f(s,0,v_1) \dd v_1 \right)^2 \leq \left( \int_{\R} |v_1|^3 \MMM_1(v_1) \dd v_1 \right) \left( \int_{\R} \vert v_1 \vert f^2(s,0,v_1)\frac{\dd v_1}{\MMM_1(v_1)} \right) ,
  \end{equation*}
  and since~\eqref{eq:hypoco} implies that
  \begin{equation*}
    \int_0^t \int_{\Omega_v} \vert v_1 \vert f^2(s,0,v_1)\frac{\dd v_1}{\MMM_1(v_1)} \le \| f_\init \|_{\sfL^{2}_{\omega_1}}^2,
  \end{equation*}
  we deduce
  \begin{equation}
    \label{eq:mmt-bd-1}
    \left| \int_{\Omega} v_1 f(t,x_1,v_1) \dd x_1 \dd v_1 \right| \lesssim \| |v_1| f_\init \|_{\sfL^{1}} +  \| f_\init \|_{\sfL^{2}_{\omega_1}},
  \end{equation}
  which bounds the first term in~\eqref{eq:V(s)contribution}.

  To control the last boundary term in~\eqref{eq:V(s)contribution}, we observe that
  \begin{equation*}
    \int_0^\infty \int_{v_1 \leq -u(s)} |v_1|^3 \MMM_1(v_1)\, \mathrm{d}v_1\, \mathrm{d}s 
    \leq \left( \int_0^\infty \MMM_1(u(s))^{\frac{1}{2}} \, \mathrm{d}s \right) \left( \int_{v_1 \leq 0} |v_1|^3 \MMM_1(v_1)^{\frac{1}{2}} \, \mathrm{d}v_1 \right) < +\infty,
  \end{equation*}
  where the finiteness follows from the assumption~\eqref{eq:hyp-u} on $u$. Thus, we estimate:
  \begin{align}
    \nonumber
    & \int_0^t \int_{\R} \left( v_1^2 + u(s) v_1 \right) f(s,0,v_1) \dd v_1 \dd s \\
    \nonumber
    & \le \int_0^t \int_{\R} {\bf 1}_{v_1+u(s) \le 0} \left( v_1^2 + u(s) v_1 \right) f(s,0,v_1) \dd v_1 \dd s \\ 
    \nonumber
    & \le \left( \int_0^t  \int_{v_1 \le - u(s)} | v_1 |^3 \MMM_1(v_1)\dd v_1 \dd s \right)^{\frac12} \left( \int_0^t \int_{\R} |v_1| f(s,0,v_1)^2 \MMM_1^{-1} \dd v_1 \dd s  \right)^{\frac12} \\
    \label{eq:mmt-bd-2}
    & \lesssim \left( \int_0^t \int_{\R} |v_1| f(s,0,v_1)^2 \MMM_1^{-1} \dd v_1 \dd s \right)^{\frac12} \lesssim \|f_\init\|_{\sfL^{2}_{\omega_1}},
  \end{align}
  by using successively that $f_s =0$ at $x_1=0$ for any $v_1 \ge 0$, the Cauchy-Schwarz inequality, the previous inequality and~\eqref{eq:hypoco} integrated in time. 

  One deduces
  \begin{equation}
    \label{eq:lem:moment1}
    \| x_1 f_t \|_{\sfL^{1}} \lesssim
    \| (1+ x_1 + |v_1|) f_\init \|_{\sfL^{1}} +  \| f_\init \|_{\sfL^{2}_{\omega_1}} + 
    \int_0^t   u'(s) \| f_s \|_{\sfL^{1}} \dd s,
  \end{equation}
  by combining~\eqref{eq:V(s)contribution} with~\eqref{eq:mmt-bd-1}-\eqref{eq:mmt-bd-2}. The last term in~\eqref{eq:lem:moment1} can be bounded by
  \begin{equation*}
    \int_0^t   u'(s) \| f_s \|_{\sfL^1} \dd s \le u(t) \| f_\init \|_{\sfL^1},
  \end{equation*}
  thanks to~\eqref{eq:mass-decay} and \( u \ge 0 \) and \( u' \ge 0 \). We thus finally deduce~\eqref{eq:lem:moment2}.
\end{proof}

\begin{proof}[Proof of~\eqref{eq:decay-int-1}]
  We employ a bootstrap argument, proceeding in two steps.

  \smallskip\noindent
  {\# \sl Step~1.}  We define $u(t) := \wangle{t}^{\eps}$ with $\eps > 0$ small, and we observe that  it satisfies the stated requirement in Lemma~\ref{lem:moments_unbdd_vel}. Then~\eqref{eq:lem:moment2} implies the control
  \begin{equation*}
    \forall \, t \in \R^+, \quad \| x_1 f_t\|_{\sfL^1} \lesssim C_0 \wangle{t}^{\eps} \quad \mbox{ with } \quad C_0 :=  \| f_\init \|_{\sfL^{1}_{1}} +  \| f_\init\|_{\sfL^2_{\omega_1}}.
  \end{equation*}
  Combined with~\eqref{eq:lastdiffineq}, it implies by ODE comparison techniques
  \begin{equation}
    \label{eq:NotOptimalDecay}
    \ZZ_t \lesssim C_0^2 \,  \wangle{t}^{ - \frac{3}2 + 2 \eps}.
  \end{equation}

  \smallskip\noindent
  {\# \sl Step~2.} The interpolation inequality~\eqref{eq:interpol-L1L2L11} and the first estimate in~\eqref{eq:borne-rho&iota} imply
  \begin{equation*}
    \| f_t \|_{\sfL^{1}}  = \| \rho_t \|_{\sfL^{1}(\Omega_x)} \le \| \rho_t \|_{\sfL^{2}(\Omega_x)}^{\frac{2}{3}} \| x_1 \rho_t\|_{\sfL^{1}(\Omega_x)}^{\frac{1}{3}} \lesssim \|f_t\|_{\sfL^2_{\omega_1}}^{\frac{2}{3}} \| x_1 f_t \|_{\sfL^1}^{\frac{1}{3}},
  \end{equation*}
  which implies, together with the estimates of the first step, the following decay on the mass
  \begin{align*}
    \| f_t \|_{\sfL^{1}} \lesssim \big( C_0 \wangle{t}^{\eps} \big)^{\frac{1}{3}} \left( C_0^2  \wangle{t}^{-\frac{3}2 + 2 \eps} \right)^{\frac{1}{3}} \lesssim C_0  \wangle{t}^{- \frac{1}2 +  \eps},
  \end{align*}
  which decays by choosing $0 < \eps < 1/2$. Using now~\eqref{eq:lem:moment1}, we deduce
  \begin{equation}
    \label{eq:uniformL11}
    \| x_1 f_t \|_{\sfL^{1}} \lesssim C_0 + C_0\int_0^t s^{\eps-1} \wangle{s}^{- \frac{1}2 +  \eps} \dd s \lesssim C_0,
  \end{equation}
  which is bounded by choosing $\eps < 1/4$. We finally proceed as in the case of bounded velocities in Section~\ref{subsec:Moment-bddvelocities}: we combine~\eqref{eq:uniformL11} with~\eqref{eq:lastdiffineq} to deduce by ODE comparison techniques
  \begin{equation}\label{eq:weakZ2}
    \ZZ_t \lesssim C_0^2 \,  \wangle{t}^{- \frac{3}2},
  \end{equation}
  which is nothing but~\eqref{eq:decay-int-1}.
\end{proof}

\subsection{The case of unbounded velocities in  dimension $d \ge 2$} \label{subsec:Moment-unbddvelocities}

Let $\mathcal L$ be the linear relaxation operator or the Fokker-Planck operator in $\R^d$ with $d \ge 2$. We define the marginal
\begin{equation}
  \label{eq:def-marginal}
  g (t,x_1,v_1) := \int_{\R^{2d-2}} f(t,x,v) \dd x'\dd v'
\end{equation}
and we observe that $g$ satisfies
\begin{equation*}
  \partial_t g + v_1 \partial_{x_1} g = \LL_1 g \quad \hbox{on}\quad (0,\infty) \times \R_+ \times \R,
\end{equation*}
where $\LL_1$ is the relaxation operator or the Fokker-Planck operator in dimension $d=1$. Therefore we apply the one-dimensional estimate~\eqref{eq:uniformL11} to $g$:
\begin{equation*}
  \| x_1 g_t \|_{\sfL^{1}} 
  \lesssim \| (1+ x_1 + |v_1|) g_{\init} \|_{\sfL^{1}} +  \| g_{\init} \|_{\sfL^{2}_{\bar \omega_1}},
\end{equation*}
with $\bar \omega_1 := \MMM_1(v_1)^{-1/2}$ is the weight in the first variable only, which implies, coming back to $f$:
\begin{equation}
  \label{eq:uniformL11g}
  \| x_1 f_t \|_{\sfL^{1}} 
  \lesssim \| f_{\init} \|_{\sfL^{1}_1} +  \| f_{\init} \|_{\sfL^{2}_{\bar \omega_1}(\Omega_1; \sfL^{1}(\Omega')))}.
\end{equation}
We then proceed exactly as in the case of bounded velocities in Section~\ref{subsec:Moment-bddvelocities} and we deduce
\begin{equation*}
  \ZZ_t \lesssim \left( \| f_{\init} \|_{\sfL^{1}_{1}}  +  \| f_{\init} \|_{\sfL^{2}_{\omega_1}} +  \| f_{\init} \|_{\sfL^{2}_{\bar \omega_1}(\Omega_1; \sfL^{1}(\Omega'))} 
  \right)^{2} \wangle{t}^{-1-\frac{d}2}.
\end{equation*}
which is nothing but~\eqref{eq:decay-int-add}.

\section{Semigroup factorisation estimates and pointwise decay}
\label{sec:kin-pw}

In this section, we prove the stronger integral decay estimate~\eqref{eq:decay-int} in all cases, and we demonstrate that the integral decay estimate~\eqref{eq:decay-int} leads to the pointwise decay estimate~\eqref{eq:decay-pw}. We refine the semigroup extension techniques introduced in~\cite{MR3779780} and further developed in \cite{MisMou,MR3489637,MR3465438,MR4265692}. The new idea as compared to these works is the use of an \emph{approximate duality argument} based on the factorisation formula.

\subsection{The dual estimate}
\label{subsec:pw-dual}

%We first reformulate the estimate~\eqref{eq:decay-int} \EB{check} in semigroup language to deduce a dual estimate. 
We denote by $S(t)$ the semigroup associated with the kinetic equation~\eqref{eq:main}–\eqref{eq:mainBoundary} in various spaces specified below (the spaces used will be indicated in the bounds, and we omit them from the notation). Using the notation from~\eqref{eq:defLpkomega},~\eqref{eq:mixed-space}, and~\eqref{eq:notation-weights}, and letting $\bar \omega_1 := \MMM_1(v_1)^{-1/2}$ denote the weight in the first velocity variable only, we can restate the estimate~\eqref{eq:decay-int-1} from the previous section as
\begin{equation}
  \label{eq:SL1->L2}
  \| S(t) \|_{\sfL^{1}_{1} \cap \sfL^{2}_{\omega_1} \to \sfL^{2}_{\omega_1}} \lesssim \wangle{t}^{-\frac{1}{2} - \frac{d}{4}}.
\end{equation}
Similarly, the estimate~\eqref{eq:decay-int-add} for the case of unbounded velocity and dimension $d \geq 2$ can be reformulated as
\begin{equation}
  \label{eq:SX->L2-bis}
  \| S(t) \|_{\sfL^{1}_{1} \cap \sfL^{2}_{\omega_1} \cap \sfL^{2}_{\bar \omega_1}(\Omega_1; \sfL^{1}(\Omega')) \to \sfL^{2}_{\omega_1}} \lesssim \wangle{t}^{-\frac{1}{2} - \frac{d}{4}}.
\end{equation}

We denote by $S^*(t)$ the adjoint semigroup  associated to the $\sfL^2_{\omega_1}$ duality product. It is associated to the dual problem
\begin{equation}
  \label{eq:main-dual}
  \left\{
    \begin{array}{l}
      \partial_t g - v \cdot \nabla_x g = \mathcal L g \quad\hbox{in}\quad (0,\infty) \times \Omega,\\[2mm]
      g = 0  \quad\hbox{in}\quad (0,\infty) \times \partial\Omega_-,
    \end{array}
  \right.
\end{equation}
where
\begin{equation*}
  \partial\Omega_- := (\{0\} \times \R^{d-1}) \times (\Omega_v \cap \{v_1 < 0 \}).
\end{equation*}
We have utilized the antisymmetry of $v \cdot \nabla_x$ and the symmetry of $\LL$ in the space $\sfL^2_{\omega_1}$. Since equation~\eqref{eq:main-dual} shares the same structure as~\eqref{eq:main}-\eqref{eq:mainBoundary}, the analysis presented in Sections~\ref{sec:kin-int} and \ref{sec:local} remain applicable. Consequently, the $\sfL^2_{\omega_1}$ solution to~\eqref{eq:main-dual} with initial data $g_{\rm in} \in \sfL^{1}_{1,\omega_0} \cap \sfL^{2}_{\omega_1}$ exhibits the same decay properties. In terms of semigroup notation, this means that $S^*$ satisfies~\eqref{eq:SL1->L2} or~\eqref{eq:SX->L2-bis} in the case of unbounded velocities in dimension $d \ge 2$.

By $\sfL^2_{\omega_1}$ duality, we deduce
\begin{equation}
  \label{eq:SL2->Linfty}
  \| S(t) \|_{\sfL^{2}_{\omega_1} \, \to \, \sfL^{\infty}_{-1,\omega_2} \, + \, \sfL^{2}_{\omega_1}} \lesssim \wangle{t}^{-\frac12-\frac{d}4}
\end{equation}
in the case of bounded velocities or dimension $d=1$, or
\begin{equation}
  \label{eq:SL2->Xprim-bis}
  \| S(t) \|_{\sfL^{2}_{\omega_1} \, \to \, \sfL^{\infty}_{-1,\omega_2} \, + \, \sfL^{2}_{\omega_1} \, + \, \sfL^{2}_{\bar \omega_1}(\Omega_1 ;\sfL^{\infty}_{\bar \omega'}(\Omega'))}  \lesssim \wangle{t}^{-\frac12-\frac{d}4},
\end{equation}
in the case of unbounded velocities in dimension $d \ge 2$, where $\bar \omega':=\MMM_{d-1}^{-1}$ in the variables $v'\in \R^{d-1}$. We have used the following standard norm on a sum $X+Y$ of two normed vector spaces:
\begin{equation*}
  \|\zeta \|_{X+Y} := \inf \Big\{ \| x \|_X + \| y \|_Y  \ : \ \zeta=x+y, \ x \in X, \ y \in Y  \Big\},
\end{equation*}
and the following standard norm on a sum $X+Y+Z$ of three normed vector spaces:
\begin{equation*}
  \|\zeta \|_{X+Y+Z} := \inf \Big\{ \| x \|_X + \| y \|_Y + \| z \|_Z \ : \ \zeta=x+y+z, \ x \in X, \ y \in Y, \ z \in Z \Big\}.
\end{equation*}

\subsection{The extension argument}
\label{subsec:kin-pw-splitting}

The main idea is to introduce a decomposition
\begin{equation}
  \label{eq:splitting}
  \mathcal{L} - v \cdot \nabla_x = \mathcal{A} + \mathcal{B},
\end{equation}
where the semigroup $S_\BB$ generated by $\mathcal{B}$, together with the boundary condition~\eqref{eq:mainBoundary}, exhibits rapid decay as $t \to \infty$, and the time-dependent operator $\AA S_\BB(t) \AA$ possesses regularization properties. The specific form of this splitting depends on the operator $\mathcal{L}$ under consideration; precise definitions will be provided in Subsections~\ref{subsec-split-relax}, \ref{subsec-split-FP}, and~\ref{subsec-split-LB}.

The Duhamel formula can be expressed in semigroup notation as follows:
\begin{equation*}
  S =  S_{\mathcal B} + (S_{\mathcal B}\mathcal A) \star S =  S_{\mathcal B} + S \star (\mathcal A  S_{\mathcal B}).
\end{equation*}
Here, the symbol $\star$ denotes the time convolution of two functions defined on the non-negative real axis. By iteration, we deduce (see for instance \cite[(1.28)]{MR4265692}):
\begin{equation}
  \label{eq:decompo}
  \forall \, n \ge 1, \quad S = \underbrace{\sum_{j=0}^{2n-1} (S_{\mathcal B} \mathcal A)^{\star j} \star S_{\mathcal B}}_{S_1} + \underbrace{(S_{\mathcal B} \mathcal A)^{\star n} \star S \star (\mathcal A S_{\mathcal B})^{\star n}}_{S_2} =: S_1 + S_2.
\end{equation}
It is important to note that, when considered individually, the terms $S_1$ and $S_2$ in this decomposition do not satisfy the semigroup property. Additionally, observe that $S_2$ possesses a symmetric arrangement of prefactors $S_{\mathcal B} \mathcal A$ and postfactors $\mathcal A S_{\mathcal B}$. The core decay estimates are as follows:
\begin{proposition}
  \label{prop:estSB}
  Given $\omega$ an admissible weight function according to the Definition~\ref{def:admissible},  there are $\mathfrak{a}_0>0$  and $n \ge 2$ such that for all $t \ge 0$:
  \begin{align}
    \label{eq:dec-B}
    & \left\|  S_{\mathcal B}(t) \right\|_{\sfL^{2}_{\omega} \to \sfL^2_{\omega}}  + \| S_{\mathcal B}(t) \|_{\sfL^{1}_{1} \to \sfL^{1}_{1}} + \| S_{\mathcal B}(t) \|_{\sfL^{1}_{1,\omega} \to \sfL^{1}_{1,\omega}} + \| S_{\mathcal B}(t) \|_{\sfL^{\infty}_{-1,\omega} \to \sfL^{\infty}_{-1,\omega}} \lesssim  e^{-\mathfrak{a}_0t}, \\
    \label{eq:dec-A}
    & \| \AA \|_{\sfL^{1}_{1} \to \sfL^{1}_{1,\omega}} + \| \AA \|_{\sfL^{\infty}_{-1,\omega} \to \sfL^{\infty}_{-1,\omega}} \lesssim 1, \\
    \label{eq:dec-iterated-gal}
    & \left\| (\AA S_{\mathcal B}  )^{\star n} \right\|_{\sfL^{1} \to \sfL^{2}_{\omega_1}} +
      \left\| (S_{\mathcal B} \mathcal A)^{\star n} \right\|_{\sfL^{2}_{\omega} \to \sfL^2_{\omega}} + \left\| (S_{\mathcal B} \mathcal A)^{\star n}\right\|_{ \sfL^{2}_{\omega_1} \to \sfL^{\infty}_{\omega}} \lesssim e^{- \mathfrak{a}_0 t},
  \end{align}
  with also finally, in the case of unbounded velocities with dimension $d \ge 2$:
  \begin{equation}
    \label{eq:dec-iterated-gal-1DprimeBIS}
      \left\| (\AA S_{\mathcal B}  )^{\star n} \right\|_{\sfL^{1}  \to \sfL^{2}_{\bar \omega_1}(\Omega_1 ; \sfL^{1}(\Omega'))}  + \left\| (S_{\mathcal B} \mathcal A)^{\star n} \right\|_{\sfL^{2}_{\bar \omega_1}(\Omega_1;\sfL^{\infty}_{\bar\omega'}(\Omega')) \to \sfL^{\infty}_{\omega}} \lesssim e^{- \mathfrak{a}_0 t}.
  \end{equation}
\end{proposition}

The proof of this proposition is presented case by case in the following subsections. We now proceed to prove Theorem~\ref{theo:main}, assuming the validity of Proposition~\ref{prop:estSB}. The proof does not depend on the explicit structures of $\mathcal{A}$ and $\mathcal{B}$ but only on the above estimates.

\begin{proof}[{\bf Proof of Theorem~\ref{theo:main} assuming Proposition~\ref{prop:estSB}}]
  \text{ } \medskip
  
  \noindent
  {\bf (I) Case of bounded velocities or unbounded velocities in dimension $d=1$.} Let us first prove~\eqref{eq:decay-int}. Using the first estimate in~\eqref{eq:dec-B} and the second estimate in~\eqref{eq:dec-iterated-gal} gives 
  \begin{equation*}
    \forall \, t \ge 0, \quad \big\| S_1 (t) \big\|_{\sfL^{2}_\omega \to \sfL^{2}_\omega} \lesssim e^{- \mathfrak{a}_0 t}.
  \end{equation*}
  The estimates~\eqref{eq:dec-B}--\eqref{eq:dec-A} and the first estimate in~\eqref{eq:dec-iterated-gal} yield
  \begin{equation*}
    \forall \, t \ge 0, \quad \big\| (\AA S_\BB)^{\star n} (t) \big\|_{\sfL^{1}_{1} \to \sfL^{1}_{1,\omega} \cap \sfL^{2}_{\omega_1}} \lesssim e^{-\frac{\mathfrak{a}_0}{2}t}.
  \end{equation*}
  Combining this to the second estimate in~\eqref{eq:dec-iterated-gal}, and \eqref{eq:SL1->L2}, one gets
  \begin{equation*}
    \forall \, t \ge 0, \quad \big\| S_2 (t) \big\|_{\sfL^{1}_{1} \to \sfL^{2}_\omega}^2 \lesssim \wangle{t}^{-1-\frac{d}2}.
  \end{equation*}
  We have thus obtained, 
  \begin{equation*}
    \forall \, t \ge 0, \quad \| f_t \|_{\sfL^{2}_\omega}^2 \lesssim \wangle{t}^{-1-\frac{d}2} \| f_{\init} \|_{\sfL^{1}_{1}}^2 + e^{-\mathfrak a_0 t} \|  f_{\init} \|_{\sfL^{2}_\omega}^2.
  \end{equation*}
  which implies~\eqref{eq:decay-int}.
  \smallskip
  
  Let us now prove that~\eqref{eq:decay-int} implies~\eqref{eq:decay-pw}. We first write
  \begin{align}
    \nonumber
    S (t) & = S\left(\tfrac{t}{2}\right) S\left(\tfrac{t}{2}\right)  \\
    \label{eq:decomposition}
          & =  S_1\left(\tfrac{t}{2}\right) S_1\left(\tfrac{t}{2}\right) +  S_1\left(\tfrac{t}{2}\right) S_2\left(\tfrac{t}{2}\right) +  S_2\left(\tfrac{t}{2}\right) S_1\left(\tfrac{t}{2}\right) +  S_2\left(\tfrac{t}{2}\right) S_2\left(\tfrac{t}{2}\right)
  \end{align}
  thanks to~\eqref{eq:decompo}. We then estimate successively the four terms in the right-hand side of~\eqref{eq:decomposition}.
  
  \smallskip\noindent
  {\# \sl Step~1 (first term in~\eqref{eq:decomposition}).} Combining the third estimate in~\eqref{eq:dec-B} and the second estimate in~\eqref{eq:dec-A}, we get
  \begin{equation}
    \label{eq:S1estim1}
    \forall \, t \ge 0, \quad \big\| S_1(t) \big\|_{\sfL^{1}_{1,\omega} \to \sfL^{1}_{1,\omega}} +
    \big\| S_1(t) \big\|_{\sfL^{\infty}_{-1,\omega}\to\sfL^{\infty}_{-1,\omega}} \lesssim e^{-\frac{\mathfrak{a}_0}{2}t}.
  \end{equation}
  
  \smallskip\noindent
  {\# \sl Step 2 (second term in~\eqref{eq:decomposition}).} The estimates~\eqref{eq:dec-B}--\eqref{eq:dec-A} and the first estimate in~\eqref{eq:dec-iterated-gal} yield
  \begin{equation}
    \label{eq:S2estimASBL1->L1capL2}
    \forall \, t \ge 0, \quad \big\| (\AA S_\BB)^{\star n} (t) \big\|_{\sfL^{1}_{1} \to \sfL^{1}_{1,\omega} \cap \sfL^{2}_{\omega_1}} \lesssim e^{-\frac{\mathfrak{a}_0}{2}t}.
  \end{equation}
  Combining this estimate~\eqref{eq:S2estimASBL1->L1capL2} with~\eqref{eq:SL1->L2} and the second  estimate in~\eqref{eq:dec-iterated-gal}, we get
  \begin{equation*}
    \forall \, t \ge 0, \quad \big\| S_2 (t) \big\|_{\sfL^{1}_{1} \to  \sfL^{\infty}_{\omega}} \lesssim \wangle{t}^{-\frac12-\frac{d}4}.
  \end{equation*}
  Since $\sfL^{\infty}_{-1,\omega} \lesssim \sfL^{\infty}_{\omega}$, we deduce
  \begin{equation}
    \label{eq:S2estim1}
    \forall \, t \ge 0, \quad \big\| S_2 (t) \big\|_{\sfL^{1}_{1} \to  \sfL^{\infty}_{-1,\omega}} \lesssim \wangle{t}^{-\frac12-\frac{d}4}.
  \end{equation}
  This last estimate~\eqref{eq:S2estim1} combined with the second estimate in~\eqref{eq:S1estim1} yields
  \begin{equation*}
    \forall \, t \ge 0, \quad \big\| S_1\left(\tfrac{t}{2}\right) S_2\left(\tfrac{t}{2}\right)\big\|_{\sfL^{1}_{1} \to  \sfL^{\infty}_{-1,\omega}} \lesssim e^{-\frac{\mathfrak{a}_0}{2}t}.
  \end{equation*}
  
  \smallskip\noindent
  {\# \sl Step 3 (third term in~\eqref{eq:decomposition}).}
  A similar argument based on~\eqref{eq:S2estim1} and the first estimate in~\eqref{eq:S1estim1} gives
  \begin{equation*}
    \forall \, t \ge 0, \quad \big\| S_2\left(\tfrac{t}{2}\right) S_1\left(\tfrac{t}{2}\right) \big\|_{\sfL^{1}_{1}  \to \sfL^{\infty}_{-1,\omega}} \lesssim e^{-\frac{\mathfrak{a}_0}{2}t}.
  \end{equation*}
  
  \smallskip\noindent
  {\# \sl Step 4 (fourth term in~\eqref{eq:decomposition}).}
  Here, we use the estimate~\eqref{eq:SL2->Linfty} and we take advantage of the ``duality trick'' that was crucial in Nash's computation in~\cite{Nash}.
  Observing that $\omega_1$ is an admissible weight function, we get from~\eqref{eq:dec-iterated-gal} (for any $n \ge 1$)
  \begin{equation}
    \label{eq:ASBn-L2L2}
    \forall \, t \ge 0, \quad \big\| (S_\BB\AA )^{\star n} (t) \big\|_{\sfL^{2}_{\omega_1}  \to \sfL^{2}_{\omega_1}} + \big\| (\AA S_\BB)^{\star n} (t) \big\|_{\sfL^{2}_{\omega_1}  \to \sfL^{2}_{\omega_1}} \lesssim e^{- \frac{\mathfrak{a}_0}{2} t}.
  \end{equation}
  First, we deduce by combining the first estimate in~\eqref{eq:ASBn-L2L2} with~\eqref{eq:S2estimASBL1->L1capL2} and~\eqref{eq:SL1->L2}
  \begin{equation}
    \label{eq:S2estim2}
    \forall \, t \ge 0, \quad \big\| S_2 (t) \big\|_{\sfL^{1}_{1} \to \sfL^{2}_{\omega_1}} \lesssim \wangle{t}^{-\frac12-\frac{d}4}.
  \end{equation}
  Second, combining $\omega \lesssim \omega_2$ with the second estimates in respectively~\eqref{eq:dec-B} and~\eqref{eq:dec-A}, we have
  \begin{equation*}
    \forall \, t \ge 0, \quad \big\| (S_\BB \AA)^{\star n} (t) \big\|_{\sfL^{\infty}_{-1,\omega_2} \to \sfL^{\infty}_{-1,\omega}} \lesssim e^{- \frac{\mathfrak{a}_0}{2} t}.
  \end{equation*}
  Third, we combine it with the second estimate in~\eqref{eq:dec-iterated-gal} to obtain
  \begin{equation}
    \label{eq:SBAn-Linfty+L2Linfty}
    \forall \, t \ge 0, \quad \big\| (S_\BB\AA)^{\star n} (t) \big\|_{\sfL^{\infty}_{-1,\omega_2}  + \sfL^{2}_{\omega_1} \to  \sfL^{\infty}_{-1,\omega}} \lesssim e^{-\frac{\mathfrak{a}_0}{2}  t}.
  \end{equation}
  The preceding two estimates for $(\AA S_\BB)^{\star n} (t)$ and $(S_\BB \AA)^{\star n} (t)$, combined with~\eqref{eq:SL2->Linfty}, yield:
  \begin{equation}
    \label{eq:S2estim3}
    \forall \, t \ge 0, \quad \big\| S_2 (t) \big\|_{\sfL^{2}_{\omega_1} \to \sfL^{\infty}_{-1,\omega}} \lesssim \wangle{t}^{-\frac12-\frac{d}4}.
  \end{equation}
  Combining~\eqref{eq:S2estim2} and~\eqref{eq:S2estim3}, we can estimate the fourth and last term in~\eqref{eq:decomposition}:
  \begin{equation*}
    \forall \, t \ge 0, \quad \big\| S_2 \left(\tfrac{t}{2}\right) S_2 \left(\tfrac{t}{2}\right) \big\|_{\sfL^{1}_{1}  \to \sfL^{\infty}_{-1,\omega}} \lesssim \wangle{t}^{-1-\frac{d}2}.
  \end{equation*}
  
  \medskip
  Altogether, we have thus proved~\eqref{eq:decay-pw}.
  \medskip
  
  \noindent
  {\bf (II) Case of unbounded velocities in dimension $d \ge 2$.} It is enough to prove~\eqref{eq:decay-int} using the slightly weaker bound~\eqref{eq:decay-int-add}. Once this is established, the preceding argument shows that~\eqref{eq:decay-int} implies~\eqref{eq:decay-pw}, completing the proof. To prove~\eqref{eq:decay-int}, we use again the decomposition~\eqref{eq:decomposition}. Then~\eqref{eq:S1estim1} implies by interpolation
  \begin{equation*}
    \forall \, t \ge 0, \quad \big\| S_1(t) \big\|_{\sfL^2_{\omega} \to \sfL^2_{\omega}} \lesssim e^{-\frac{\mathfrak{a}_0}{2}t}.
  \end{equation*}
  Finally the combination of~\eqref{eq:dec-iterated-gal},~\eqref{eq:dec-iterated-gal-1DprimeBIS} and~\eqref{eq:decay-int-add} easily implies that the term $S_2$ satisfies 
  \begin{equation*}
    \forall \, t \ge 0, \quad \big\| S_2(t) \big\|_{\sfL^1 _1 \to \sfL^2_{\omega}} \lesssim \wangle{t}^{-\frac12-\frac{d}4}
  \end{equation*}
  which proves~\eqref{eq:decay-int} and concludes the proof.
\end{proof}

\subsection{Proof of Proposition~\ref{prop:estSB} for the linear relaxation equation}
\label{subsec-split-relax}

In this subsection, we denote by $\mathcal{L}$ the linear relaxation operator as introduced in Subsection~\ref{subsec:setting}. We assume that either $\Omega_v$ is bounded, or that $\Omega_v = \mathbb{R}^d$ with $d \ge 1$. The decomposition~\eqref{eq:splitting} is defined by
\begin{equation}
  \label{eq:dec-relaxSB}
  {\mathcal A} f  :=  \rho[f] M,
\end{equation}
and $S_\BB$ is then the semigroup associated to the evolution equation
\begin{equation}
  \label{eq:dual-B}
  \left\{
    \begin{array}{l}
      \partial_t f = -  v \cdot \nabla_x  f - f   \quad\hbox{in}\quad (0,\infty) \times \Omega,\\[2mm]
      f = 0  \quad\hbox{in}\quad (0,\infty) \times \partial\Omega_+,  \qquad f(0,\cdot,\cdot) = f_\init \quad \hbox{in}\quad \Omega.
    \end{array}
  \right.
\end{equation}
We establish the estimates in Proposition~\ref{prop:estSB} through a sequence of lemmas.
\begin{lemma}
  \label{lem:SBrelaxation-estim1}
  For any $p \in [1,+\infty]$, $k \in \{ -1, 0, 1 \}$ and any admissible weight function $\omega : \R^d \to \R_+$,
  then any solution $f$ to the equation~\eqref{eq:dual-B} in $\sfL^{p}_{k,\omega}$ satisfies the decay estimate
  \begin{equation}
    \label{eq:dec-B-relax}
    \fa t \ge 0, \quad \| f_t \|_{\sfL^{p}_{k,\omega}} = \| S_\BB(t) f_\init \|_{\sfL^{p}_{k,\omega}} \lesssim e^{- \kappa_{p,k} t} \| f_\init \|_{\sfL^{p}_{k,\omega}}
  \end{equation}
  with $\kappa_{p,-1} := 1-\frac{1}{p}$, $\kappa_{p,0} := 1$, $\kappa_{p,1} := \frac{1}{p}$. This proves~\eqref{eq:dec-B} for the linear relaxation equation.
\end{lemma}

\begin{proof}[{\bf Proof of Lemma~\ref{lem:SBrelaxation-estim1}}]
  \emph{Case $k=0$.} Let $f \ge 0$ be a solution to~\eqref{eq:dual-B} and $p \in (1,\infty)$. We compute
  \begin{align*}
    \frac1p  \dt \int_\Omega  f_t ^p   \omega^p \dd x \dd v
    &= \int_\Omega     f_t ^{p-1} \left( - v \cdot \nabla_x f_t  - f_t \right) \omega^p \dd x \dd v \\
    &= -\frac1p   \int_\Omega \nabla_x \cdot \left( v   f_t ^{p} \right) \omega^p \dd x \dd v -  \int_\Omega  f_t ^p \omega^p \dd x \dd v \\
    &= \frac1p \int_{\partial \Omega_x \times \Omega_v} v_1 f_t ^{p} \omega^p \dd x ' \dd v - \int_\Omega f_t ^p \omega^p \dd x \dd v \le - \int_\Omega f_t ^p \omega^p \dd x \dd v,
  \end{align*}
  by integrating by parts and using the boundary condition, which implies~\eqref{eq:dec-B-relax} for $k=0$.

  \smallskip \noindent
  \emph{Case $k=1$.} We first estimate
  \begin{equation*}
    \frac{v \cdot \nabla_x \wangle{x_1}^{p}}{p}
    = \left( v_1 \cdot x_1 \right) \wangle{x_1}^{p-2} \le
    \frac1p \wangle{v_1}^p + \left( 1 - \frac{1}{p} \right) \wangle{x_1}^p
  \end{equation*}
  and we deduce
  \begin{align*}
    & \frac1p\dt \int_\Omega \wangle{x_1}^p f^p \omega^p \dd x \dd v
      = - \int_\Omega \wangle{x_1}^{p} \frac{v \cdot \nabla_x \left( f^{p} \right)}{p} \omega^p \dd x \dd v -  \int_\Omega \wangle{x_1}^{p}  f^p \omega^p \dd x \dd v \\
    & \qquad = \frac{1}{p} \int_{\partial \Omega_x \times \Omega_v} v_1 \wangle{x_1}^p f^p \omega^p \dd x' \dd v + \int_\Omega \frac{v \cdot \nabla_x\left( \wangle{x_1}^p \right)}{p} f^{p} \omega^p \dd x \dd v  - \int_\Omega \wangle{x_1}^p  f^p \omega^p \dd x \dd v \\
    & \qquad \le \frac1p \int_\Omega \wangle{v_1}^p f^p \omega^p \dd x \dd v - \frac1p \int_\Omega \wangle{x_1}^p f^p \omega^p \dd x \dd v,
  \end{align*}
  where we have used the boundary conditions to ensure the boundary term has the correct sign. By the Gronwall lemma and the previous estimate in the case $k=0$ with the weight $\wangle{v} \omega$, we deduce
  \begin{equation*}
    \int_\Omega \wangle{x_1}^p f_t^p \omega^p \dd x \dd v  \le
    e^{-t} \int_\Omega \wangle{x_1}^p f_\init^p \omega^p \dd x \dd v + e^{-t} \left(\frac{1-e^{(1-p)t}}{p-1} \right) \int_\Omega \wangle{v}^p f_\init^p \omega^p \dd x \dd v.
  \end{equation*}
  Together with the estimate in the case $k=0$ with the weight $\wangle{v} \omega$, we deduce~\eqref{eq:dec-B-relax} for $k=1$.

  \smallskip \noindent
  \emph{Case $k=-1$.} Since the previous estimate holds for any weight function $\omega$, we obtain the case $k=-1$ by duality from the case $k=1$. 
 \end{proof}

We point out that the estimate~\eqref{eq:dec-B-relax} enables us to establish the existence of the associated semigroup $S_\BB$ in the spaces $\sfL^{2}_{k,\omega}$ for any admissible weight function $\omega$ and for $k \in \{-1, 0, 1\}$, by following the argument presented in the proof of Lemma~\ref{lem:exist}.

\begin{lemma}
  \label{lem-relax-AAf}
  Consider either $p = 1$ and two weight functions $\omega$ and $\tilde{\omega}$ so that $\wangle{\cdot} \omega ^{-1} \in \sfL^{\infty}(\Omega_v)$ and $\tilde \omega \wangle{\cdot} M \in \sfL^{2}(\Omega_v)$, or $p=2$ and two weight functions $\omega$ and $\tilde{\omega}$ so that   $\wangle{\cdot} \omega ^{-1} \in \sfL^{2}(\Omega_v)$ and $\tilde \omega  M \in \sfL^{2}(\Omega_v)$. Then
  \begin{equation*}
    \text{$\AA$ is bounded from $\sfL^{p}_{1,\omega}$ to $\sfL^{p}_{1,\tilde \omega}$.}
  \end{equation*}

  Consider either $p= 2$ and two weight functions $\omega$ and $\tilde{\omega}$ so that $\wangle{\cdot} \,\omega^{-1} \in \sfL^{2}(\Omega_v)$ and $\tilde \omega M \in \sfL^{2}(\Omega_v)$, or $p=\infty$ and two weight functions $\omega$ and $\tilde{\omega}$ so that $\omega^{-1} \in \sfL^{1}(\Omega_v)$ and $\tilde \omega \wangle{\cdot} M \in \sfL^{\infty}(\Omega_v)$. Then
  \begin{equation*}
    \text{$\AA$ is bounded from $\sfL^{p}_{-1,\omega}$ to $\sfL^{p}_{-1,\tilde \omega}$.}
  \end{equation*}
\end{lemma}

\begin{proof}[{\bf \bf Proof of Lemma~\ref{lem-relax-AAf}}]
  Given $p \in \lbrace 1 , 2 \rbrace$, we compute
  \begin{equation*}
    \| \AA f \|_{\sfL^{p}_{1,\tilde \omega}}
    = \big\| \rho[f] M \wangle{x_1,v} \tilde \omega \big\|_{\sfL^{p}} \le 2 \big\| M \wangle{v}\tilde \omega \big\|_{\sfL^{p}_v}  \big\| \wangle{x_1} \rho[f] \big\|_{\sfL^{p}_x}
    \le \|   f \|_{\sfL^{p}_{1,\omega }},
  \end{equation*}
  so that the first claim is proved.

  Given $p \in \left\lbrace 2, \infty \right\rbrace$, we compute
  \begin{align*}
    \| \AA f \|_{\sfL^{p}_{-1,\tilde \omega}}
    &= \left\Vert \frac{M \tilde \omega}{\wangle{x_1,v}} \rho[f]  \right \Vert_{\sfL^{p}} \leq \left\Vert \frac{M \tilde \omega}{\wangle{x_1,v}} \times \big\Vert \omega^{-1} \wangle{x_1,v}\big\Vert_{\sfL_{v}^{p'}} \times \big\Vert f \omega \wangle{x_1,v}^{-1}\big\Vert_{\sfL^p_v} \right\Vert_{\sfL^{p}}\\
    &\le \|  M \tilde \omega   \|_{\sfL^{p}_v} \| \wangle{\cdot} \omega^{-1} \|_{\sfL^{p'}} \| f \|_{\sfL^{p}_{-1,\omega}} \lesssim \| f \|_{\sfL^{p}_{-1,\omega}}
  \end{align*}
  so that the second claim is proved.
\end{proof}

We now follow \cite[Lemma 3.2]{MR3591133}, adapting a classical dispersion result of Bardos and Degond \cite{MR794002}, to establish a gain of integrability property (see also \cite{MR259662}, where such a trick is used). Then, based on the following Lemma and Lemmas~\ref{lem:SBrelaxation-estim1} and~\ref{lem-relax-AAf}, the two estimates in~\eqref{eq:dec-iterated-gal} follow classically, see \cite{MR3779780}, \cite{MisMou}, or \cite[Proposition 2.5]{MR3465438}.

\begin{lemma}
  \label{lem:dispersionAS_BA}
  There is $n_0 \in \N^*$ so that for $t > 0$
  \begin{eqnarray}
    \label{eq:AS_Btwice}
    && \big\| (\AA S_\BB)^{\star n_0}(t) \big\|_{\sfL^{1} \to \sfL^{2}_{\omega_1}} \lesssim e^{-t}, \\
    \label{eq:S_BAtwice}
    && \big\| (S_\BB \AA)^{\star n_0}(t) \big\|_{\sfL^{2}_{\omega_1} \to \sfL^{\infty}_{\omega_2}} \lesssim e^{-t}.
  \end{eqnarray}
\end{lemma}

\begin{proof}[{\bf Proof of Lemma~\ref{lem:dispersionAS_BA}}]
  The main idea of the proof is that each term in the convolution product exhibits both improved integrability and enhanced decay at large velocities, while also displaying a singularity at $t=0$. Nevertheless, this polynomial singularity remains integrable within the convolution, provided that only a sufficiently small portion of the gain is utilized at each step. The full desired gain can then be recovered by employing a sufficiently large number of convolution factors; that is, for some sufficiently large $n_0$.

  We start by showing that
  \begin{equation}
    \label{eq:estimASBA-relax}
    \forall \, t > 0, \quad \big\| \AA S_\BB(t) \AA \big\|_{\sfL^1  \to \sfL_{\omega_1}^2} \lesssim t^{-\frac{d}{2}}e^{-t}.
  \end{equation}
  Given $0 \le f_\init \in \sfL^{1}(\Omega)$, we have
  \begin{align*} 
    \AA S_\BB(t) \AA f_\init(x,v)
    &= M(v) \int_{\Omega_v} \big(S_\BB \AA f_\init\big)(t,x,w) \dd w\\
    &= M(v) \int_{\Omega_v} e^{-t} \rho[f_\init] (x-wt) M(w) \dd w\\
    &= M(v) \int_{x-\Omega_vt} e^{-t} \rho[f_\init] (y) M\left(\frac{x-y}{t}\right) t^{-d} \dd y\\
    &\leq M(v) \int_{\Omega_x} e^{-t} \rho[f_\init] (y) M\left(\frac{x-y}{t}\right) t^{-d} \dd y,
  \end{align*}
  or equivalently
  \begin{equation}
    \label{eq:estimASBA-relax2}
    \AA S_\BB(t) \AA f_\init(x,v) \leq e^{-t} t^{-d}  M(v) \big( \bar \rho[f_\init] \star M_t \big)(x),
  \end{equation}
  with $\bar \rho[f_\init] = \rho[f_\init] {\bf 1}_{\Omega_x}$ extended to $\R^d$ and $M_t := M(\cdot/t)$. On the one hand, we deduce
  \begin{align*}
    \left\Vert \AA S_\BB(t) \AA f_\init \right\Vert_{\sfL_{\omega_1}^2}
    &\leq t^{-d}e^{-t} \| \bar\rho[f_\init] \star M_t \|_{\sfL^2(\R^d)}  \\
    &\leq t^{-d}e^{-t}  \| M_t \|_{\sfL^2(\R^d)} \| \bar\rho[f_\init]  \|_{\sfL^1(\R^d)}  \\
    &\leq t^{-d}e^{-t} t^{\frac{d}{2}}   \left\Vert M \right\Vert_{\sfL^2} \left\Vert f_\init \right\Vert_{\sfL^1(\Omega)},
  \end{align*}
  what is nothing but~\eqref{eq:estimASBA-relax}. By interpolation, since we have
  \[
    \forall\, t > 0,\quad \big\| \AA S_\BB(t) \AA \big\|_{\sfL^1 \to \sfL^1} \lesssim e^{-t}, 
  \]
  we obtain
  \[
    \forall\, t > 0,\quad \big\| \AA S_\BB(t) \AA \big\|_{\sfL^1 \to \sfL_{\omega_1}^p} \lesssim t^{-\frac{\theta d}{2}} e^{-t},
  \]
  for $p \in (1,2)$ and $\theta := 2- 2/p \in (0,1)$. By choosing $p$ sufficiently close to $1$, we ensure that $\theta d/2 \in (0,1)$, so the singularity as $t \to 0$ is integrable. Combining this with~\eqref{eq:dec-B-relax}, we deduce
  \[
    \forall\, t > 0,\quad \big\| (\AA S_\BB) \star (\AA S_\BB)(t) \big\|_{\sfL^1 \to \sfL_{\omega_1}^p} \lesssim e^{-t}.
  \]
  Proceeding as above, we similarly obtain
  \[
    \forall\, t > 0,\quad \big\| (\AA S_\BB) \star (\AA S_\BB)(t) \big\|_{\sfL_{\omega_1}^{1 + k(p-1)} \to \sfL_{\omega_1}^{1 + (k+1)(p-1)}} \lesssim e^{-t},
  \]
  for $k = 1, \dots, k_0 := \lceil \frac{1}{p-1} - 1 \rceil$, where $k_0$ is the smallest integer greater than or equal to $\frac{1}{p-1} - 1$. By iterating this estimate $k_0 + 1$ times, that is, by convolving $n_0 := 2(k_0 + 1)$ times the function $t \mapsto \AA S_\BB(t)$, we obtain~\eqref{eq:AS_Btwice}.

  \smallskip
  We now show that
  \begin{equation}
    \label{eq:estimASBA-relax-bis}
    \forall \, t > 0, \quad \big\| \AA S_\BB(t) \AA \big\|_{\sfL^2_{\omega_1}  \to \sfL_{\omega_2}^\infty} \lesssim t^{-\frac{d}{2}}e^{-t}.
  \end{equation}
  Using~\eqref{eq:estimASBA-relax2} again and the first estimate in~\eqref{eq:borne-rho&iota}, we have
  \begin{align*}
    \left\Vert \AA S_\BB(t) \AA f_\init \right\Vert_{\sfL_{\omega_2}^\infty}
    &\leq t^{-d}e^{-t} \| \bar\rho[f_\init] \star M_t \|_{\sfL^\infty(\R^d)}  \\
    &\leq t^{-d}e^{-t}  \| M_t \|_{\sfL^2(\R^d)} \| \bar\rho[f_\init]  \|_{\sfL^2(\R^d)}  \\
    &\leq t^{-d}e^{-t} t^{\frac{d}{2}}   \left\Vert M \right\Vert_{\sfL^2(\R^d)} \| \omega_1^{-1} \|_{\sfL^2(\R^d)} \left\Vert f_\init \right\Vert_{\sfL^2_{\omega_1}},
  \end{align*}
  what is nothing but~\eqref{eq:estimASBA-relax-bis}. It implies~\eqref{eq:S_BAtwice} by combining this final estimate with~\eqref{eq:dec-B-relax} and performing a convolution in time along intermediate interpolated steps as above. 
\end{proof}

We now present the proof of~\eqref{eq:dec-iterated-gal-1DprimeBIS} in Proposition~\ref{prop:estSB} for the linear relaxation equation in the case $\Omega_v = \mathbb{R}^d$ with $d \geq 2$. We first establish that
\begin{equation}
  \label{eq:estimASBA-relax-theta1}
  \left\| \AA S_\BB(t) \AA \right\|_{\sfL^1  \to \sfL^{2}_{\bar \omega_1}(\Omega_1 ; \sfL^{1}(\Omega'))} \lesssim t^{-\frac{1}{2}}e^{-t}.
\end{equation}
For a function $f \in \sfL^1(\Omega)$, we define its marginal $g$ on the first components of $x$ and $v$ as in~\eqref{eq:def-marginal}:
\begin{equation*}
  g (x_1,v_1) = \fM[f] := \int_{\R^{2d-2}}  f(x,v) \dd x'\dd v' \in \sfL^1(\R_+ \times \R).
\end{equation*}
We observe that
\begin{align*}
  & \fM [\AA f] = \AA_1 g := \MMM_1(v_1) \rho[g], \quad \text{with } g := \fM (f), \\[2mm]
  & \fM \left[S_\BB(t) f_\init \right] = S_{\BB_1}(t) g_\init, \quad \text{with } g_\init := \fM (f_\init),
\end{align*}
where $S_{\BB_1}(t)$ is the semigroup associated with the following transport equation in dimension $1$:
\begin{align*}
  & \partial_t g + v_1 \partial_{x_1} g = - g \ \text{in} \ (0,\infty) \times \mathbb{R}_+ \times \R \\[2mm]
  & g = 0 \ \text{on} \ (0,\infty) \times \{ 0 \} \times \mathbb{R}_+.
\end{align*}

As a consequence, we have
\begin{equation}
  \label{eq:prev-id}
  \fM \left[ \AA S_\BB(t) \AA f_\init \right] = \AA_1  S_{\BB_1}(t) \AA_1 g_\init.
\end{equation}
Given $0 \le f_\init \in \sfL^1(\Omega)$, we then compute
\begin{eqnarray*}
  \left\| \AA S_\BB(t) \AA f_\init \right\|_{\sfL^{2}_{\bar \omega_1}(\Omega_1 ; \sfL^{1}(\Omega'))}
  &= & \left\| \fM ( \AA S_\BB(t) \AA f_\init) \right\|_{\sfL^2_{\bar \omega_1}}  
  = \left\| \AA_1  S_{\BB_1}(t) \AA_1 g_\init  \right\|_{\sfL^2_{\bar \omega_1}} \\
  &\lesssim& t^{-\frac{1}{2}}e^{-t}
        % \left\Vert M \right\Vert_{\sfL^2(\R)}
        \left\Vert g_\init \right\Vert_{\sfL^1(\R_+ \times \R)} 
  = t^{-\frac{1}{2}}e^{-t}   % \left\Vert M \right\Vert_{\sfL^2(\R)}
      \left\Vert f_\init \right\Vert_{\sfL^1(\Omega)} ,
 \end{eqnarray*}
 where we have used the identity~\eqref{eq:prev-id}, the non-negativity of the operators $\AA$, $S_\BB(t)$, $\AA_1$ and $S_{\BB_1}(t)$, and the estimate~\eqref{eq:estimASBA-relax} for the case $d=1$. The estimate~\eqref{eq:estimASBA-relax-theta1} is thus proved.
 
In a similar fashion, since Lemma~\ref{lem:SBrelaxation-estim1} yields
\begin{equation*}
  \left\| S_{\BB_1}(t) \right\|_{\sfL^2_{\bar \omega_1} \to \sfL^2_{\bar \omega_1}} \lesssim e^{-\kappa_{2,0} t},
\end{equation*}
we deduce
\begin{equation}
  \label{eq:estimSB-relax-theta1}
  \left\| S_\BB(t) \right\|_{\sfL^{2}_{\bar \omega_1}(\Omega_1 ; \sfL^{1}(\Omega')) \to \sfL^{2}_{\bar \omega_1}(\Omega_1 ; \sfL^{1}(\Omega'))} \lesssim e^{-\kappa_{2,0} t}.
\end{equation}
 
It is then straightforward to verify that $\AA_1$ is a bounded operator on $\sfL^2_{\bar \omega_1}$, either by interpolating the two estimates established in Lemma~\ref{lem-relax-AAf} or by repeating its proof. Therefore, as before, we deduce
\begin{equation}
  \label{eq:estimA-relax-theta1}
  \left\| \AA \right\|_{\sfL^{2}_{\bar \omega_1}(\Omega_1 ; \sfL^{1}(\Omega')) \to \sfL^{2}_{\bar \omega_1}(\Omega_1 ; \sfL^{1}(\Omega'))} \lesssim 1.
\end{equation}
By combining~\eqref{eq:estimASBA-relax-theta1},~\eqref{eq:estimSB-relax-theta1}, and~\eqref{eq:estimA-relax-theta1}, and applying time convolution along the intermediate interpolated steps as before, we establish the first estimate in~\eqref{eq:dec-iterated-gal-1DprimeBIS}.

\smallskip
To derive the second estimate in~\eqref{eq:dec-iterated-gal-1DprimeBIS}, we note that the bound~\eqref{eq:estimASBA-relax-theta1} also applies to the adjoint operator $\AA^* S_\BB^* \AA^*$. By duality, this yields
\begin{align}
  \label{eq:estimASBA-relax-theta1-bis}
  \left\| \AA S_\BB(t) \AA \right\|_{\sfL^{2}_{\bar \omega_1}(\Omega_1; \sfL^{\infty}_{\bar \omega'}(\Omega')) \to \sfL^{\infty}_{0,\omega}} \lesssim t^{-\frac{1}{2}} e^{-t}.
\end{align}
Furthermore, by combining Lemma~\ref{lem-relax-AAf}, Lemma~\ref{lem:SBrelaxation-estim1}, and interpolation, we obtain
\begin{equation*}
  \left\| \AA \right\|_{\sfL^{\infty}_{0,\omega} \to \sfL^{\infty}_{0,\omega}} \lesssim 1, \qquad \text{and} \qquad \left\| S_\BB(t) \right\|_{\sfL^{\infty}_{0,\omega} \to \sfL^{\infty}_{0,\omega}} \lesssim e^{-\frac{t}{2}}.
\end{equation*}
Combining these estimates with~\eqref{eq:estimASBA-relax-theta1-bis}, and applying time convolution along the intermediate interpolated steps as before, we conclude that the second bound in~\eqref{eq:dec-iterated-gal-1DprimeBIS} holds.

\subsection{Proof of Proposition~\ref{prop:estSB} for the kinetic Fokker-Planck equation}
\label{subsec-split-FP}

In this subsection, $\mathcal{L}=\nabla_v \cdot ( \nabla_v + v )$ is the Fokker-Planck operator as defined in Subsection~\ref{subsec:setting}, with $\Omega_v=\R^d$ and $d \ge 1$. The decomposition~\eqref{eq:splitting} is defined by
\begin{equation}
  \label{eq:KFP-defA}
  \mathcal{A} f(t,x,v) := \mathfrak{C} \chi_R(v) f(t,x,v),
\end{equation}
where $\mathfrak{C}>0$ is a constant, $\chi_R := \chi(\cdot/R)$ with $\chi \in \mathcal{C}^\infty(\R^d)$ so that ${\bf 1}_{B(0,1)} \le \chi \le {\bf 1}_{B(0,2)}$, and $R>0$ some constant to be chosen later. The semigroup $S_\mathcal{B}$ is associated to the kinetic equation
\begin{equation}
  \label{eq:KFP-eqSB}
  \left\{
    \begin{array}{l}
      \partial_t f + v \cdot \nabla_x f = \Delta_v f + \nabla_v \cdot (v f) - \mathfrak{C} \chi_R f,  \quad \hbox{on}\quad (0,\infty) \times \Omega,   \\[2mm]
      f = 0 \quad \hbox{on}\quad (0,\infty) \times \partial\Omega_+, 
      \qquad f(0,\cdot,\cdot) = f_\init \quad \hbox{in}\quad \Omega.
    \end{array}
  \right.
\end{equation}

We prove the estimates in Proposition~\ref{prop:estSB} through a series of lemmas. Observe that with this choice of $\mathcal{A}$, estimates \eqref{eq:dec-A} are clearly satisfied.
\begin{lemma}
  \label{lem:estSBLpkKFP}
  There are $\mathfrak{C}, R>0$ large enough so that for any $p \in [1,\infty]$, $k \in \{ -1,0,1 \}$ and admissible weight function $\omega$, the solution $f$ to the modified kinetic Fokker-Planck equation~\eqref{eq:KFP-eqSB} in $\sfL^{p}_{k,\omega}$ satisfies the decay estimate
  \begin{equation}
    \label{eq:dec-B-KFP}
    \forall \, t \ge 0, \quad \| f_t \|_{\sfL^p_{k,\omega}} = \| S_\BB(t) f_\init \|_{\sfL^p_{k,\omega}}  \lesssim   e^{- 2t} \| f_\init\|_{\sfL^{p}_{k,\omega}}.
  \end{equation}
  This proves~\eqref{eq:dec-B} for the kinetic Fokker-Planck equation.
\end{lemma}

\begin{proof}[{\bf Proof of Lemma~\ref{lem:estSBLpkKFP}}]
   We restrict our attention to the case $k=1$, as the cases $k=0$ and $k=-1$ can be handled analogously. Let $f \ge 0$ be a solution to~\eqref{eq:KFP-eqSB} and let $p \in (1,\infty)$. We compute
  \begin{equation*}
    \frac{1}{p}  \dt \| f_t \|_{\sfL^{p}_{1,\omega}}^p
    \le \int_{\Omega} \big(- v \cdot \nabla_x  f_t + \LL f_t - \mathfrak{C}\chi_R f_t \big) \, f_t ^{p-1}  \, \wangle{x_1}^p \wangle{v}^p \omega^p  \dd x \dd v.
  \end{equation*}
  On the one hand, we have
  \begin{align*}
    \int_{\Omega} \big(- v \cdot \nabla_x  f_t \big) \, f_t^{p-1}  \, \wangle{x_1}^p \wangle{v}^p \omega^p \dd x \dd v \le \frac{1}{p}
    \int_{\Omega} \frac{x_1 v_1}{\wangle{x_1}^2} f_t^p \wangle{x_1}^p \wangle{v}^p \omega^p  \dd x \dd v,
  \end{align*}
  by performing an integration by parts and using the boundary condition. On the other hand, from \cite[(3.18)]{MisMou}, we have (denoting $\tilde \omega := \wangle{v} \omega$)
  \begin{equation*}
    \int_{\R^d} \big[\big(\LL - \mathfrak{C}\chi_R\big) \, f_t \big]\, f_t ^{p-1}  \, \tilde \omega^p  \dd v = - (p-1) \int_{\R^d}  |\nabla_v (f_t \tilde \omega)|^2 \, (f_t \tilde \omega)^{p-2} \dd v + \int_{\R^d} (f_t \tilde \omega)^p \, \mathsf w \dd v
  \end{equation*}
  with the weight $\mathsf w$ defined as follows in terms of $\tilde \omega$ and $p$:
  \begin{equation}
    \label{eq:Def-psi-mp}
    \mathsf w(v) := \left( 1 - \frac{2}{p'} \right) \frac{\Delta_v \tilde \omega}{\tilde \omega} + \frac{2 }{p'} \frac{|\nabla_v \tilde \omega|^2}{\tilde \omega^2} + \frac{d }{p'} - v \cdot \frac{\nabla_v \tilde \omega }{\tilde \omega} - \mathfrak{C} \,
    \chi_R,
  \end{equation}
  and where $p' := p/(p-1)$ stands for the conjugate exponent. 
  Furthermore, computing
  \begin{align*}
    & \frac{\nabla_v \tilde \omega}{\tilde \omega} (v)=  v \wangle{v}^{-2} + rs v  \wangle{v}^{s-2}, 
    \\[2mm]
    & \frac{\Delta_v \tilde \omega}{\tilde \omega} (v) = d \wangle{v}^{-2} +  r sd \wangle{v}^{s-2} - |v|^2 \wangle{v}^{-4} +  rs^2  |v|^2 \wangle{v}^{s-4} +  (rs)^2  |v|^2 \wangle{v}^{2s-4} ,
  \end{align*}
  we observe that
  \begin{equation*}
    \mathsf w =  (rs)^2 \wangle{\cdot}^{2s-2} - (rs) \wangle{\cdot}^s + \OO(1) - \mathfrak{C} \chi_R.
  \end{equation*}
  Gathering the above pieces of information, we find
  \begin{equation}
    \label{eq:dec-B-KFP-ODE}
    \frac{1}{p}\dt \| f_t \|_{\sfL^{p}_{1,\omega}}^p
    \le - \frac{4(p-1)}{p^2} \big\| \vert \nabla_v (f_t\wangle{x_1} \tilde \omega) \vert^{\frac{p}{2}} \big\|_{\sfL^{2}} ^{2} + \int_\Omega  (f_t\wangle{x_1} \tilde \omega)^p \, \big( |v_1|  + \mathsf w \big) \dd x \dd v.
  \end{equation}
  By choosing $\mathfrak{C}$ and $R > 0$ large enough, we can impose
  \begin{equation*}
    \forall \, v \in \Omega_v, \quad |v_1|  + \mathsf w(v) \le -2 - \nu \wangle{v}^{s},
  \end{equation*}
  where the constant $\nu >0$ is defined by
  \begin{equation*}
    \begin{cases}
      & \nu := \frac{r-1}{2}> 0 \text{ when } s=1, \\[2mm]
      & \nu := \frac{rs}{2} > 0 \text{ when } s \in (1,2), \\[2mm]
      & \nu := r(1-2r) > 0 \text{ when } s=2.
    \end{cases}
  \end{equation*}
  We deduce~\eqref{eq:dec-B-KFP} when $p \in (1,\infty)$ thanks to the Gronwall lemma. Finally,~\eqref{eq:dec-B-KFP} holds when $p \in [1,\infty]$ by a continuity argument as $p \to 1$ and $p\to\infty$.
\end{proof}

We point out again that the estimate~\eqref{eq:dec-B-KFP} enables us to establish the existence of the associated semigroup $S_\BB$ in the spaces $\sfL^2_{k,\omega}$, for any admissible weight function $\omega$ and $k \in \{-1,0,1\}$, by following the argument presented in the proof of Lemma~\ref{lem:exist}.

\begin{lemma}
  \label{lem:estSBL1LinftyKFP}
  The semigroup $S_{\mathcal B}$ satisfies the ultracontractivity property
  \begin{equation}
    \label{eq:dec-U-KFP}
    \forall \, t >0,\quad \big\| S_\BB (t) \big\|_{\sfL^{1}_{\omega} \to \sfL^{\infty}_{\omega}}  
    \lesssim \frac{e^{-t}}{t^\Theta},
  \end{equation}
  for some $\Theta >0$ and for any admissible weight function $\omega$. We deduce that~\eqref{eq:dec-iterated-gal} holds for the kinetic Fokker-Planck equation, for $n$ large enough depending on $\Theta$.
\end{lemma}

The estimate~\eqref{eq:dec-U-KFP} is a variant of~\cite[Theorem 1.3]{carrapatoso2024KR}. In this previous work, \cite[Theorem 1.3]{carrapatoso2024KR} is a decay estimate on the whole semigroup, not the modified semigroup $S_\BB$ (and assuming that $\Omega_x$ is bounded), which explains the slightly weaker decay estimates. The proof is based on an idea introduced in \cite{MR1166050,MR3591133,sanchez2024KR} and developed in \cite{carrapatoso2024KFP,carrapatoso2024landau,carrapatoso2024KR,sanchez2024voltage} motivated by the gains of integrability established in~\cite{MR2068847,MR3923847} for Kolmogorov equations with rough coefficients. We present a sketch of the proof and we refer to~\cite{carrapatoso2024KFP,carrapatoso2024landau,carrapatoso2024KR,sanchez2024voltage} for more details. The fact that~\eqref{eq:dec-U-KFP} and $\AA$ is bounded imply~\eqref{eq:dec-iterated-gal} is now well understood and has been widely discussed in previous works on the extension method introduced in \cite{MR2197542} and developed in \cite{MR3779780,MisMou,MR3465438}.

\begin{proof}[{\bf Proof of Lemma~\ref{lem:estSBL1LinftyKFP}}] We do not optimize all exponents in this proof to prioritize conciseness and clarity; such optimization can be easily done by the reader if necessary. We first prove~\eqref{eq:dec-U-KFP}. Let $f \ge 0$ be a solution to~\eqref{eq:KFP-eqSB}, $\omega=\omega(v)$ an admissible weight function, $T>0$ and  $\varphi=\varphi(t)$ a non-negative function in $\mathcal D((0,T))$. We define the \emph{distance to the boundary}
\begin{equation*}
  \delta(x) := \min(1,x_1)
\end{equation*}
We now assert that, by replicating the proof of Lemma~\ref{lem:estSBLpkKFP} with $p=2$ and the weight
\begin{equation*}
  \omega \left( 1 - \frac{v_1}{2\wangle{v}^2} \delta(x)^{\frac{1}{2}} \right)  \varphi,
\end{equation*}
as well as utilizing
\begin{equation*}
  \begin{cases}
    & \frac{1}{2} \le \left( 1 - \frac{v_1}{2\wangle{v}^2} \delta(x)^{\frac{1}{2}} \right)  \le 1, \\[2mm]
    \displaystyle
    & - v \cdot \nabla_x \left( 1 - \frac{v_1}{2\wangle{v}^2} \delta(x)^{\frac{1}{2}} \right)  = \frac{1}{4 \delta^{\frac{1}{2}}} \frac{v_1^2}{\wangle{v}^2}  {\bf 1}_{x_1 \in (0,1)}, \\[2mm]
    & v \cdot \nabla_v \omega = r s |v|^2 \wangle{v}^{s-2} \omega,
  \end{cases}
\end{equation*}
we get
\begin{equation}
  \int_0^T \!\! \int_\Omega \left[  f^2  \omega^2  \left(  \frac{1}{\delta^{\frac{1}{2}}}  \frac{v_1^2}{\wangle{v}^2} + \wangle{v}^s \right) +  |\nabla_v(f \omega )|^2
  \right] \varphi^2 \dd x \dd v \dd t 
  \label{eq:Ultracontractivity1}
  \lesssim \left\| f \omega (\varphi + |\varphi'|) \right\|_{\sfL^2((0,T) \times \Omega)} ^2.
\end{equation}
For further details, we refer the reader to~\cite[Section~3]{carrapatoso2024KFP}.

Then we prove the following interpolation inequality directly inspired from~\cite[Lemma~4.2]{carrapatoso2024KR}:
\begin{equation}
  \label{eq:EstimL2poids}
  \int_{\Omega}  \frac{F^2}{\delta^{\beta}} \dd x \dd v \lesssim
  \int_{\Omega} |\nabla_v F|^2 \dd x \dd v + \int_{\Omega} F^2 \left(  \frac{1}{\delta^{\frac{1}{2}}} \frac{v_1^2}{\wangle{v}^2}   +   \wangle{v}^s \right) \dd x \dd v
\end{equation}
with $\beta>0$ small enough in terms of $d$ and $1/s$ (the exact condition can easily be calculated below). It follows from splitting the integral into three regions:
\begin{enumerate}
\item a region $\mathcal{R}_1$ determined by $\wangle{v} \delta^{\frac{\beta}{s}} \ge 1$ on which
  \begin{equation*}
    \int_{\mathcal{R}_1} \frac{F^2}{\delta^{\beta}} \dd x \dd v \le \int_{\Omega}  \frac{F^2}{\delta^{\beta}} \left( \wangle{v} \delta^{\frac{\beta}{s}} \right)^s \dd x \dd v = \int_{\Omega} F^2 \wangle{v}^s \dd x \dd v,
  \end{equation*}
\item a region $\mathcal{R}_2$ determined by $\wangle{v} \delta^{\frac{\beta}{s}} < 1$ and $|v_1| > \delta^{\frac{1}{4}-\beta(\frac{1}{2}+\frac{1}{s})}$ on which
  \begin{equation*}
    \int_{\mathcal{R}_2}  \frac{F^2}{\delta^{\beta}} \dd x \dd v \le \int_{\Omega} \frac{F^2}{\delta^{\beta}} \left( |v_1| \delta^{\beta(\frac{1}{2}+\frac{1}{s})-\frac{1}{4}} \right)^2 \left( \wangle{v} \delta^{\frac{\beta}{s}} \right)^{-2} \dd x \dd v = \int_{\Omega} F^2 \frac{1}{\delta^{\frac{1}{2}}} \frac{v_1^2}{\wangle{v}^2} \dd x \dd v,
  \end{equation*}
\item a region $\mathcal{R}_3$ determined by $\wangle{v} \delta^{\frac{\beta}{s}} < 1$ and $|v_1| \le \delta^{\frac{1}{4}-\beta(\frac{1}{2}+\frac{1}{s})}$ on which
  \begin{align*}
    \int_{\mathcal{R}_3}  \frac{F^2}{\delta^{\beta}} \dd x \dd v
    & \le \int_{\Omega_x} \delta^{-\beta} \left\| F(x,\cdot) \right\|_{\sfL^{\frac{2d}{d-2}}(\R^d)}^2   \left| \left\{ \wangle{v} \delta^{\frac{\beta}{s}} < 1 \text{ and } |v_1| \le \delta^{\frac{1}{4}-\beta(\frac{1}{2}+\frac{1}{s})} \right\} \right|^{\frac{2}{d}} \dd x \\
    & \lesssim \int_{\Omega_x}  \delta^{-\beta} \left( \int_{\R} \big|\nabla_v F(x,\cdot)\big|^2 \dd v \right) \left| \left\{ \wangle{v} \delta^{\frac{\beta}{s}} < 1 \text{ and } |v_1| \le \delta^{\frac{1}{4}-\beta(\frac{1}{2}+\frac{1}{s})} \right\} \right|^{\frac{2}{d}} \dd x \\
    & \lesssim \int_{\Omega}  \big|\nabla_v F\big|^2 \dd x \dd v,
  \end{align*}
\end{enumerate}
where we have used the Gagliardo-Nirenberg-Sobolev inequality in this last estimate (when $d=1$ or $d=2$, the $\sfL^{2d/(d-2)}$ norm on $F(x,\cdot)$ is replaced by $\sfL^\infty$) and the smallness condition on $\beta$.

By applying~\eqref{eq:EstimL2poids} to $F:=f\omega \varphi$, we deduce from~\eqref{eq:Ultracontractivity1}:
\begin{equation}
  \label{eq:hypo1}
  \int_0^T\!\!\int_\Omega f^2 \omega^2 \varphi^2 \Bigl( \frac{1}{\delta^{\beta}} +  \wangle{v}^s  \Bigr) \dd t \dd x \dd v 
  \lesssim \left\| f \omega (\varphi +|\varphi '|) \right\|^2 _{\sfL^2((0,T) \times \Omega)},
\end{equation}

We then use the hypoelliptic gain of integrability in the form~\cite{MR2068847,MR3923847}, but after a localization \emph{away from the boundary}. Consider the new unknown
\begin{equation*}
  F := \wangle{v}^{-2} \delta \varphi \omega f.
\end{equation*}
It satisfies the equation
\begin{equation*}
  \partial_t F + v \cdot \nabla_x F - \Delta_v F = S
\end{equation*}
over $\R^{2d+1}$ thanks to the localization, with a source term $S$ that satisfies 
\begin{equation*}
  \left\| S \right\|_{\sfL^2(\R^{2d+1})} \lesssim \left\| f \omega (\varphi + |\varphi'|) \right\|_{\sfL^2((0,T) \times \Omega)}
\end{equation*}
thanks to the following properties (the third one follows from Definition~\ref{def:admissible})
\begin{equation*}
   \wangle{v}^{-2} \left| v \cdot \nabla_x \delta \right| \lesssim 1 \quad \text{ and } \quad \left| v \cdot \nabla_v \wangle{v}^{-2} \right| \lesssim 1 \quad \text{ and } \wangle{v}^{-2} \left| v \cdot \nabla_v \omega \right| \lesssim \omega
\end{equation*}
and the energy estimate on $f$. Then standard formula for the fundamental solution to the Kolmogorov equation $\partial_t F + v \cdot \nabla_x F = \Delta_v F$ (see~\cite{MR2068847} or more recently~\cite[Lemma~10]{MR4453413}) yields
\begin{equation}
  \label{eq:hypo2}
  \left\| \wangle{v}^{-2} \delta \omega \varphi f \right\|_{\sfL^{p}(\R_+ \times \R^{2d})} = \left\| F \right\|_{\sfL^{p}(\R^{2d+1})} \lesssim_p \left\| f \omega (\varphi + |\varphi'|) \right\|_{\sfL^2((0,T) \times \Omega)}
\end{equation}
for $p \in (2,2+1/d)$.

By interpolating between the two inequalities~\eqref{eq:hypo1} and~\eqref{eq:hypo2}, we deduce
\begin{equation}
  \label{eq:hypo3}
  \left\| f \omega \varphi \right\|_{\sfL^{\tilde p}(\R_+ \times \Omega)}  \lesssim
  \left\| f \omega (\varphi + |\varphi'|) \right\|_{\sfL^2((0,T) \times \Omega)}
\end{equation}
for some $\tilde p \in (2,p) \subset (2,2+1/d)$. We then use the Hölder inequality
\begin{align*}
  \left\| f \omega \varphi \right\|_{\sfL^2((0,T) \times \Omega)}
  & \lesssim \left\| f \omega \varphi (\varphi + |\varphi'|)^{-1} \right\|_{\sfL^{\tilde p}((0,T) \times \Omega)} ^{\frac{\tilde p}{2(\tilde p-1)}}  \left\| f \omega \varphi (\varphi + |\varphi'|)^{\frac{\tilde p-2}{\tilde p}} \right\|_{\sfL^1((0,T) \times \Omega)} ^{1-\frac{\tilde p}{2(\tilde p-1)}} \\
  & \lesssim \left\| f \omega \varphi \right\|_{\sfL^2((0,T) \times \Omega)} ^{\frac{\tilde p}{2(\tilde p-1)}}  \left\| f \omega \varphi (\varphi + |\varphi'|)^{\frac{\tilde p-2}{\tilde p}} \right\|_{\sfL^1((0,T) \times \Omega)} ^{1-\frac{\tilde p}{2(\tilde p-1)}}
\end{align*}
which implies 
\begin{equation*}
  \left\| f \omega \varphi \right\|_{\sfL^2((0,T) \times \Omega)}
  \lesssim \left\| f \omega \varphi (\varphi + |\varphi'|)^{\frac{\tilde p-2}{\tilde p}} \right\|_{\sfL^1((0,T) \times \Omega)}.
\end{equation*}
Then we use the decay estimate~\eqref{eq:dec-B-KFP} twice to control earlier times by later times from above ($T$ is bounded say by $1$):
\begin{align*}
  \left\| f_T \right\|_{\sfL^2_\omega}
  & \lesssim  \left\| f \omega \varphi \right\|_{\sfL^2((0,T) \times \Omega)} \\
  & \lesssim \left\| f \omega \varphi (\varphi + |\varphi'|)^{\frac{\tilde p-2}{\tilde p}} \right\|_{\sfL^1((0,T) \times \Omega)} \\
  & \lesssim \left( \sup_{t \in (0,T)} \varphi (\varphi + |\varphi'|)^{\frac{\tilde p-2}{\tilde p}} \right) \left\| f_\init  \right\|_{\sfL^1_{\omega}} \\
  & \lesssim T^{-\frac{\Theta}{2}} \left\| f_\init \right\|_{\sfL^1_{\omega}}
\end{align*}
for some $\Theta >0$. Using then the same duality argument as in  Subsection~\ref{subsec:pw-dual}, we deduce~\eqref{eq:dec-U-KFP} for $t \in (0,1]$, which in combination with the decay estimate~\eqref{eq:dec-B-KFP} proves~\eqref{eq:dec-U-KFP}. 

Finally, since $\AA$ is bounded from $\sfL^p_\omega$ to $\sfL^p_{\tilde \omega}$ for any admissible weight functions $\omega$ and $\tilde \omega$ and for all $p \in [1, \infty]$, we can deduce~\eqref{eq:dec-iterated-gal} from~\eqref{eq:dec-U-KFP} and~\eqref{eq:dec-B-KFP} by choosing $\mathfrak{a} := \frac{1}{2}$ and $n := \lfloor \Theta \rfloor + 2$. This follows by an argument analogous to the one used at the end of the proof of Lemma~\ref{lem:dispersionAS_BA}: namely, the $n$-fold convolution of the function $t \mapsto t^{-\frac{\Theta}{n}}$ remains bounded on $(0,1]$, and we employ $n$ intermediate interpolation steps to bridge the integrability gap between $\sfL^1$ and $\sfL^\infty$. For similar arguments, see also~\cite{MR3779780,MR3465438}.
\end{proof}

We now present the proof of~\eqref{eq:dec-iterated-gal-1DprimeBIS} in Proposition~\ref{prop:estSB} for the kinetic Fokker-Planck equation in the case $\Omega_v = \mathbb{R}^d$ with $d \geq 2$. The decomposition~\eqref{eq:splitting} is slightly modified compared to~\eqref{eq:KFP-defA}-\eqref{eq:KFP-eqSB}. While the overall formulas remain the same, the localization function $\chi_R$ is now defined as
\begin{equation*}
  \chi_R(v) := \chi\left(\frac{v_1}{R}\right) \chi\left(\frac{|v'|}{R}\right),
\end{equation*}
where $\chi \in \mathcal{C}^\infty(\R)$ satisfies ${\bf 1}_{[-1,1]} \leq \chi \leq {\bf 1}_{[-2,2]}$. This choice of $\chi_R$ serves to decouple the first component from the remaining components. It is straightforward to verify that all the properties and estimates established in this Subsection~\ref{subsec-split-FP} continue to hold with this modified definition of $\AA$ and $\BB$. We define the weight function $\tilde \omega(v') := e^{\wangle{v'}}$ and introduce the corresponding weighted marginal operator $\tilde \fM$ by
\begin{equation*}
  \tilde \fM : f \in \sfL^1_{\tilde \omega}(\Omega) \mapsto \tilde \fM f(x_1,v_1) := \int_{\R^{d-1} \times \R^{d-1}} f(x,v)\, \tilde \omega(v')\, \dd x'\dd v' \in \sfL^1(\R_+ \times \R).
\end{equation*}
Let $f \in \sfL^1_{\tilde \omega}(\Omega)$ and denote its weighted marginal by $g := \tilde \fM f$. We observe that
\begin{align*}
  \tilde \fM [\AA f](x_1,v_1)
  &= \mathfrak{C} \int_{\R^{d-1} \times \R^{d-1}} \chi\left(\frac{v_1}{R}\right) \chi\left(\frac{|v'|}{R}\right) f(x,v)\, \tilde \omega(v')\, \dd x'\dd v' \\
  &\leq \mathfrak{C} \int_{\R^{d-1} \times \R^{d-1}} \chi\left(\frac{v_1}{R}\right) f(x,v)\, \tilde \omega(v')\, \dd x'\dd v' \\
  &= \mathfrak{C} \chi\left(\frac{v_1}{R}\right) g(x_1,v_1) =: (\AA_1 g)(x_1,v_1).
\end{align*}

\smallskip
We define the weight functions $m(v_1) := e^{\wangle{v_1}}$ and $\omega(v) := m(v_1) \tilde \omega(v')$. The function $f := S_\BB(t) f_\init \in C([0,\infty); \sfL^1_\omega)$ satisfies~\eqref{eq:KFP-eqSB}. Consequently, the function $g(t,x_1,v_1) := \tilde \fM S_\BB(t) f_\init \in C([0,\infty); \sfL^1_m)$ satisfies
\begin{equation*}
  \partial_t g + v_1 \partial_{x_1} g = \partial^2_{v_1v_1} g + \partial_{v_1} (v_1 g) + \int_{\R^{d-1} \times \R^{d-1}} f \tilde \omega \mathsf w \dd x'\dd v',
\end{equation*}
where the weight $\mathsf w(v_1,v')$ is given by
\begin{equation*}
  \mathsf w(v_1,v') := \frac{\Delta_{v'} \tilde \omega}{\tilde \omega} - \frac{v' \cdot \nabla_{v'} \tilde \omega}{\tilde \omega} - \mathfrak C \chi_R(v_1,v').
\end{equation*}
We compute
\begin{eqnarray*}
  \mathsf w(v_1,v')
  &=& \frac{d}{\wangle{v'}} - \frac{|v'|^2}{\wangle{v'}^3} + \frac{|v'|^2}{\wangle{v'}^2} - \frac{|v'|^2}{\wangle{v'}} - \mathfrak{C} \chi_R(v_1,v') \\
  &\le& d+1 - \wangle{v'} + \mathfrak{C} {\bf 1}_{|v'| \ge R} - \mathfrak{C} \chi\left( \frac{v_1}{R} \right) \\
  &\le& d+1 - \mathfrak{C} \chi\left( \frac{v_1}{R} \right),
\end{eqnarray*}
by choosing $R > \mathfrak{C} > 0$ sufficiently large. We deduce that $g$ satisfies
\begin{equation*}
  \partial_t g + v_1 \partial_{x_1} g \le \partial^2_{v_1v_1} g + \partial_{v_1} (v_1 g) + \left[d+1 - \mathfrak{C} \chi\left(\frac{v_1}{R}\right) \right] g.
\end{equation*}
The maximum principle then implies $g \le S_{\BB_1}(t) g_{\text{init}}$, where $S_{\BB_1}$ is the semigroup associated with the one-dimensional modified kinetic Fokker-Planck equation
\begin{equation}
  \label{eq:KFP-eqSB-dim1}
  \left\{
    \begin{array}{l}
      \partial_t g + v_1 \partial_{x_1} g  =\partial^2_{v_1v_1}g + \partial_{v_1} (v_1 g) + \left[ d+1 - \mathfrak{C} \chi\left(\frac{v_1}{R}\right)\right] g \quad \text{on}\quad (0,\infty) \times \R_+ \times \R,   \\[3mm]
      g = 0 \quad \text{on}\quad (0,\infty) \times \{0 \} \times (0,\infty).
    \end{array}
  \right.
\end{equation}
It is straightforward to verify that the conclusions of Lemma~\ref{lem:estSBLpkKFP} and Lemma~\ref{lem:estSBL1LinftyKFP} remain valid for this analogous semigroup: if $R > \mathfrak{C} > 0$ are sufficiently large, there exists $\Theta > 0$ such that
\begin{equation*}
  \forall \, t > 0, \quad \| S_{\BB_1} (t) \|_{\sfL^{p}_{m} \to \sfL^{p}_{m}}  +  t^\Theta \| S_{\BB_1} (t) \|_{{\sfL^{1}_{m}  \to \sfL^{\infty}_{m}}}  \lesssim e^{-t}
\end{equation*}
for all $p \in [1,\infty]$. By analogous reasoning as above, there exists an integer $n_0 \ge 1$ such that
\begin{equation*}
  \left\| (\AA_1 S_{\mathcal B_1})^{\star n_0}(t) \right\|_{\sfL^{1} \to \sfL^{2}_{m}} \lesssim e^{- \mathfrak{a} t}.
\end{equation*}
Using the pointwise inequality
\begin{equation*}
  \tilde \fM (\AA S_{\mathcal B})^{\star n_0} f_\init \le (\AA_1 S_{\mathcal B_1})^{\star n_0} \tilde \fM f_\init \quad \text{for any } 0 \le f_\init \in \sfL^1_\omega(\Omega),
\end{equation*}
we obtain
\begin{equation*}
  \left\| (\AA S_{\mathcal B})^{\star n}(t) \right\|_{\sfL^{1} \to \sfL^{2}_{\bar \omega_1}(\Omega_1 ; \sfL^{1}(\Omega'))} \lesssim e^{- \mathfrak{a} t},
\end{equation*}
with $n := n_0+1$, which is the first estimate in~\eqref{eq:dec-iterated-gal-1DprimeBIS}, and where we recall that $\bar \omega_1(v_1) = \MMM_1(v_1)^{-1/2}$. The second estimate in~\eqref{eq:dec-iterated-gal-1DprimeBIS} follows by duality, as before.

\subsection{Proof of Proposition~\ref{prop:estSB} for the kinetic Laplace-Beltrami equation}
\label{subsec-split-LB}

In this subsection, $\mathcal{L} = \Delta_v ^{\text{LB}}$ is the Laplace-Beltrami operator as defined in Subsection~\ref{subsec:setting}, $d \ge 2$ and $\Omega_v=\mathbb S^{d-1}$. We recall that for this operator, no weight function in $v$ is needed since the velocity domain is bounded. The decomposition~\eqref{eq:splitting} is defined by
\begin{equation*}
  {\mathcal A} f  := \mathfrak{C}f
 \end{equation*}
with $\mathfrak{C}>0$ to be chosen later. The semigroup $S_\mathcal{B}$ is associated to the kinetic equation
 \begin{equation}
    \label{eq:LB-eqSB}
    \left\{
      \begin{array}{l}
      \partial_t f = - v \cdot \nabla_x f + \Delta_v^{\text{LB}} f - \mathfrak{C} f,  \quad \hbox{on}\quad (0,\infty) \times \Omega,   \\[3mm]
       f = 0 \quad \hbox{on}\quad (0,\infty) \times \partial\Omega_+, 
       \qquad f(0,\cdot,\cdot) = f_\init \quad \hbox{in}\quad \Omega.
      \end{array}
      \right.
    \end{equation}
We prove the estimates in Proposition~\ref{prop:estSB} through a series of lemmas. Observe again that with this choice of $\mathcal{A}$, equation~\eqref{eq:dec-A} is clearly satisfied.

\begin{lemma}\label{lem:estSBLpkLB}
  There exists $\mathfrak{C} > 0$ such that for any $p \in [1,\infty]$, $k \in \{ -1,0,1\}$ and any solution $f$ to the modified kinetic Laplace-Beltrami  equation~\eqref{eq:LB-eqSB} in $\sfL^p_k$, the following decay holds
  \begin{equation}
    \label{eq:dec-B-LB}
    \forall \, t \ge 0, \quad \| f_t \|_{\sfL^{p}_{k}} = \| S_\BB(t) f_\init \|_{\sfL^{p}_{k}} \lesssim   e^{- 2t} \| f_\init \|_{\sfL^{p}_{k}},
  \end{equation}
  This proves~\eqref{eq:dec-B} for the kinetic Laplace-Beltrami equation.
\end{lemma}

\begin{proof}[{\bf Proof of Lemma~\ref{lem:estSBLpkLB}}]
  We compute,
  \begin{align*}
    \dt \int_\Omega f_t ^p \wangle{x_1}^{kp} \dd x \dd v
    & = \int_\Omega p f_t^{p-1} \left( - v \cdot \nabla_x f_t + \Delta_v^{\text{LB}} f_t - \mathfrak{C}f_t \right) \wangle{x_1}^{kp} \dd x \dd v \\
    & \le \int_\Omega f_t ^p \left( v \cdot \nabla_x (\wangle{x_1}^{kp}) - p \mathfrak{C} \wangle{x_1}^{kp}\right) \dd x \dd v  \\
    & \le \int_\Omega f_t ^p \left( - p \mathfrak{C} + |k|\wangle{x_1}^{-1}  \right)  \wangle{x_1}^{kp} \dd x \dd v  \\
    &\le \left(|k| - p \mathfrak{C}\right)  \int_\Omega f_t ^p \wangle{x_1}^{kp} \dd x \dd v,
  \end{align*}
  where we have used that the Laplace-Beltrami operator is non-positive in this space by integration by parts. The Gronwall lemma then implies
  \begin{equation*}
    \| f_t \|_{\sfL^{p}_{k}} \lesssim e^{- \left( \mathfrak{C} - \frac{|k|}{p} \right) t} \| f_\init \|_{\sfL^{p}_{k}},
  \end{equation*}
  which concludes the proof by choosing $\mathfrak{C}$ large enough.
\end{proof}

\begin{remark}
  The estimate~\eqref{eq:dec-B-LB} allows to deduce the existence of the associated semigroup $S_\BB$ in the spaces $\sfL^{2}_k$ for any $k \in \{-1,0,1\}$ by following the proof of Lemma~\ref{lem:exist}.
\end{remark}

\begin{lemma}
  \label{lem:estSBL1LinftyLB}
  The semigroup $S_{\mathcal B}$ satisfies the ultracontractivity property
  \begin{equation}
    \label{eq:dec-U-LB}
    \forall \, t >0,\quad \big\| S_\BB (t) \big\|_{\sfL^{1} \to \sfL^{\infty}}
    \lesssim \frac{e^{-t}}{t^\Theta},
  \end{equation}
  for some $\Theta >0$. We deduce that~\eqref{eq:dec-iterated-gal} holds for the kinetic Laplace-Beltrami equation, for $n$ large enough depending on $\Theta$.
\end{lemma}

\begin{proof}[{\bf Proof of Lemma~\ref{lem:estSBL1LinftyLB}}]
  The proof follows exactly the same steps as in the case of the kinetic Fokker-Planck equation, i.e., the proof of Lemma~\ref{lem:estSBL1LinftyKFP}. A few adjustments are needed since $\Omega_v$ is a proper submanifold of $\mathbb{R}^d$. We only discuss the few differences. We use the distance to the boundary function $\delta$ again, but the weight for the first estimate is now
  \begin{equation*}
    \left( 1 - \frac{v_1}{2} \delta(x)^{\frac{1}{2}} \right)  \varphi,
  \end{equation*}
  and the resulting estimate is 
  \begin{equation}
    \label{eq:Ultracontractivity1LB}
    \int_0^T \!\! \int_\Omega \left[  f^2 \left(  \frac{v_1^2}{\delta^{\frac{1}{2}}} + 1 \right) +  |\nabla_v^{\mathbb{S}^{d-1}} f|^2 \right] \varphi^2 \dd t \dd x \dd v \lesssim \left\| f (\varphi + |\varphi'|) \right\|^2 _{\sfL^2(\R_+ \times \Omega)}.
  \end{equation}

  The second key estimate, the interpolation estimate, is 
  \begin{equation*}
    \int_{\Omega} \frac{F^2 }{\delta^{\beta}} \dd x \dd v\lesssim
    \int_\Omega |\nabla_v^{\mathbb{S}^{d-1}} F|^2  \dd x \dd v + \int_\Omega F^2 \left(  \frac{v_1^2 }{\delta^{\frac{1}{2}} } + 1\right) \dd x \dd v,
  \end{equation*}
  with any $\beta \in (0,1/2d)$. It follows from splitting the integral into three regions:
  \begin{enumerate}
    \item a region $\mathcal{R}_1$ determined by $\delta \ge \alpha := e^{-(\frac{1}{4}-\frac{\beta}{2})^{-1}}$ on which
      \begin{equation*}
        \int_{\mathcal{R}_1} \frac{F^2}{\delta^{\beta}} \dd x \dd v \le \frac{1}{\alpha^\beta} \int_{\Omega} F^2,
      \end{equation*}
    \item a region $\mathcal{R}_2$ determined by $\delta < \alpha$ and $|v_1| > \delta^{\frac{1}{4}-\frac{\beta}{2}}$ on which
      \begin{equation*}
        \int_{\mathcal{R}_2}  \frac{F^2}{\delta^{\beta}} \dd x \dd v \le \int_{\Omega} F^2 \frac{v_1^2}{\delta^{\frac{1}{2}}} \dd x \dd v,
      \end{equation*}
    \item a region $\mathcal{R}_3$ determined by $\delta < \alpha$ and $|v_1| \le \delta^{\frac{1}{4}-\frac{\beta}{2}}$ on which
      \begin{align*}
        \int_{\mathcal{R}_3} \frac{F^2}{\delta^{\beta}} \dd x \dd v
        & \le \int_{\Omega_x}  \delta^{-\beta} \left\| F(x,\cdot) \right\|_{\sfL^{\frac{2(d-1)}{d-3}}}^2 % \left( \int_{\mathbb S^{d-1}} F^{\frac{2d}{d-2}} \dd v \right)^{\frac{d-2}{d}}
          \left| \left\{ \delta < 
          \alpha\text{ and } |v_1| \le \delta^{\frac{1}{4}-\frac{\beta}{2}} \right\} \right|^{\frac{2}{d-1}} \dd x\\
        & \lesssim \int_{\Omega_x}  \delta^{-\beta} \left( \int_{\mathbb S^{d-1}} \big|\nabla_v ^{\mathbb S^{d-1}} F(x,\cdot) \big|^2 \dd v \right) \left| \left\{ \delta < \alpha \text{ and } |v_1| \le \delta^{\frac{1}{4}-\frac{\beta}{2}} \right\} \right|^{\frac{2}{d-1}} \dd x \\
        & \lesssim \int_{\Omega}  \big|\nabla_v ^{\mathbb S^{d-1}}F\big|^2 \dd x \dd v,
      \end{align*}
    \end{enumerate}
    We have once again used the Gagliardo-Nirenberg-Sobolev inequality, together with the smallness condition on $\beta$. In the case $d=2$, the $\sfL^{\frac{2(d-1)}{d-3}}$ norm of $F(x,\cdot)$ is replaced by the $\sfL^\infty$ norm; for $d=3$, it is replaced by the $\sfL^q$ norm, where $q \gg 1$ is finite but can be taken arbitrarily large. 
  
    We deduce the third key estimate 
    \begin{eqnarray*}
      \int_0^T\!\!\int_\Omega  f^2 \Bigl( \frac{1}{\delta^{\beta}} +  1  \Bigr)  \varphi^2 \dd t \dd x \dd v 
      \lesssim \left\| f (\varphi + |\varphi'|) \right\|^2 _{\sfL^2(\R_+ \times \Omega)}.
    \end{eqnarray*}
    
    In the case of the kinetic Fokker-Planck equation, the next step is a gain of integrability based on the integrability estimates on the fundamental solution to the Kolmogorov equation. We do not know explicit formulas for the fundamental solutions to the Laplace-Beltrami equation so we replace this step by a regularization estimate on the solutions to the kinetic Laplace-Beltrami equation proven in~\cite[Proposition~4.7]{Bouin_Frouvelle}:
    \begin{equation*}
      %\forall \, t \in [0,1], \quad
      \Vert f_t \Vert_{\sfL^2} ^2 + \gamma_1 t \Vert\nabla_v f_t \Vert_{\sfL^2}^2+2 \gamma_2 t ^2\left\langle\nabla_v f_t, \nabla_x f_t\right\rangle_{\sfL^2} 
      + \gamma_3 t ^3\Vert P_{v}^\perp\nabla_x f_t\Vert^2_{\sfL^2}+ \gamma_4 t ^4\Vert\nabla_x f_t\Vert_{\sfL^2}^2 \lesssim\Vert f_\init \Vert_{\sfL^2}^2,
    \end{equation*}
    for any $t \in [0,1]$ and for explicit constants $\gamma_1$, $\gamma_2$, $\gamma_3$, $\gamma_4 > 0$ such that $\gamma_2 < \gamma_1 \gamma_3$, and where $P_v^\perp$ is the orthogonal projection on $v^\perp$. This implies, combined with Sobolev inequalities that
    \begin{equation*}
      \Vert f_t \Vert_{\sfL^p} \lesssim t^{-\Theta_0} \Vert f_\init \Vert_{\sfL^2}
    \end{equation*}
    for some $\Theta_0>0$ and some $p>2$. The rest of the argument is exactly similar and leads to
    \begin{equation*}
      \forall \, t \in (0,1], \quad \| f (t) \|_{\sfL^\infty}   \lesssim t^{-\Theta} \| f_\init \|_{\sfL^{1}}
    \end{equation*}
    for some $\Theta>0$.
    
    Finally, since $\AA$ is bounded from $\sfL^p$ to $\sfL^p$ for any $p \in [1,\infty]$, we can deduce~\eqref{eq:dec-iterated-gal} from~\eqref{eq:dec-U-LB} and~\eqref{eq:dec-B-LB} with $\mathfrak a := \frac{1}{2}$ and $n := \lfloor  \Theta \rfloor + 2$ by the same argument as for the kinetic Fokker-Planck equation.
\end{proof}

\section{Localization of the mass in dimension one}
\label{sec:loc-mass}

In this section, we focus on the special case $d=1$ and prove Theorem~\ref{theo:main2}. We consider only the relaxation operator and the Fokker-Planck operator, as the Laplace-Beltrami operator is not defined in one dimension.
 
 \subsection{Proof of~\eqref{eq:loc-mass-int} in Theorem~\ref{theo:main2}}
\label{subsec:loc-mass-relax}

Let $d=1$, let $\omega$ be an admissible weight function, and let the initial data $f_{\init}$ belong to $\sfL^{1}_{1,\omega} \cap \sfL^{\infty}_{-1,\omega}$, with $\mathcal{N}_0 > 0$, where $\mathcal{N}_0$ is defined in~\eqref{eq:defMMM0}.
% \subsubsection{Proof of~\eqref{eq:loc-mass-int} for the linear relaxation equation}
% \label{ss:int}
Recall the notation $\rho_t := \rho[f_t]$ and $\iota_t := \iota[f_t]$ as defined in~\eqref{eq:rf}. By applying either~\eqref{eq:bddrhox1-casborne} in the case of bounded velocities, or~\eqref{eq:uniformL11} when velocities are unbounded, we obtain
\begin{equation}
  \label{eq:estim-rhoL11}
  \forall\, t \ge 0, \quad \int_{\R_+} x\, \rho_t \dd x \lesssim 1.
\end{equation}
 
Hölder's inequality and~\eqref{eq:decay-int} imply 
\begin{equation}
  \label{eq:estim-kinetic-rho}
  \forall\, t \ge 0, \quad   \| \rho_t \|_{\sfL^{2}(\mathbb{R}_+)}^2  \lesssim  t^{-\frac{3}{2}}.
\end{equation}
Together with \eqref{eq:estim-rhoL11} and the interpolation inequality  \eqref{eq:interpol-L1L2L11}, we deduce 
\begin{equation}
  \label{eq:rho-L1decay-kinetic}
  \forall \, t \ge 0, \quad \| \rho_t \|_{\sfL^{1}(\R_+)} \lesssim  t^{-\frac{1}{2}},
\end{equation}
which matches the decay rate for the heat equation; see~\eqref{eq:rateL1-heat}. We have the same decay estimate on $\iota_t$ and on the first velocity moment of $f_t$. 
In order to prove that claim, we proceed in a similar way as for $\rho_t$. 

We assume $\langle v \rangle^3 \omega^{-2}$ integrable,  and for any $A > 0$, we write 
\begin{align*}
\int_\Omega f |v| \dd x \, \dd v
  &= \int_\Omega f |v| {\bf 1}_{|x| \le A |v|} \dd x \, \dd v  +  \int_\Omega f |v| {\bf 1}_{|x| > A |v|} \dd x \, \dd v 
  \\
    &\le \Bigl( \int_\Omega |v|^2  \omega^{-2} {\bf 1}_{|x| \le A |v|} \dd x \, \dd v \Bigr)^{\frac{1}{2}}  \Bigl( \int_\Omega f^2 \omega^2 \dd x \, \dd v \Bigr)^{\frac{1}{2}}    + A^{-1}  \| \rho \|_{\sfL^{1}_1(\R_+)} 
  \\
    &\le A^{\frac{1}{2}} \left( \int_{\Omega_v} |v|^3  \omega^{-2} \dd v \right)^{\frac{1}{2}}  \left( \int_\Omega f^2 \omega^2 \dd x \, \dd v \right)^{\frac{1}{2}}    + A^{-1}  \| \rho \|_{\sfL^{1}_1(\R_+)}.    
\end{align*}

Optimising on $A$, we  get 
\begin{equation*}
  \| v f  \|_{\sfL^{1}(\Omega)}  \lesssim \| f \|_{\sfL_\omega^{2}(\Omega)} ^{\frac23} \| \rho \|_{\sfL^{1}_1(\R_+)} ^{\frac13}.
\end{equation*}
Together with \eqref{eq:decay-int} and \eqref{eq:estim-rhoL11}, we have established  the claim
\begin{equation}
  \label{eq:1stmoment-L1decay-kinetic}
  \forall \, t \ge 0, \quad  \| v f  \|_{\sfL^{1}(\Omega)} \lesssim  t^{-\frac12}, 
\end{equation}
what in turn implies 
\begin{equation}
  \label{eq:j-L1decay-kinetic}
  \forall \, t \ge 0, \quad \| \iota_t \|_{\sfL^{1}(\R_+)} \lesssim  t^{-\frac12}.
\end{equation}

On the other hand, we introduce the weight function
\begin{equation*}
  \mathsf{w}(x,v) := x^2 + 2 x v + 4 v^2 \gtrsim   x^2 + v^2.
\end{equation*}
Utilising the identities $\LL^\sharp v = -v$ and $\LL^\sharp v^2 = 1 - v^2$, we compute
\begin{align*}
  \frac{\mathrm{d}}{\mathrm{d}t} \int_\Omega f_t \, \mathsf{w} \dd x \dd v
  &= \int_\Omega f_t v \, \partial_x \mathsf{w} \dd x \dd v + \int_{\partial\Omega_-} f_t v^3 \dd x \dd v + \int_\Omega f_t \, \LL^\sharp \mathsf{w} \dd x \dd v \\
  &\leq \int_\Omega f_t \left(2 x v + 2 v^2\right) \dd x \dd v + \int_\Omega f_t \left(2 x \LL^\sharp v + 4 \LL^\sharp v^2\right) \dd x \dd v \\
  &\lesssim \int_\Omega f_t \dd x \dd v % - 2 \int_\Omega f_t v^2 \dd x \dd v
    \lesssim \frac{1}{t^{\frac{1}{2}}},
\end{align*}
where we have used that the boundary term is non-positive, as well as $\LL^\sharp v = -v$, and $\LL^\sharp v^2 = 1-v^2$ for the linear relaxation operator and $\LL^\sharp v^2 = 2-2v^2$ for the Fokker-Planck operator. Integrating in time this differential inequality,  we obtain
\begin{equation}
  \label{eq:rho-L12estimate-kinetic}
  \int_\Omega f_t \left( x^2+v^2 \right) \dd x \dd v \lesssim  \int_\Omega f_t \mathsf w \dd x \dd v \lesssim  t^{\frac12}.
\end{equation}

Note also that the total first-component average 
\begin{equation*}
  \NN(t):= \int_\Omega f_t \left( x_1 + v_1 \right) \dd x \dd v
\end{equation*}
is non-decreasing along time:
\begin{equation*}
  \NN'(t) = \int_{\partial \Omega_-} f_t v_1^2 \dd x' \dd v \ge 0.
\end{equation*}

Together with  the Cauchy-Schwarz inequality, as in the proof of~\eqref{eq:heat-rho-estimate-from-below}, we thus obtain
\begin{align*}
  \forall\, t \ge T_0, \quad \NN_0^2
  & \le \NN(t)^2 \le \left[ \int_\Omega f_t \left( |x_1| + |v_1| \right) \dd x \dd v \right]^2 \\
  & \lesssim \| \rho_t \|_{\sfL^{1}(\R_+)} \left( \int_\Omega f_t \mathsf w \dd x \dd v \right)  \lesssim t^{\frac{1}{2}} \| \rho_t \|_{\sfL^{1}(\R_+)}.
\end{align*}
This yields
\begin{equation}
  \label{eq:rho1}
  \forall\, t \ge T_0, \quad \| \rho_t \|_{\sfL^{1}(\R_+)} \gtrsim t^{-\frac{1}{2}},
\end{equation}
which is precisely the same lower bound as in~\eqref{eq:heat-rho-estimate-from-below} for the heat equation.

This latter estimate can now be localized as follows. Applying the Cauchy–Schwarz inequality together with the $\sfL^2(\R_+)$ bound on $\rho_t$, we obtain
\begin{equation}
  \label{eq:rho2}
  \int_0^{a \sqrt{t}} \rho_t(x)\, \mathrm{d}x \leq \left( at^{\frac{1}{2}} \right)^{\frac{1}{2}} \| \rho_t \|_{\sfL^2(\mathbb{R}_+)} \lesssim a^{\frac{1}{2}} t^{-\frac{1}{2}},
\end{equation}
whereas, using the estimate~\eqref{eq:estim-rhoL11}, we find
\begin{equation}
  \label{eq:rho3}
  \int_{b \sqrt{t}}^\infty \rho_t(x)\, \mathrm{d}x \leq \left(bt^{\frac{1}{2}}\right)^{-1} \int_0^\infty x \rho_t(x)\, \mathrm{d}x \lesssim b^{-1} t^{-\frac{1}{2}}.
\end{equation}
Combining~\eqref{eq:rho-L1decay-kinetic}–\eqref{eq:rho1}–\eqref{eq:rho2}–\eqref{eq:rho3}, and choosing $a > 0$ sufficiently small, $b > 0$ sufficiently large, and $T_0$ large enough, we deduce that
\begin{equation}
  \label{eq:kinetic-rho-estimate-localization}
  \forall\, t \geq T_0, \qquad \frac{c_0}{t^{\frac{1}{2}}} \leq \int_{a \sqrt{t}}^{b \sqrt{t}} \rho_t(x)\, \mathrm{d}x \leq \frac{c_1}{t^{\frac{1}{2}}},
\end{equation}
for some explicit constants $c_0, c_1 > 0$.

This estimate can be localized in the $v$ variable. From \eqref{eq:1stmoment-L1decay-kinetic}, for any $V_0 > 0$, we have
\begin{equation*}
  \int_0^\infty \int_{|v| > V_0} f_t(x,v) \, \mathrm{d}x \, \mathrm{d}v \leq \frac{1}{V_0} \int_0^\infty \int_{\mathbb{R}} |v| f_t(x,v) \, \mathrm{d}x \, \mathrm{d}v \leq \frac{C}{V_0 t^{1/2}}.
\end{equation*}
Combining this with \eqref{eq:rho1}, we deduce that in the case where $\Omega_v = \mathbb{R}$ is unbounded,
\begin{equation}
  \label{eq:loc-loc}
  \frac{c_0'}{t^{1/2}} \leq \int_{a \sqrt{t}}^{b \sqrt{t}} \int_{|v| \leq V_0} f_t(x,v) \, \mathrm{d}x \, \mathrm{d}v \leq \frac{c_1}{t^{1/2}},
\end{equation}
for $V_0 > 0$ sufficiently large and some constant $c_0' > 0$. (In the case where $\Omega_v$ is bounded, this estimate remains formally valid by choosing $V_0$ large enough so that $\Omega_v \subset [-V_0, V_0]$.)

\subsection{Proof of~\eqref{eq:loc-mass} in Theorem~\ref{theo:main2} for the linear relaxation equation}

As a consequence of  \eqref{eq:loc-loc},
 by the pigeonhole principle, there must exist an interval
\begin{equation}
  \label{eq:defI}
  I_t := \left[x_t - \tfrac12,x_t+ \tfrac12 \right] \subset \left[at^{\frac{1}{2}},bt^{\frac{1}{2}}\right]
\end{equation}
such that
\begin{equation}
  \label{eq:relax-rhoIt}
  \int_{I_t} \int_{\VV} f_t(y,v) \dd y \dd v \gtrsim  \frac{1}{t}
\end{equation}
for some
\begin{equation*}
  \VV := [-V_0,V_0] \subset \Omega_v \quad \text{with} \quad V_0>0.
\end{equation*}

We utilize the transport properties to convert bounds on local averages of $\rho$ into pointwise bounds for $\rho$. The solution $f$ to~\eqref{eq:main} satisfies the following Duhamel formula for all $x \in (0, \infty)$ and $v \in \Omega_v$:
\begin{equation}
  \label{eq:characteristic-relax}
  \forall \, s, t \ge 0, \quad f(t+s, x, v) = e^{-s} f(t, x - v s, v) + M(v) \int_0^s e^{-\sigma} \rho_f(t + s - \sigma, x - \sigma v) \, \mathrm{d}\sigma,
\end{equation}
where we set $f(t, x, v) = 0$ and $\rho(t,x)=0$ whenever $x < 0$. Let us denote
\begin{equation*}
  J_t := \left[x_t-\tfrac12 -V_0,x_t+\tfrac12+ V_0 \right].
\end{equation*}
Neglecting the last term in~\eqref{eq:characteristic-relax} and integrating with respect to both $x$ and $v$, we obtain
\begin{align*}
  \int_{J_t} \rho (t+s,y) \dd y
  &= \int_{J_t} \int_{\Omega_v} f(t+s,y,w)  \dd y \dd w \\
  &\ge e^{-s} \int_{J_t} \int_{\Omega_v} f(t,y-sw,w) \dd y \dd w
  \\
  &\ge e^{-s} \int_{I_t} \int_{\VV} f(t,z,w) \dd z \dd w   \gtrsim \frac{1}{t},
\end{align*}
for any $s \in [0,1]$, where the last line follows from applying~\eqref{eq:relax-rhoIt} and the observation that $J_t - sw$ contains the interval $I_t$ defined in~\eqref{eq:defI} for all $s \in [0,1]$ and $w \in [-V_0, V_0]$. By integrating with respect to time, we obtain
\begin{equation*}
  \int_{0}^{\tfrac13} \!\! \int_{J_t} \rho (t+s,y) \dd s  \dd y \gtrsim \frac{1}{t}.
\end{equation*}
By the pigeonhole principle, there exists some $\tilde{x}_t \in J_t$ such that
\begin{equation*}
  \int_{0}^{\frac{1}{3}} \!\! \int_{\tilde J_t} \rho (t+s,y) \dd y  \dd s\gtrsim \frac{1}{t} \qquad \text{where} \qquad 
  \tilde J_t := \left[\tilde x_t- \tfrac{V_0}{6}, \tilde x_t + \tfrac{V_0}{6}\right] = \tilde x_t + \tfrac16 \VV.
\end{equation*}
We also observe that for any $x \in \tilde{J}_t$, we have
\begin{equation*}
  \tilde{J}_t \ \subset \ x + \frac13 \VV.
\end{equation*}
Neglecting the first term in~\eqref{eq:characteristic-relax}, integrating with respect to $v$, and utilizing the preceding estimate, we obtain
\begin{align*}
  \forall \, s \in \left[\frac{2}{3},1\right], \ x \in \tilde J_t, \quad \rho(t+s,x)
  &\ge \int_0^s \! \int_\VV e^{-  \sigma } \rho_f(t+s- \sigma ,x- \sigma  v) M(v) \dd v\dd \sigma \\
  &\gtrsim  \int_{\tfrac13}^{s} \! \int_{x+\tfrac13\VV}  \rho_f(t+s- \sigma ,y) \dd y \dd \sigma \\
  & \gtrsim  \int_0^{s-\tfrac13} \!\!\! \int_{x+\tfrac13\VV}  \rho_f(t+\sigma,y) \dd y \dd \sigma \gtrsim  \frac{1}{t}.
\end{align*}
Finally, discarding the first term in~\eqref{eq:characteristic-relax} once more and utilizing the previous estimate, we obtain
\begin{equation*}
  \forall \, s \in \left[ \frac{5}{6},1 \right], \ x \in \tilde x_t + \frac1{12}\VV, \quad f(t+s,x,v) \ge \left( \int_0^{\tfrac{1}{12}} e^{- \sigma } \rho_f(t+s- \sigma ,x-\sigma v) \dd \sigma \right) M(v) \gtrsim \frac{M(v)}{t}.
\end{equation*}
Therefore, for any $t \geq T_0 + 1$, there exists $\bar{x}_t \sim t^{\frac{1}{2}}$ such that
\begin{equation*}
  \forall\, s \in \left[0, \frac{1}{6}\right],\ \forall\, x \in \bar{x}_t + \frac{\VV}{12}, \quad f(t+s, x, v) \gtrsim \frac{M(v)}{t},
\end{equation*}
which completes the proof (the upper bound immediately follows from~\eqref{eq:decay-pw}).

 \subsection{Proof of~\eqref{eq:loc-mass} in Theorem~\ref{theo:main2} for  the kinetic Fokker-Planck equation}
\label{subsec:loc-mass-kFP}

% Let again $d=1$, let $\omega$ be an admissible weight function, and let the initial data $f_{\mathrm{init}}$ belong to $\sfL^{1}_{1,\omega} \cap \sfL^{\infty}_{-1,\omega}$, with $\mathcal{N}_0 > 0$ where $\mathcal{N}_0$ is defined in~\eqref{eq:defMMM0}. Observe that the specific form of the collision operator in Subsection~\ref{ss:int} was not utilized; therefore, the same argument yields once more
% \begin{equation}
%   \label{eq:kFP-rhoIt}
%   \frac{c_0}{t} \le \int_{I_t} \rho(t,y) \dd y \le   \int_{a \sqrt{t}}^{b \sqrt{t}} \rho_t \dd x \le  \frac{c_1}{t^{\frac{1}{2}}},
% \end{equation}
% for some $0<a<b<\infty$, $0<c_0<c_1<\infty$, $I_t := \left[x_t - 1,x_t+ 1 \right] \subset [at^{\frac{1}{2}},bt^{\frac{1}{2}}]$ and any $t \ge T_0$, which proves~\eqref{eq:loc-mass-int} for the kinetic Fokker-Planck equation. For the same reasons, the optimality of \eqref{eq:decay-int} and \eqref{eq:decay-pw} also hold in this case.

 Proceeding exactly as during the proof of~\eqref{eq:loc-loc} and ~\eqref{eq:relax-rhoIt}, we deduce that 
 % \begin{equation}
 %   \label{eq:kinetic-rho-estimate-localizationBIS}
 %   \frac{c'_0}{t^{\frac{1}{2}}} \le  \int_{a t^{\frac{1}{2}}}^{bt^{\frac{1}{2}}}\int_{-r}^r f_t \, \dd v \dd x  \le  \frac{c'_1}{t^{\frac{1}{2}}}, 
 % \end{equation}
 % and 
 \begin{equation}
   \label{eq:kFP-rhoItBIS}
   \frac{c'_0}{t} \le \int_{I_t} \! \int_{\VV} f_t \dd x \dd v \le   \int_{a \sqrt{t}}^{b \sqrt{t}}  \! \int_{\VV} f_t \dd x \dd v \le  \frac{c'_1}{t^{\frac{1}{2}}},
 \end{equation}
 for some $V_0>0$, $0<a<b<\infty$, $0<c'_0<c'_1<\infty$, $I_t := \left[x_t - 1,x_t+ 1 \right] \subset [at^{\frac{1}{2}},bt^{\frac{1}{2}}]$ and any $t \ge T_0$.

 Harnack inequalities are well established for this class of equations; see, for example,~\cite[Theorem~5.1]{bony_1969} for a seminal result and~\cite[Theorem~5.1]{Lanconelli_1994} for a comprehensive statement. Specifically, for any $r > 0$, there exist constants $C \geq 1$ and $\delta \in (0,1)$ such that
 \begin{equation*}
   \sup_{|x-x_t| < 1,\, |v| < V_0,\, s \in (0,\delta)} f(s+t,x,v) \leq C
   \inf_{|x-x_t| < 1,\, |v| < V_0,\, s \in (1-\delta,1)} f(s+t,x,v).
 \end{equation*}
 In particular, this yields
 \begin{equation*}
   f(t,x,v) \geq \frac{c''_0}{t}
 \end{equation*}
 for all $t \geq T_0 + 1$, $x \in [x_t - 1, x_t + 1]$, and $|v| \leq V_0$, which concludes the proof of~\eqref{eq:loc-mass} for the kinetic Fokker-Planck equation (again the upper bound immediately follows from~\eqref{eq:decay-pw}).

\subsection{Proof of the optimality of ~\eqref{eq:decay-int}  and \eqref{eq:decay-pw}  in dimension $d=1$}
\label{ss:ioptimality}

On the one hand, combining the first estimate in \eqref{eq:kinetic-rho-estimate-localization} with the Cauchy–Schwarz inequality and the first bound in \eqref{eq:borne-rho&iota}, we obtain
\begin{equation}
  \label{eq:kinetic-optimal1}
  \frac{c_0}{t^{\frac{1}{2}}} \leq \int_{a\sqrt{t}}^{b\sqrt{t}} \rho_t(x)\, \mathrm{d}x \leq (b-a)^{\frac{1}{2}} t^{\frac{1}{4}} \| \rho_t \|_{L^2(\mathbb{R}_+)} \lesssim t^{\frac{1}{4}} \| f_t \|_{\sfL^2_{\omega_1}}.
\end{equation}
This demonstrates that the exponent in the decay estimate~\eqref{eq:decay-int} is optimal. Similarly, applying the first estimate in \eqref{eq:kinetic-rho-estimate-localization} together with Hölder's inequality twice yields
\begin{equation}
  \label{eq:kinetic-optimal2}
  \frac{c_0}{t^{\frac{1}{2}}} \leq \int_{a\sqrt{t}}^{b\sqrt{t}} \rho_t(x)\, \mathrm{d}x \leq \left( \frac{b^2 - a^2}{2} \right) t \| \rho_t \|_{L^\infty_{-1}(\mathbb{R}_+)} \lesssim t \| f_t \|_{\sfL^\infty_{-1,\omega_1}}.
\end{equation}
This confirms that the exponent in the decay estimate~\eqref{eq:decay-pw} is also optimal.
 
\hspace{1cm}

\bibliographystyle{acm}
\bibliography{BMM}

\end{document}